\documentclass[preprint,11pt]{amsart}

\usepackage{psfrag}
\usepackage{color}
\usepackage{amsmath}
\usepackage{amsfonts}
\usepackage{graphicx}
\usepackage{amssymb}
\usepackage{rotating}
\usepackage{multirow}
\usepackage{subfigure}
\usepackage{xspace}
\usepackage{amsthm}
\usepackage{tikz}
\usepackage[percent]{overpic}
\newtheorem{problem}{Problem}

\newtheorem{remark}{Remark}

\DeclareMathOperator*{\divergence}{div}

\DeclareMathOperator*{\trace}{tr}
\DeclareMathOperator*{\Grad}{\boldsymbol\nabla}
\DeclareMathOperator*{\Grads}{\boldsymbol\nabla_{\rm sym}}

\DeclareMathOperator{\grads}{\boldsymbol\nabla_s}
\newcommand\RE{\mathbb{R}}
\newcommand{\derivative}[2]{\frac{\partial #1}{\partial #2}}
\newcommand\Huo{H^1_0(\Omega)^d}
\newcommand\Ldo{L^2_0(\Omega)}
\newcommand\Hub{H^1(\B)^d}

\newcommand\Oft{\Omega^f_t}
\newcommand\Ost{\Omega^s_t}
\newcommand\B{\mathcal B}

\renewcommand\u{\mathbf{u}}
\renewcommand\v{\mathbf{v}}

\renewcommand\c{\mathbf{c}}

\newcommand\n{\mathbf{n}}
\newcommand\w{\mathbf{w}}
\newcommand\x{\mathbf{x}}

\newcommand\s{\mathbf{s}}
\renewcommand\S{\mathbf{S}}
\newcommand\V{\mathbf{V}}
\newcommand\X{\mathbf{X}}
\newcommand\Y{\mathbf{z}}

\newcommand\LL{\boldsymbol{\Lambda}}
\newcommand\ssigma{\boldsymbol\sigma}
\newcommand\llambda{\boldsymbol\lambda}
\newcommand\mmu{\boldsymbol\mu}

\newcommand\F{\mathbb{F}}

\renewcommand\P{\mathbb{P}}

\newcommand\ds{\mathrm{d}\s}

\newcommand\dt{\Delta t}
\newcommand\dr{\delta\rho}
\newcommand\T{\mathcal{T}}

\newcommand\qv{\mathbf{q}}
\newcommand\E{E} %Element of a mesh

\newcommand\blue{\begin{color}{black}}
\newcommand\noblue{\end{color}}
\renewcommand\lg{\begin{color}{black}}
\newcommand\gl{\end{color}}

\newcommand\fc{\begin{color}{black}}
	\newcommand\cf{\end{color}}

\newcommand\rv{\begin{color}{black}}
	\newcommand\vr{\end{color}}

\begin{document}

%\begin{frontmatter}

%\title{A parallel solver for the solution of fluid structure interaction problems with a distributed Lagrange multipliers approach}

\title[A parallel solver for FSI problems with Lagrange multiplier]{\lg A parallel solver for fluid structure interaction problems with Lagrange multiplier \gl}

\author{Daniele Boffi}
\address{Computer, electrical and mathematical sciences and engineering division, King Abdullah University of Science and Technology, Thuwal 23955, Saudi Arabia and Dipartimento di Matematica ``F. Casorati'', Universit\`a degli Studi di Pavia, Via Ferrata 5, 27100, Pavia, Italy}
\email{daniele.boffi@kaust.edu.sa}
\urladdr{kaust.edu.sa/en/study/faculty/daniele-boffi}
\author{Fabio Credali}
\address{Istituto di Matematica Applicata e Tecnologie Informatiche ``E. Magenes'', CNR, Via Ferrata 5, 27100, Pavia, Italy}
\email{fabio.credali@imati.cnr.it}
%\urladdr{cemse.kaust.edu.sa/amcs/people/person/fabio-credali}
%
\author{Lucia Gastaldi}
\address{Dipartimento di Ingegneria Civile, Architettura,Territorio, Ambiente e di Matematica, Universit\`a degli Studi di Brescia, via Branze 43, 25123, Brescia, Italy}
\email{lucia.gastaldi@unibs.it}
\urladdr{lucia-gastaldi.unibs.it}
\author{Simone Scacchi}
\address{Dipartimento di Matematica ``F. Enriques'', Universit\`a degli Studi di Milano, Via Saldini 50, 20133 Milano, Italy}
\email{simone.scacchi@unimi.it}
\urladdr{mat.unimi.it/users/scacchi/}

%\subjclass{65N30, 65N12, 74F10, 65F08}

\begin{abstract}
The aim of this work is to present a parallel solver for a formulation of fluid-structure interaction (FSI) problems
which makes use of a distributed Lagrange multiplier in the spirit of the fictitious domain method. The
fluid subproblem, consisting of the non-stationary Stokes equations, is discretized in
space by $\mathcal{Q}_2$--$\mathcal{P}_1$ finite elements, whereas the structure subproblem, consisting of the linear
or finite incompressible elasticity equations, is discretized in space by $\mathcal{Q}_1$ finite elements. A first order
semi-implicit finite difference scheme is employed for time discretization. The resulting
linear system at each time step is solved by a parallel GMRES solver, accelerated by block
diagonal or triangular preconditioners. The parallel implementation is based on the PETSc library. Several numerical
tests have been performed on Linux clusters to
investigate the effectiveness of the proposed FSI solver.\medskip


\end{abstract}

\maketitle

%\end{frontmatter}

\section*{Introduction}
Fluid-structure interaction (FSI) problems are widely studied in literature and several mathematical models and computational methods have been introduced during the past decades, such as the Arbitrary Lagrangian--Eulerian formulations \cite{hirt1974arbitrary,donea1977lagrangian,hughes1981lagrangian,donea2004arbitrary}, the level set formulation \cite{chang1996level} and the Nitsche--XFEM method \cite{alauzet2016nitsche,burman2014unfitted}.

In this work we are going to consider a fictitious domain approach with distributed Lagrange multiplier \cite{BCG15} born as evolution of the finite element immersed boundary method \cite{pesk}, which is still under study from several points of view \blue\cite{ROY20151167,ROWLATT201629,nestola2019immersed,KRISHNAN2020104500}\noblue. In particular, in our formulation, the evolution of the structure is represented via Lagrangian description on a reference domain, which is mapped, at each time step, to the actual position. The fluid is described by an Eulerian model and, in particular, its mesh is extended also in the region occupied by the solid body: the coupling between the two models is enforced through a bilinear form associated with a Lagrange multiplier.
Several numerical and theoretical aspects of this formulation have been recently studied \cite{wolf,boffi2022existence}. %Due to the complexity of this type of problems, primarily in terms of computational cost, we need to balance the accuracy of the method with its efficiency: the design of preconditioners for efficient and robust parallel solver plays a crucial role in this sense.

\fc
Solvers for fluid-structure interaction problems are commonly divided into two families: \textit{partitioned solvers} \cite{BADIA20087027,DEGROOTE2010446,LANGER201547,LI2016272,BANKS2018455} are based on alternate iterations between fluid and solid with interface conditions to tackle the exchange of information between the involved entities. Since several disadvantages (such as high complexity, loss of robustness, cost of the iteration process) have been observed, \textit{monolithic solvers} \cite{heil2008solvers,dunne2010numerical,mayr2015temporal,LANGER2018186} have been studied in recent years. In this case, all fluid and solid variables are computed at the same time by solving a single linear system at each time step: this kind of design usually produces extremely large systems, which are also ill-conditioned so that the use of direct solvers becomes prohibitive. For this reasons, a key task is to balance accuracy and efficiency through the design of preconditioners for robust parallel solvers.

\cf
In the last decade, several studies have focused on the construction and analysis of effective parallel algorithms for FSI, see e.g. \cite{cai2010,crosetto2011,cai2014,balzani2016,deparis2016,heinlein2016,langer2019,MUDDLE20127315}. However, all these works considered ALE formulations of FSI. To the best of our knowledge, the development of \fc monolithic \cf parallel solvers for the fictitious domain approach proposed in \cite{BCG15} is still lacking in the literature. \fc
Indeed, in \cite{BCG15} and the related works \cite{pesk,boffi2011finite,boffi2022interface}, it was observed that the serial implementation of this methodology is computationally very expensive due to the need of reassembling the coupling terms at each time step. Thus, in order to extend the numerical technique introduced in \cite{BCG15} to more complex multiphysics  problems, the design of a parallel solver is mandatory.
\cf

The novelty of the present investigation consists of constructing and analyzing numerically parallel block preconditioners for the linear systems arising from finite element discretizations of the FSI fictitious domain approach developed in \cite{BCG15}. For the fluid, we consider the time dependent Stokes equations, whereas for the structure we consider both linear and nonlinear constitutive laws. \fc The investigation is performed on two academic examples helpful to understand the computational challenges posed by the approach under consideration. \cf
The fluid subproblem is discretized in
space by $\mathcal{Q}_2$--$\mathcal{P}_1$ finite elements, whereas the structure subproblem
is discretized in space by $\mathcal{Q}_1$ finite elements.
%A first order semi-implicit monolithic finite difference scheme is employed for time discretization.
\fc A modified Backward Euler scheme is employed for time discretization, where the coupling terms are treated semi-implicitly. Thanks to the formulation with the Lagrange multiplier for imposing the kinematic condition, the scheme results to be unconditionally stable \cite{BCG15}, while previous approaches based on pointwise update of the structure required a CFL condition to be satisfied in order to have stability \cite{boffi2007numerical}. \cf
The resulting
linear system is solved at each time step by a parallel GMRES solver, accelerated by coupled block
diagonal or triangular preconditioners, where the diagonal blocks are inverted exactly by parallel direct methods. 
The parallel implementation is based on the PETSc library \cite{petsc-web-page,petsc-user-ref,petsc-efficient}. Several numerical
tests have been performed on Linux clusters to
investigate the scalability and robustness of the proposed FSI solver, \fc extending the discussion presented in our preliminary work~\cite{mca28020059}. \cf

The outline of the paper is the following. In Section~\ref{sec:setting}, we describe the mathematical models governing the fluid-structure interaction problems under examination, recalling the derivation of the variational formulation in the setting of the fictitious domain approach with distributed Lagrange multiplier. In Section~\ref{sec:num_methods}, we describe the main features of the numerical methods used for the simulations: finite element spaces, time marching scheme and assembly techniques for the coupling matrix. After this introductory part, we propose, in Section~\ref{sec:para_precon}, two preconditioners for the efficient solution of the linear system arising from the discretization of our problem; several numerical tests are then presented in Section~\ref{sec:numerical_res} in the case of linear and nonlinear solid models: we study robustness with respect to mesh refinement, strong scalability and other properties of interest.

We are going to use classical notation in functional analysis \cite{lions2012non}. Given an open and bounded domain $D$, we denote by $L^2(D)$ the space of square integrable functions, while $L^2_0(D)$ denotes its subspace of zero mean valued functions. Sobolev spaces are denoted by $W^{s,q}(D)$, where $s\in\RE$ represents the differentiability and $q\in [1,\infty]$ is the intergrability exponent. In particular, when $q=2$, we adopt the notation $H^s(D)=W^{s,2}(D)$: norm and scalar product in $H^s(D)$ are denoted by $\|\cdot\|_s$ and $(\cdot,\cdot)_{s,D}$, respectively and, when $s=0$ they refer to $L^2(D)$. Moreover, $H^1_0(D)\subset H^1(D)$ is the space of functions with zero trace on $\partial D$. In the case of vector valued spaces, the dimension is explicitly indicated. %The indication of the space in norms and scalar products is omitted when no confusion arises and, in particular, for the fluid domain $\Omega$.

\section{Mathematical models}
\label{sec:setting}
The object of this study is the numerical approximation of a fluid-structure interaction (FSI) problems, where a visco-elastic solid immersed in a fluid is initially distorted from its equilibrium configuration, so that it tends to return to its equilibrium position. 
Both the fluid and the solid are assumed to be incompressible. 

Let $\Omega\subset\RE^d$, with $d=2,3$, be a connected, open, and bounded domain with Lipschitz continuous boundary $\partial\Omega$; 
for simplicity, let $\Omega$ be a polygon or a polyhedron. 
The domain $\Omega$ is split into two non-intersecting open domains $\Oft$ and $\Ost$, 
representing the regions occupied at time $t$ by fluid and solid, respectively. 
Hence we have $\overline\Omega=\overline{\Oft}\cup\overline{\Ost}$. 
We denote by $\Gamma_t$ the interface between $\Oft$ and $\Ost$, that we assume to have empty intersection with the exterior boundary $\partial\Omega$.
Let $\B$ be the reference domain of $\Ost$, and let $\X:\B\to\Ost$ represent the corresponding deformation mapping. 
Hence a point $\x\in\Ost$ is the image at time $t$ of a point $\s\in\B$, that is $\x=\X(\s,t)$. For simplicity, we assume that $\B=\Omega^s_0$ is the initial position of the solid. We denote by $\F=\Grad_s\X$ the deformation gradient and by $J=\det(\F)$ its Jacobian.
We are going to use the following notation: $\u_f$, $p_f$, $\ssigma_f$, and $\rho_f$ denote, respectively, velocity, pressure, stress tensor, and density of the fluid. In particular, we consider a Newtonian fluid characterized by the usual Navier--Stokes stress tensor
\begin{equation}
	\label{eq:NSstresstensor}
	\ssigma_f=-p_f\mathbb{I}+\nu_f\Grads\u_f,
\end{equation}
where $\nu_f$ is the fluid viscosity and $\Grads\u=(1/2)\left(\Grad\u_f+(\Grad\u_f)^\top\right)$.
We use an Eulerian description so that the material derivative is given by
\[
\dot\u_f=\frac{\partial\u_f}{\partial t}+\u_f\cdot\Grad\u_f.
\]
Conversely, $\u_s$, $p_s$, and $\rho_s$ stand, respectively, for velocity, pressure, and density of the solid body. 
In this case, the Lagrangian framework is preferred, and the spatial description of the material velocity reads
\begin{equation}
	\label{eq:materialvel}
	\u_s(\x,t)=\frac{\partial\X(\s,t)}{\partial t}\Big|_{\x=\X(\s,t)}
\end{equation}
so that $\dot\u_s=\partial^2\X/\partial t^2$.
Moreover, we assume that the solid material is viscous-hyperelastic, so that the Cauchy stress tensor is given by the sum of a viscous part 
\begin{equation}
	\label{eq:viscous}
	\ssigma_s^v=-p_s\mathbb{I}+\nu_s\Grads\u_s,
\end{equation}
where $\nu_s$ is the viscosity, and an elastic part $\ssigma_s^e$, which can be expressed in terms of the Piola--Kirchhoff stress tensor $\P$:
\begin{equation}
	\label{eq:elastic}
	\P(\F(\s,t))=J\,\ssigma_s^e(\x,t)\,\F^{-\top}(\s,t)\quad\text{for }\x=\X(\s,t).
\end{equation}
The Piola--Kirchhoff stress tensor is related to the positive energy density $W(\F)$, which characterizes hyperelastic materials as:
\begin{equation}
	\label{eq:PK}
	%(\P(\F(\s,t))_{\alpha i}=
	%\frac{\partial W}{\partial\F_{\alpha i}}(\F(\s,t))
	%=\left(\frac{\partial W }{\partial \F}(\F(\s,t))\right)_{\alpha i},
	\P(\F(\s,t)) = \derivative{W(\F(\s,t))}{\F}.
\end{equation}
%
%where $i = 1,\ldots,m$ and $\alpha=1,\ldots,d$. 
%The elastic potential energy of
%the body is given by:
%%
%\begin{equation}
%\label{eq:potenergy}
%E\left(\X(t)\right)=\int_\B W(\F(s,t))\ds.
%\end{equation}
%
Assuming that both the fluid and the solid materials are incompressible, we have the following mathematical model for the fluid-structure system: 
\begin{equation}
	\label{eq:model}
	\left\{
	\aligned
	&\rho_f\dot\u_f=\divergence\ssigma_f
	&&\text{in }\Oft\\
	&\divergence\u_f=0 &&\text{in }\Oft\\
	&\rho_s\frac{\partial^2\X}{\partial t^2}
	={\divergence} _s(J\,\ssigma_s^v\,\F^{-\top}+\P(\F))\ 
	&&\text{in }\B\\
	&\divergence\u_s=0&&\text{in }\Ost\\
	&\u_f=\frac{\partial\X}{\partial t}&&\text{on }\Gamma_t\\
	&\ssigma_f\n_f=-(\ssigma_s^v+J^{-1}\,\P\,\F^\top)\n_s
	&&\text{on }\Gamma_t.
	\endaligned
	\right.
\end{equation}
The last two equations in~\eqref{eq:model} represent the transmission conditions at the interface $\Gamma_t$.
The model is completed with initial and boundary conditions:
\begin{equation}
	\label{eq:init+bc}
	\aligned
	&\u_f(0)=\u_{f0} &&\text{in }\Omega_0^f\\
	&\u_s(0)=\u_{s0} &&\text{in }\Omega_0^s\\
	&\X(0)=\X_0 &&\text{in }\B\\
	&\u_f(t)=0 &&\text{on }\partial\Omega.
	\endaligned
\end{equation}
Following~\cite{BG17}, we apply a fictitious domain approach by extending the
first equation in~\eqref{eq:model} to the whole domain $\Omega$ and by using
the following new unknowns:
\begin{equation}
	\label{eq:fictitious}
	\u=\left\{
	\begin{array}{ll}
		\u_f&\text{ in } \Oft\\
		\u_s&\text{ in } \Ost
	\end{array}
	\right.\quad
	p=\left\{
	\begin{array}{ll}
		p_f&\text{ in } \Oft\\
		p_s&\text{ in } \Ost.
	\end{array}
	\right.
\end{equation}
We also set
\begin{equation}
	\label{eq:nu}
	\nu=\left\{\begin{array}{ll}
		\nu_f&\text{ in } \Oft\\
		\nu_s&\text{ in } \Ost.
	\end{array}
	\right.
\end{equation}
\lg In this setting, we need to impose $\u$ to be equal to the material velocity $\u_s$ in $\Ost$, so that \eqref{eq:materialvel} is satisfied, that is
\begin{equation}\label{eq:kinematic}
	\u(\X(\s,t),t) = \derivative{\X}{t}(\s,t)\quad\text{for }\s\in\B.
\end{equation}
This \textit{kinematic constraint} is weakly enforced making use of a Lagrange multiplier. \gl
Let then $\LL$ be a functional space to be defined later on and
$\c:\LL\times\Hub\to\RE$ a continuous bilinear form such that
\begin{equation}
	\label{eq:cprop}
	\c(\mmu,\Y)=0\quad\forall\mmu\in\LL \qquad\text{ implies }\qquad\Y=0.
\end{equation}
Possible definitions for $\LL$ and $\c$ are:
\begin{itemize}
	\item $\c$ is the duality pairing between $\Hub$ and its dual, that is $\LL=(\Hub)'$ and ${\c(\mmu,\Y)=\langle\mmu,\Y\rangle_\B}$ for $\mmu\in\LL$, $\Y\in\Hub$;
	\item $\c$ is the scalar product in $\Hub$, that is $\LL=\Hub$ and $\c(\mmu,\Y)=(\mmu,\Y)_\B+(\grads\mmu,\grads\Y)_\B$ for $\mmu,\Y\in\Hub$.
\end{itemize}
With this notation, our problem can be rewritten in the following weak form (see \cite{BCG15,BG17} for the derivation).
\begin{problem}
	\label{pb:pbvar}
	For given $\u_0 \in\Huo$ and $\X_0\in W^{1, \infty}(\B)$, find $\u(t) \in\Huo$, $p(t) \in\Ldo$, $\X(t) \in\Hub$, and $\llambda(t) \in \LL$ such that for almost all $t \in (0, T)$:
	\begin{equation}
		\label{eq:FSIvarDLM}
		\left\{
		\begin{aligned}
			&\rho_f \left(\derivative{}{t} \u(t),\v\right)_\Omega + b\left(\u(t), \u(t), \v\right) + a\left(\u(t), \v\right)  &&\\
			&\hspace{1.8cm} - \left(\divergence \v, p(t) \right)_\Omega +\c\left(\llambda(t), \v(\X(\cdot, t))\right) = 0 & &\forall \v \in\Huo \\
			&\left( \divergence \u(t), q \right)_\Omega = 0  & &\forall q \in\Ldo
			\\
			&\dr \left(\derivative{^2 \X}{t^2}(t), \Y \right)_\B +\left( \P(\F(t)),
			\nabla_s \Y\right)_\B  -\c\left( \llambda(t), \Y\right) = 0  & &\forall \Y
			\in\Hub \\
			&\c\left(\mmu, \u(\X(\cdot, t), t)- \derivative{\X}{t}(t)\right) = 0 & &\forall \mmu \in \LL \\
			&\u(\x,0) = \u_0(\x)  & &\mathrm{in }\ \Omega\\
			&\X(\s,0) = \X_0(\s) & &\mathrm{in }\ \B.%\\
			%&\lg\derivative{\X}{t}(\s,0) = \X_1(\s)\gl & &\mathrm{in }\ \B.
		\end{aligned}
		\right.
	\end{equation}
\end{problem}
In the above problem we have used the following notation: $\dr=\rho_s-\rho_f$ and
\begin{equation}
	\label{eq:notation}
	\aligned
	& a(\u,\v)=(\nu\Grads\u,\Grads\v)_\Omega\\
	&b(\u,\v,\w)=\frac{\rho_f}2\left((\u\cdot\Grad\v,\w)_\Omega-(\u\cdot\Grad\w,\v)_\Omega\right).
	\endaligned
\end{equation}
\begin{remark}
	We remark that the unknown $\llambda$ in Problem~\ref{pb:pbvar} plays the role of a Lagrange multiplier associated with the condition which enforces the kinematic constraint, that is the equality of the velocity $\u$ with the solid velocity in the region occupied by the structure, see also~\eqref{eq:materialvel}.
\end{remark}

\begin{remark}
	\lg We observe that the initial condition for the time derivative of $\X$ is given by the initial condition for $\u$ thanks to~\eqref{eq:kinematic}.\gl
\end{remark}

\begin{remark}\label{rem:noNS}
	In the following we consider the non-stationary Stokes equations instead of the Navier--Stokes equations and we assume $\rho_s=\rho_f$, 
	thus the terms 
	\[
	b\left(\u(t), \u(t), \v\right)
	\]
	and
	\[
	\dr \left(\derivative{^2 \X}{t^2}(t), \Y \right)_\B
	\]
	in the first and third equation of (\ref{eq:FSIvarDLM}), respectively, are dropped out.
\end{remark}

\rv
\begin{remark}
	In our numerical tests we will consider fluid and solid with same viscosity: as observed by Peskin and McQueen in their seminal work about the Immersed Boundary Method \cite{peskin1989three}, this is a reasonable mathematical assumption for modeling biological systems.
\end{remark}

\begin{remark}
	The choice of $\rho_s=\rho_f$ is related to the study of the added mass effect observed in~\cite{CAUSIN20054506}. In the cited work, it has been proved that a non-implicit time scheme combined with $\dr=0$ is unconditionally unstable, regardless of the discrete parameters. Since the mentioned phenomenon has been observed in the case of the ALE framework, we test our solver for fictitious domain approach for this particular case to be aware of possible instabilities. We anticipate here that no issues occurred in this sense. The reader interested in knowing how our method behaves when different fluid and solid densities are considered is referred to~\cite{boffi2011finite}.
\end{remark}
\vr

As a consequence of Remark \ref{rem:noNS}, the problem that we are going to consider in the rest of the paper is: 
\begin{problem}
	\label{pb:pbvar2}
	For given $\u_0 \in\Huo$ and $\X_0\in W^{1, \infty}(\B)$, find $\u(t) \in\Huo$, $p(t) \in\Ldo$, $\X(t) \in\Hub$, and $\llambda(t) \in \LL$ such that for almost all $t \in (0, T)$:
	
	\begin{equation}
		\label{eq:FSIvarDLM_2}
		\left\{
		\begin{aligned}
			&\rho_f \left(\derivative{}{t} \u(t),\v\right)_\Omega + a\left(\u(t), \v\right)  &&\\
			&\hspace{1cm} - \left(\divergence \v, p(t) \right)_\Omega +\c\left(\llambda(t), \v(\X(\cdot, t))\right) = 0 & &\forall \v \in\Huo \\
			&\left( \divergence \u(t), q \right)_\Omega = 0  & &\forall q \in\Ldo
			\\
			&\left( \P(\F(t)),
			\nabla_s \Y\right)_\B  -\c\left( \llambda(t), \Y\right) = 0  & &\forall \Y
			\in\Hub \\
			&\c\left(\mmu, \u(\X(\cdot, t), t)- \derivative{\X}{t}(t)\right) = 0 & &\forall \mmu \in \LL \\
			&\u(\x,0) = \u_0(\x)  & &\mathrm{in }\ \Omega\\
			&\X(\s,0) = \X_0(\s) & &\mathrm{in }\ \B.
			%&\lg\derivative{\X}{t}(\s,0) = \X_1(\s)\gl & &\mathrm{in }\ \B.
		\end{aligned}
		\right.
	\end{equation}
\end{problem}
In the rest of the paper, we assume $d=2$.

\section{Numerical methods}\label{sec:num_methods}

\subsection{Finite element discretization}
\label{se:FEM}
Let $\T_h^\Omega$ and $\T_h^\B$ be regular meshes in $\Omega$ and $\B$ respectively,
which are independent of each other. The corresponding mesh sizes are
denoted by $h_f$ and $h_s$, respectively. We consider two finite element spaces
$\V_h\subset\Huo$ and $Q_h\subset\Ldo$ such that the pair $(\V_h,Q_h)$ satisfies
the usual inf-sup condition for the Stokes equations. In particular our choice is the popular $\mathcal{Q}_2$--$\mathcal{P}_1$ pair on quadrilateral meshes
\begin{equation*}
	\begin{aligned}
		\V_h &= \{ \v\in\Huo : \v|_E \in \mathcal{Q}_2(E)\ \forall E\in\T_h^\Omega\}\\
		Q_h &= \{ q\in\Ldo : q|_E\in\mathcal{P}_1(E)\ \forall E\in\T_h^\Omega\}.
	\end{aligned}
\end{equation*}
%We assume that $\T_h^\B$ contains only simplices and introduce the space of continuous piecewise affine functions on $\T_h^\B$
\lg Note that the pressure space $Q_h$ contains functions which are discontinuous from one element to the other. \gl We assume that also $\T_h^\B$ is a quadrilateral mesh and we introduce the space of continuous bilinear functions on $\T_h^\B$ as
\begin{equation}
\label{eq:Sh}
\S_h=\{\Y\in\Hub: \Y|_E\in\mathcal{Q}_1(E)\ \forall E\in\T_h^\B\}.
\end{equation}
%
%In order to discretize $\LL$, we set $\LL_h=\S_h$. With this definition we have that when $\LL=(\Hub)'$ and $\c$ is the duality pairing, we can compute easily $\c$ using the scalar product in $L^2(\B)^d$.
We choose $\LL=(\Hub)'$ so that the coupling form $\c$ is the duality pairing in $\Hub$, and set $\LL_h=\S_h\subset\LL$. Since $\LL_h$ is a suitable subspace contained in $L^2(\B)^d$, at discrete level $\c$ can be replaced by the scalar product
\begin{equation}\label{eq:c_l2_inner}
	\c(\mmu_h,\Y_h) = (\mmu_h,\Y_h)_\B \qquad \forall\mmu_h\in\LL_h,\forall\Y_h\in\S_h.
\end{equation}

Therefore, the discrete counterpart of Problem~\ref{pb:pbvar2} reads as follows.
\begin{problem}
\label{pb:pbvarh}
For given $\u_{0h}\in\V_h$ and $\X_{0h}\in\S_h$, find $\u_h(t)\in\V_h$, $p_h(t)\in Q_h$, $\X_h(t)\in\S_h$, $\llambda_h(t)\in\LL_h$ such that for almost all $t \in (0, T)$:

\begin{equation}
	\label{eq:FSIvarDLMh}
\left\{
	\begin{aligned}
		&\rho_f \left(\derivative{}{t}\u_h(t),\v_h\right)_\Omega + a\left(\u_h(t), \v_h\right)  &&\\
		&\quad- \left(\divergence\v_h, p_h(t) \right)_\Omega +\left(\llambda_h(t), \v_h(\X_h(\cdot, t))\right)_\B = 0 & &\forall \v_h\in\V_h \\
		&\left(\divergence\u_h(t), q_h\right)_\Omega = 0  & &\forall q_h\in Q_h\\
		&\left( \P(\F_h(t)), \grads\Y_h\right)_\B  -\left( \llambda_h(t), \Y_h\right)_\B = 0  & &\forall \Y_h\in\S_h \\
		&\left(\mmu_h, \u_h(\X_h(\cdot, t), t)- \derivative{\X_h}{t}(t)\right)_\B = 0 & &\forall \mmu_h\in\LL_h \\
		&\u_h(\x,0) = \u_{0h}(\x)  & &\mathrm{in }\  \Omega\\
		&\X_h(\s,0) = \X_{0h}(\s) & &\mathrm{in }\  \B.
		%&\lg \derivative{\X_h}{t}(\s,0)=\X_{1h}(\s)\gl& &\mathrm{in }\  \B.
	\end{aligned}
\right.
\end{equation}
\end{problem}

\subsection{Semi-implicit time discretization}
We recall the Backward Euler time discretization as presented in \cite{BCG15}. We subdivide the time interval $[0,T]$ into $N$ equal parts with size $\dt=T/N$ and subdivision points $t_n=n\dt$. Moreover, for a certain function $y(t)$, we set $y^n=y(t_n)$ \fc and we approximate the time derivative at $t_{n+1}$ as \cf %and use the following finite difference in order to approximate the time derivatives:
\begin{equation}
	\label{eq:Finitediff}
	\aligned
	&\frac{\partial y(t_{n+1})}{\partial t}\approx\frac{y^{n+1} -y^{n}}{\dt}.\\
	\endaligned
\end{equation}
%
%Notice that both approximations are of first order.

The fully discrete version of Problem~\ref{pb:pbvar2} using a semi-implicit first-order scheme is the following one.
\begin{problem}
	Given $\u_{0h}\in\V_h$ and $\X_{0h}\in\S_h$, for all $n=1,\dots, N$ find $\u_h^n\in\V_h$, $p_h^n\in Q_h$, $\X_h^n\in\S_h$, and $\llambda_h^n\in\LL_h$ fulfilling:
	
	\begin{equation}
		\label{eq:FSIdiscBDF1}
		\left\{
		\begin{aligned}
			&\rho_f \left( \frac{\u_h^{n+1} - \u_h^{n}}{\dt},\v_h\right)_\Omega + a\left(\u_h^{n+1}, \v_h\right)\\
			&\quad - \left( \divergence \v_h, p_h^{n+1}\right)_\Omega+
			\left(\llambda_h^{n+1}, \v_h(\X_h^{n})\right)_\B = 0 & & \forall \v_h \in\V_h \\
			&\left( \divergence \u_h^{n+1}, q_h \right)_\Omega = 0 
			& & \forall q_h \in Q_h\\
			&\left( \P(\F_h^{n+1}), \nabla_s \Y_h\right)_\B
			- \left(\llambda_h^{n+1}, \Y_h\right)_\B = 0& & \forall \Y_h \in\S_h \\
			&\left(\mmu_h, \u_h^{\fc n+1 \cf}(\X_h^{n}) - \frac{\X_h^{n+1} -\X_h^n}{\dt}\right)_\B = 0 & & \forall \mmu_h \in \LL_h \\
			&\u_h^0 = \u_{0h}, \quad	\X_h^0 = \X_{0h}.
		\end{aligned}
		\right.
	\end{equation}
\end{problem}

\subsection{Matrix formulation}

For sake of simplifying the notation, we take for the moment $\P(\F)=\kappa\,\F=\kappa\,\grads\X$.
Problem (\ref{eq:FSIdiscBDF1}) can be presented in the following matrix form:
\begin{equation}
\label{eq:matrix}
\begin{pmatrix}
A_f		& -B^\top	& 0			& L_f(\X_h^n)^\top \\
-B		& 0	& 0 			& 0 \\
0 		& 0 	& K_s   		& -L_s^\top \\
L_f(\X_h^n) 	& 0   	& -\frac{1}{\dt}L_s 	& 0 \\
\end{pmatrix}	
\begin{pmatrix}
\u_h^{n+1} \\
p_h^{n+1} \\
\X_h^{n+1} \\
\lambda_h^{n+1} \\
\end{pmatrix}
= 
\begin{pmatrix}
g_1 \\
0 \\
0 \\
g_2 
\end{pmatrix},
\end{equation}
with
\begin{equation*}
        \begin{aligned}
        	& A_f = \frac{\rho_f}{\dt} M_f + K_f\\
        	& (M_f)_{ij} = \left( \boldsymbol{\phi}_j,  \boldsymbol{\phi}_i \right)_\Omega,
        	\quad(K_f)_{ij} = a\left(\boldsymbol{\phi}_j, \boldsymbol{\phi}_i\right) \\
        	& B_{ki} = \left( \divergence \boldsymbol{\phi}_i, \psi_k\right)_\Omega\\
        	& (K_s)_{ij} = \kappa \left( \nabla_s \boldsymbol{\chi}_j, \nabla_s \boldsymbol{\chi}_i \right)_\B\\
        	& (L_f(\X_h^n))_{lj} = \left(\boldsymbol{\zeta}_l, \boldsymbol{\phi}_j(\X_h^n)\right)_\B, 
        	\quad (L_s)_{lj} = \left(\boldsymbol{\zeta}_l, \boldsymbol{\chi}_j\right)_\B\\
        	& g_1 = \frac{\rho_f}{\dt}M_f\u_h^n,
        	\quad g_2 = -\frac{1}{\dt} L_s \X_h^n.
        \end{aligned}
\end{equation*}
Here $\boldsymbol{\phi}_i$, $\psi_k$, $\boldsymbol{\chi}_i$, and $\boldsymbol{\zeta}_l$ denote the basis functions in
$\V_h$, $Q_h$, $\S_h$, and $\LL_h$, respectively. \lg We observe that, due to our choice of spaces, $\boldsymbol{\chi}_i=\boldsymbol{\zeta}_i$. \gl

The system above is associated to a steady saddle point problem which admits a unique solution. Thanks to the theory developed in~\cite{BG17}, the finite element discretization is stable, thus giving optimal convergence rates depending on the regularity of the solution.

In the general case, when $\P(\F)$ is nonlinear, \lg the matrix $K_s$ depends on $\X$; hence, we solve the system by means of a solver for nonlinear systems of equations \gl
like a fixed point iteration or Newton like methods.
In our numerical experiments, we will use the Newton method.

The matrix $L_f(\X_h^n)$ plays a crucial role in the discretization of this type of problems because it represents the coupling between fluid and structure via the Lagrange multiplier. Possible assembly techniques for this matrix have been already discussed in~\cite{boffi2022interface} in the particular case of triangular meshes, but they can be extended also to different geometries. We recall the main features of this operation in the particular setting of our numerical methods, where quadrilateral meshes are used for both fluid and structure.

%Since we chose $\c$ to be the scalar product in $\LdB^d$, its discrete counterpart is given by
\lg The discrete counterpart of the coupling term is given by\gl
\begin{equation}
	\lg(\mmu_h,\v_h(\X_h^n))_\B = \int_{\B}\gl \mmu_h\cdot\v_h(\X_h^n)\,d\s = \sum_{\E\in\T_h^\B} \int_{\E} \mmu_h\cdot\v_h(\X_h^n)\,d\s
\end{equation}
where the variable $\mmu_h\in\LL_h$ is defined on the solid mesh, while the velocity variable $\v_h\in\V_h$ is defined on the fluid one. In addition, $\v_h$ is combined with the solid mapping at the previous time step $\X_h^n$, so that the position of the mapped solid mesh $\X_h^n(\T_h^\B)$ with respect to $\T_h^\Omega$ has to be taken into account. In particular, the action of the map $\X_h^n$ is carried out by mapping the nodes of $\T_h^\B$ and keeping straight the edges between them.

For these reasons, the assembly of the matrix under consideration is not a trivial operation, hence we discuss two possible techniques: we can compute the integrals in a coarse way working directly on each $\E\in\T_h^\B$ or, alternatively, we can compute the intersection between $\X_h^n(\T_h^\B)$ and $\T_h^\Omega$ in order to implement an exact composite quadrature rule.% based on a particular refinement of $\T_h^\B$.

In the first case, we consider the first order quadrature rule on squares, where the nodes to be considered are simply the four vertices $\qv_1,\dots,\qv_4$ and the weights are all equal to $1/4$; we have
\begin{equation}\label{eq:inexact_coupling}
	\int_{\E} \mmu_h\cdot\v_h(\X_h^n)\,d\s \approx \frac{area(\E)}{4}\sum_{k=1}^{4} \mmu_h(\qv_k)\cdot\v_h(\X_h^n(\qv_k))\quad\forall \E\in \T_h^\B.
\end{equation}
In particular, we notice that, in order to determine which fluid basis functions are involved in a single integration of this type, we have to detect in which fluid elements the mapped quadrature points $\X_h^n(\qv_k)$ are placed.

The second approach consists in computing the intersection between the mapped solid mesh and fluid one; in this way, each element $\E\in\T_h^\B$ is partitioned into a certain number of polygons $P_1,\dots,P_J$, each contained, once mapped, in a single fluid element. If a polygon is not already a triangle, it is triangulated by connecting the barycenter with each vertex. In this case, we have
\begin{equation}
	\begin{aligned}
		\int_{\E} \mmu_h\cdot\v_h(\X_h^n)\,d\s &= \sum_{j=1}^{J} \int_{P_j} \mmu_h\cdot\v_h(\X_h^n)\,d\s
		=  \sum_{j=1}^{J} \sum_{i=1}^{N_j} \bigg[\int_{T_i} \mmu_h\cdot\v_h(\X_h^n)\,d\s\bigg]\\ &= 
		\sum_{j=1}^{J} \sum_{i=1}^{N_j} area(T_i) \bigg[ \sum_{k=1}^{4} \omega_k\mmu_h(\qv_k)\cdot\v_h(\X_h^n(\qv_k)) \bigg]\quad\forall \E\in \T_h^\B.
	\end{aligned}
\end{equation}
where the chosen rule is characterized by the nodes $\qv_1=(3/5,1/5,1/5)$, $\qv_2=(1/5,3/5,1/5)$, $\qv_3=(1/5,1/5,3/5)$ in barycentric coordinates plus the barycenter $\qv_4$ with weights $\omega_1=\omega_2=\omega_3=25/48$ and $\omega_4=-9/16$.

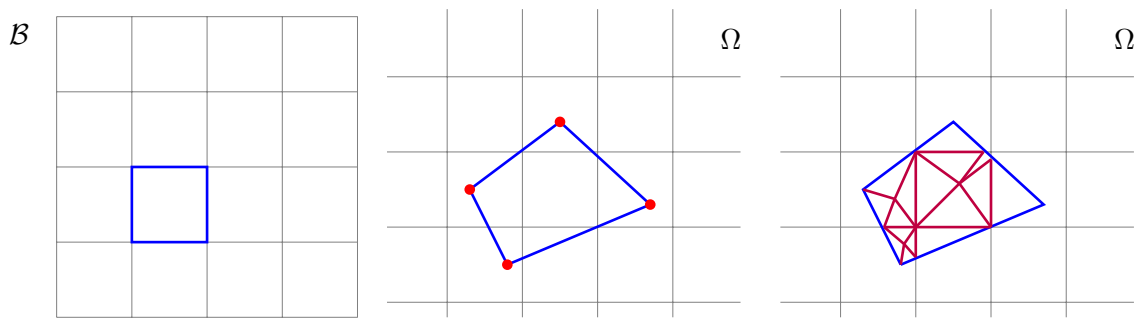
\begin{figure}[h!]
	\centering
	\begin{tikzpicture}
		\node[left=0.5cm,above=4.5cm]{$\B$};\draw[step=1cm,gray,very thin] (0,1) grid (4,5);
		\draw[blue,thick,line width=1pt] (1,2) rectangle (2,3);
		%\node[left=4cm,above=2cm]{$E$};
	\end{tikzpicture}\quad
	%\begin{tikzpicture}
	%	\draw[->] (2,7) .. controls(3,10) .. (5,7);
	%\end{tikzpicture}\quad
	\begin{tikzpicture}
		\node[above=4.5cm,right=4.5cm]{$\Omega$};\draw[step=1cm,gray,very thin] (0.2,0.8) grid (4.9,4.9);
		\draw[blue,thick,line width=1pt] (1.8,1.5) -- (3.7,2.3) -- (2.5,3.4) -- (1.3,2.5) -- (1.8,1.5);
		\draw[fill,red,thick,line width=1pt](1.8,1.5) circle (1.5pt);
		\draw[fill,red,thick,line width=1pt](3.7,2.3) circle (1.5pt);
		\draw[fill,red,thick,line width=1pt](2.5,3.4) circle (1.5pt);
		\draw[fill,red,thick,line width=1pt](1.3,2.5) circle (1.5pt);
	\end{tikzpicture}\quad
	\begin{tikzpicture}
		\node[above=4.5cm,right=4.5cm]{$\Omega$};\draw[step=1cm,gray,very thin] (0.2,0.8) grid (4.9,4.9);
		\draw[blue,thick,line width=1pt] (1.8,1.5) -- (3.7,2.3) -- (2.5,3.4) -- (1.3,2.5) -- (1.8,1.5);
		\draw[purple,thick,line width=1pt] (2,1.6) -- (2,3);
		\draw[purple,thick,line width=1pt] (3,2) -- (1.58,2);
		\draw[purple,thick,line width=1pt] (3,2) -- (3,2.9);
		\draw[purple,thick,line width=1pt] (2.9,3) -- (2,3);
		\draw[purple,thick,line width=1pt] (1.845,1.775) -- (2,1.6);
		\draw[purple,thick,line width=1pt] (1.845,1.775) -- (2,2);
		\draw[purple,thick,line width=1pt] (1.845,1.775) -- (1.58,2);
		\draw[purple,thick,line width=1pt] (1.845,1.775) -- (1.8,1.5);
		\draw[purple,thick,line width=1pt] (1.72,2.375) -- (1.58,2);
		\draw[purple,thick,line width=1pt] (1.72,2.375) -- (2,3);
		\draw[purple,thick,line width=1pt] (1.72,2.375) -- (2,2);
		\draw[purple,thick,line width=1pt] (1.72,2.375) -- (1.3,2.5);
		\draw[purple,thick,line width=1pt] (2.58,2.58) -- (2,2);
		\draw[purple,thick,line width=1pt] (2.58,2.58) -- (2,3);
		\draw[purple,thick,line width=1pt] (2.58,2.58) -- (3,2);
		\draw[purple,thick,line width=1pt] (2.58,2.58) -- (3,2.9);
		\draw[purple,thick,line width=1pt] (2.58,2.58) -- (2.9,3);
	\end{tikzpicture}
	\caption{A schematic representation of the geometric aspects of the coupling operations. From the left hand side: a portion of $\T_h^\B$ with a particular element under consideration, its immersed counterpart in the case of the coarse integration with the rule of the vertices, the same immersed element partitioned accordingly with the computation of the intersection.}
	\label{fig:coupling_geo}
\end{figure}

A schematic representation of these procedures is presented in Figure \ref{fig:coupling_geo}, where a solid element is mapped into the fluid mesh and managed with the two approaches just described. \fc We finally remark that the increase of precision for the inexact approach described in Equation \eqref{eq:inexact_coupling} could still produce a non-optimal method: this fact is discussed in \cite{boffi2022interface}. \cf

\begin{figure}[h!]
	\centering
	\includegraphics[width=6.5cm]{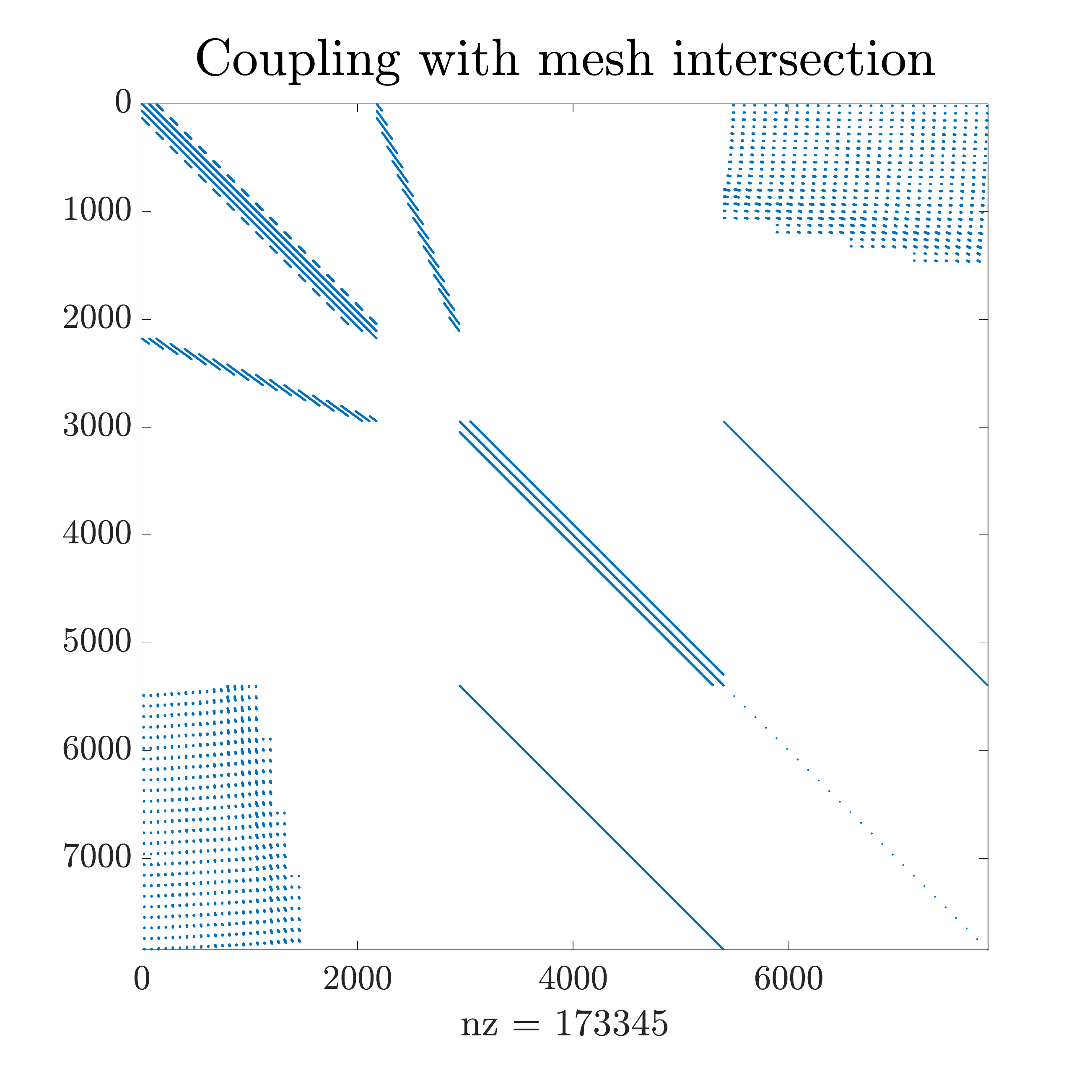}\quad%{figures/spy_fact4}
	\includegraphics[width=6.5cm]{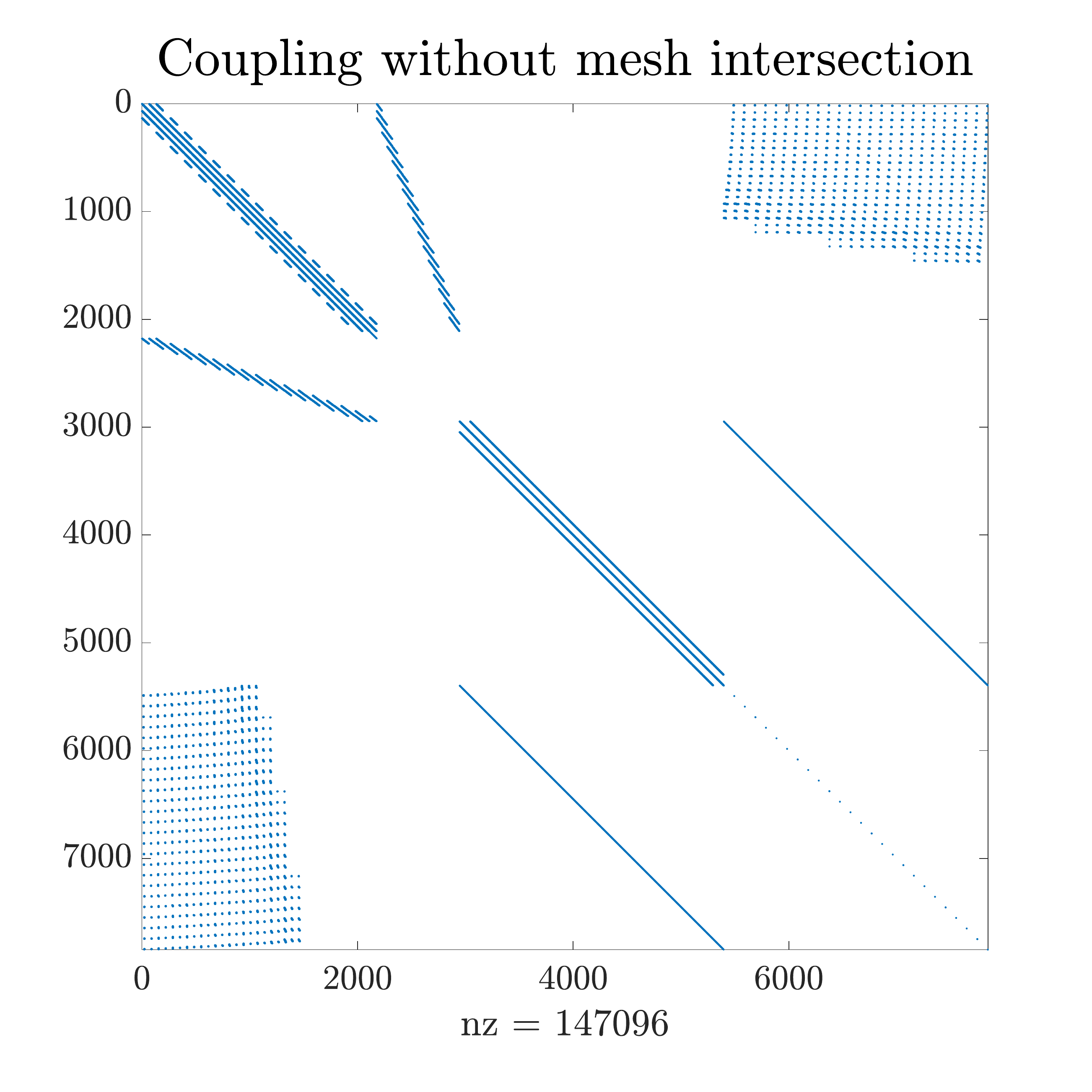}
	\caption{Sparsity patterns for the matrix \eqref{eq:matrix} when $L_f(\X_h^n)$ is computed with the two assembly processes in the case of the linear elastic model described in Section \ref{sec:linear_model} discretized with 7846 global dofs. \lg $nz$ represents the number of nonzero entries of the matrix. The matrix is more dense when the coupling term is assembled with mesh intersection. \gl}
	\label{fig:sparsity}
\end{figure}

\begin{remark}
	%The choice of the assembling technique affects also the structure of the sparsity pattern of $L_f(\X_h^n)$. Indeed, when the assembling is performed after the computation of the mesh intersection, the matrix is more dense with respect to the case without mesh intersection because an higher number of degrees of freedom are involved in the integration. An example is reported in Figure \ref{fig:sparsity}.
	The sparsity pattern of the coupling matrix $L_f(\X_h^n)$ is affected by the choice of the assembly technique. Indeed, if we assembly by computing the intersection between the two involved meshes, we have that the matrix is more dense, since a higher number of degrees of freedom contributes to the interaction (see Figure~\ref{fig:sparsity}).
\end{remark}

\section{Parallel preconditioners}
\label{sec:para_precon}

Our strategy for building an efficient parallel solver is based on the library PETSc from Argonne National Laboratory~\cite{petsc-web-page,petsc-user-ref,petsc-efficient}.
Such library is built on the MPI standard and offers advanced data structures and routines for the parallel solution of partial differential equations, from basic vector and matrix operations to more complex linear and nonlinear equation solvers.
In our Fortran90 code, vectors and matrices are built and subassembled in parallel on each processor.

Let us denote by $\mathcal{A}$ the matrix of the linear system (\ref{eq:matrix}) or the Jacobian matrix in the nonlinear case.
$\mathcal{A}$ has clearly the following block structure
\[
\mathcal{A}=
\begin{pmatrix}
\mathcal{A}_{11} & \mathcal{A}_{12}\\
\mathcal{A}_{21} & \mathcal{A}_{22}
\end{pmatrix}
\]
where
\[
\begin{array}{cc}
\displaystyle \mathcal{A}_{11} = 
\begin{pmatrix}
A_f             & -B^\top  \\
-B              & 0     
\end{pmatrix}
& \displaystyle \mathcal{A}_{12} =
\begin{pmatrix}
0                     & L_f(\X_h^n)^\top \\
0                     & 0 
\end{pmatrix}
\vspace{0.2cm}\\
\displaystyle \mathcal{A}_{21} =
\begin{pmatrix}
0               & 0     \\
L_f(\X_h^n)     & 0     \\
\end{pmatrix}
& \displaystyle \mathcal{A}_{22} =
\begin{pmatrix}
K_s                   & -L_s^\top \\
-\frac{1}{\dt}L_s     & 0 \\
\end{pmatrix}
\\
\end{array}
\]
To solve such  linear system, we use the parallel GMRES method provided by the PETSc library,
preconditioned by two types of block preconditioners:
\begin{itemize}
\item {\bf block-diag} where the preconditioner is the diagonal matrix 
\[
\begin{pmatrix}
\mathcal{A}_{11}	& 0 \\
0     			& \mathcal{A}_{22} \\
\end{pmatrix}
\]
\item {\bf block-tri} where the preconditioner is the triangular matrix
\[
\begin{pmatrix}
\mathcal{A}_{11}        & 0 \\
\mathcal{A}_{21}        & \mathcal{A}_{22} \\
\end{pmatrix}.
\]
\end{itemize}
In both cases, the action of the preconditioner consists of the exact inversion of the two diagonal blocks 
$\mathcal{A}_{11}$ and $\mathcal{A}_{22}$, that we perform by means of the parallel multifrontal direct solver Mumps~\cite{amestoy.2001,amestoy.2006}.

\section{Numerical results}\label{sec:numerical_res}

\lg In this section, we report results obtained from the simulation of fluid-structure interaction systems involving both linear and nonlinear models for the solid evolution. In particular, our code has been run on two Linux clusters. The Shaheen cluster at King Abdullah University of Science and Technology (KAUST, Saudi Arabia) is a Cray XC40 cluster constituted by 6,174 dual sockets compute nodes, based on 16 core Intel Haswell processors running at 2.3GHz. Each node has 128GB of DDR4 memory running at 2300MHz. The Eos cluster at University of Pavia (Italy) is a Linux Infiniband cluster with 21 nodes, each carrying two 16 cores Intel Xeon Gold 6130 processors running at 2.1 GHz. \gl

\subsection{Linear solid model}\label{sec:linear_model}
Let us consider the case of a linear model governing the structure; this means that the Piola--Kirchhoff stress tensor can be written as
\begin{equation}
	\P(\F) = \kappa\F = \kappa\grads\X.
\end{equation}
In this case, the energy density is given by
\begin{equation}
	W(\F) = \frac{\kappa}{2}\F:\F
\end{equation}
and the elastic potential energy takes the form
\begin{equation}
	E(\X) = \frac{\kappa}{2} \int_\B W(\F)\,\ds = \frac{\kappa}{2} \int_\B |\grads\X|^2 \,\ds.
\end{equation}

In particular, we consider the annulus \lg which, at rest, occupies the region \gl ${\{\x\in\RE^2:0.3\le|\x|\le 0.5\}}$ placed at the center of \lg the square $[-1,1]^2$ \gl filled with fluid; the annulus is stretched and its internal forces make it return to its resting configuration. \lg Thanks to the symmetry of the geometry, we can reduce the system to a quarter of the domain. Therefore, we set as reference domain \gl ${\B = \{\s\in\RE^2:s_1,s_2\ge0,\,0.3\le|\s|\le 0.5\}}$.

We impose null velocity $\u(\x,t)=0$ on the upper and right edges of the fluid domain, while on the other part of the boundary we allow fluid and structure to move along the tangential direction, setting to zero the normal component. Furthermore, we impose also the following initial conditions
\begin{equation}
	\u(\x,0) = 0, \qquad \X(\s,0) = \bigg(\frac{s_1}{1.4},1.4\, s_2\bigg)\quad \lg\text{for }\s=(s_1,s_2)\gl.
\end{equation}

We assume that densities and viscosities of solid and fluid are equal, therefore we set
\begin{equation}
	\rho_f=\rho_s=1 \quad\text{and}\quad \nu_f = \nu_s = 0.1
\end{equation}
and, moreover, $\kappa = 10$.

\begin{figure}
	\begin{center}
		\includegraphics[trim = 30 10 20 0,clip, width=4cm]{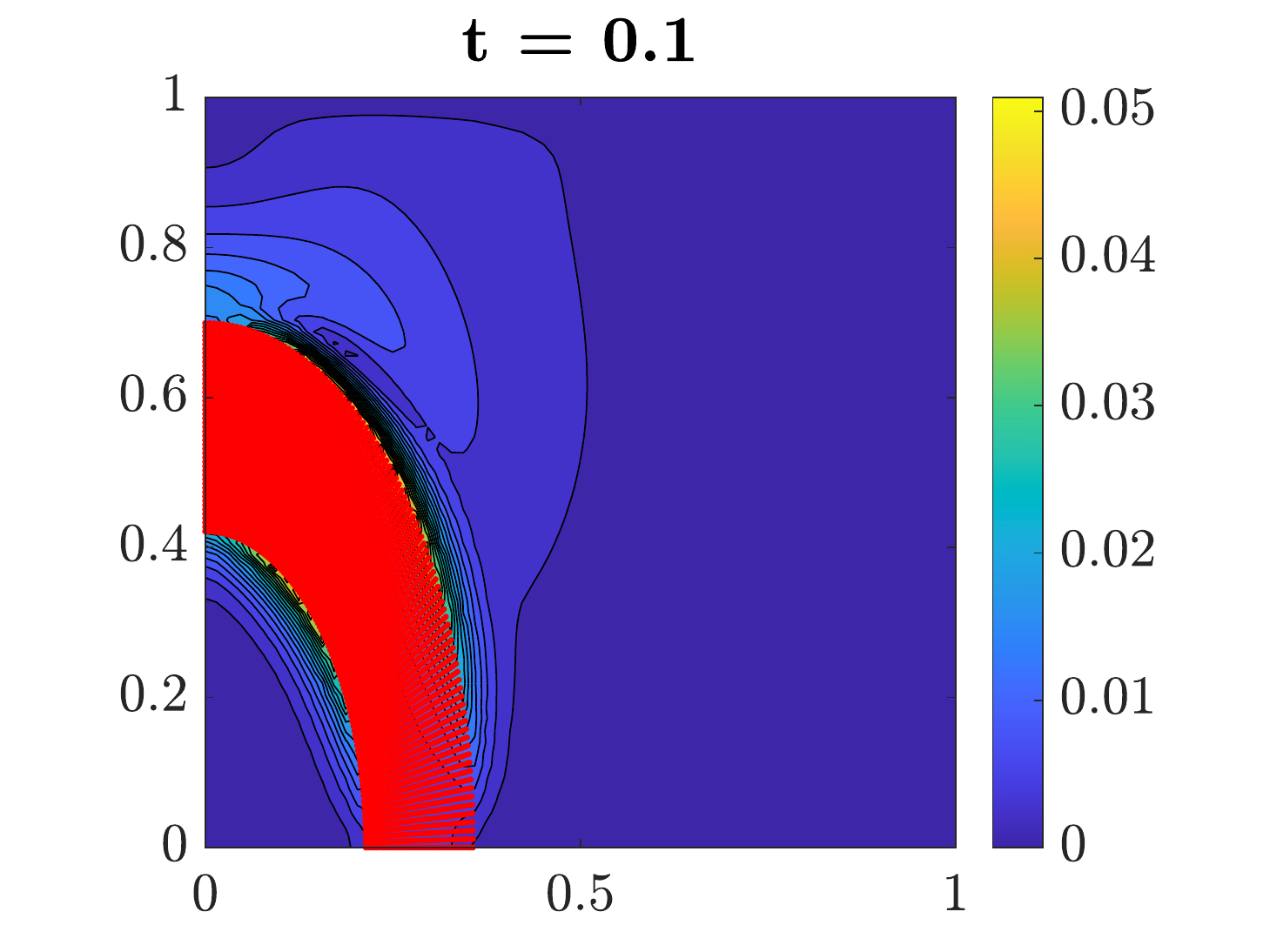}
		\includegraphics[trim = 30 10 20 0,clip, width=4cm]{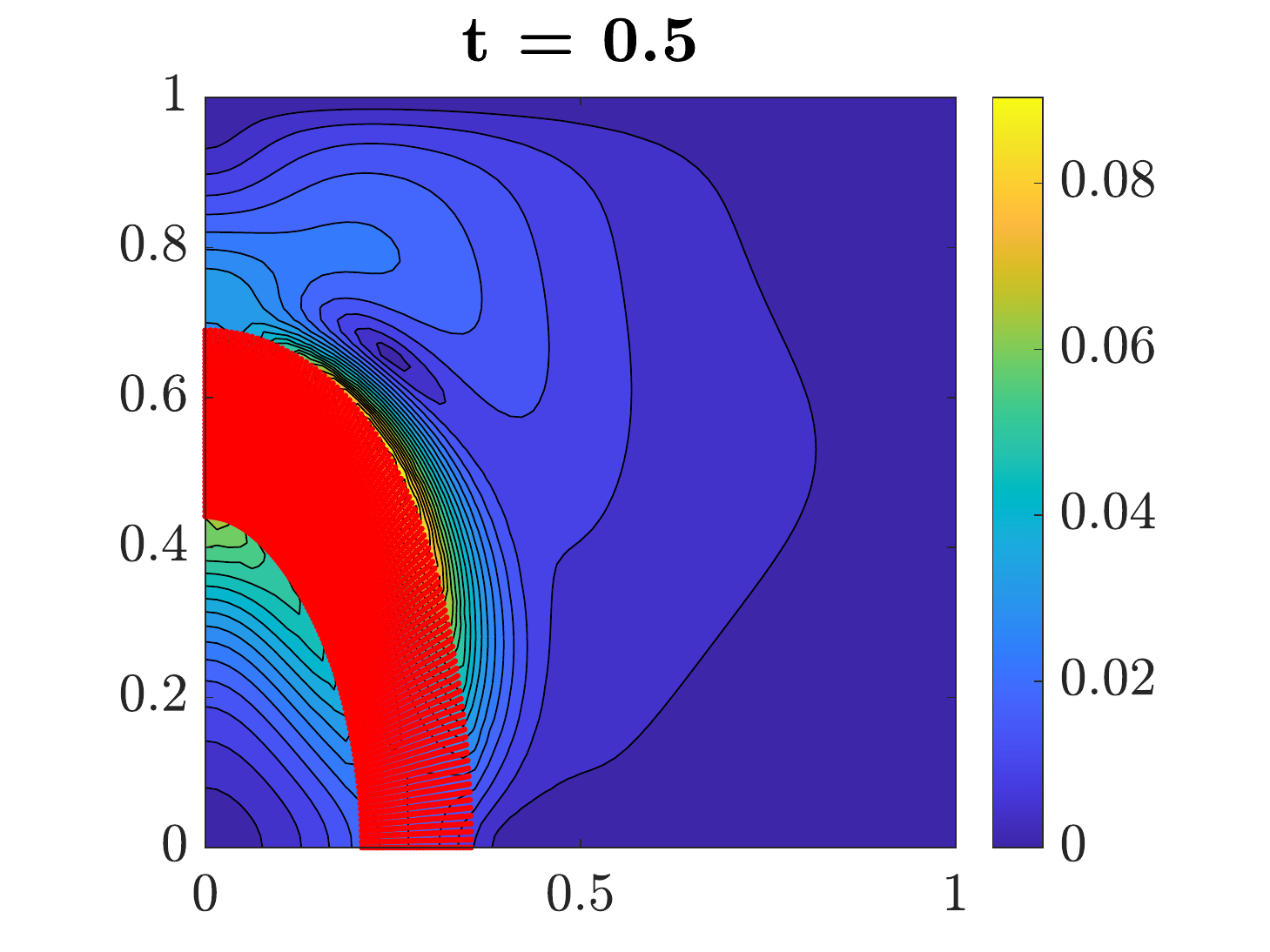}
		\includegraphics[trim = 30 10 20 0,clip, width=4cm]{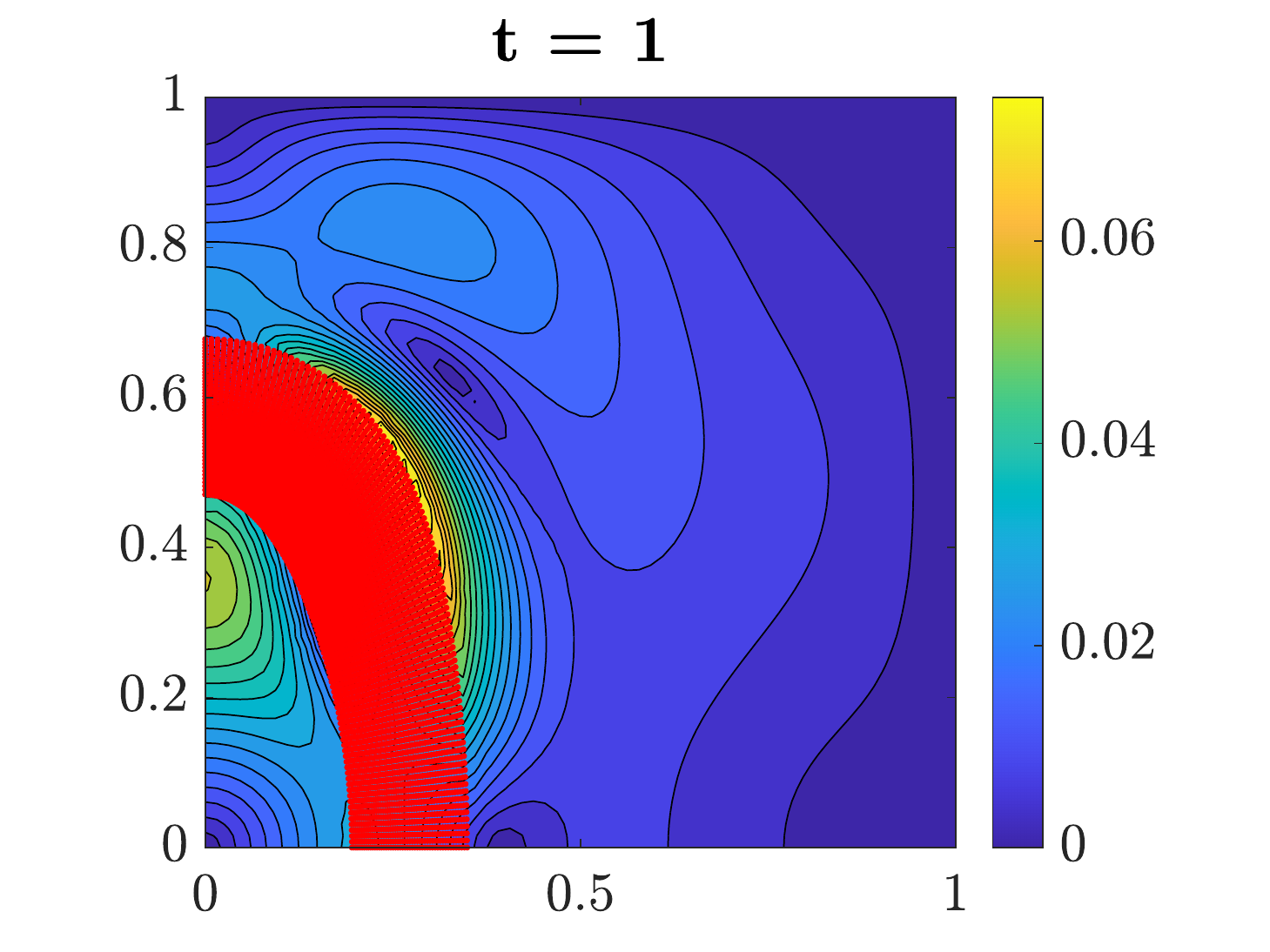}
		\includegraphics[trim = 30 10 20 0,clip, width=4cm]{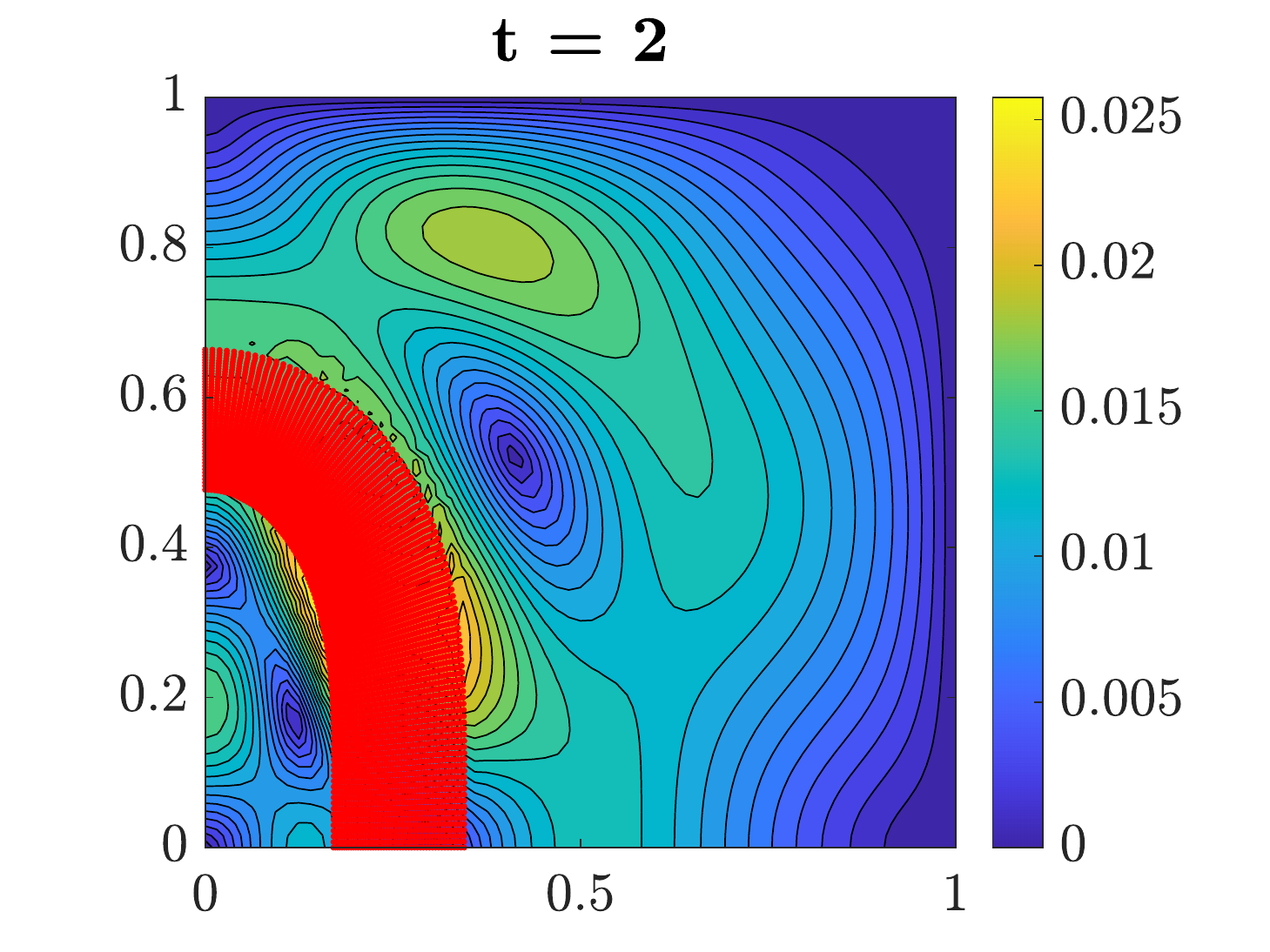}\\
		\medskip
		\includegraphics[trim = 30 10 20 0,clip, width=4cm]{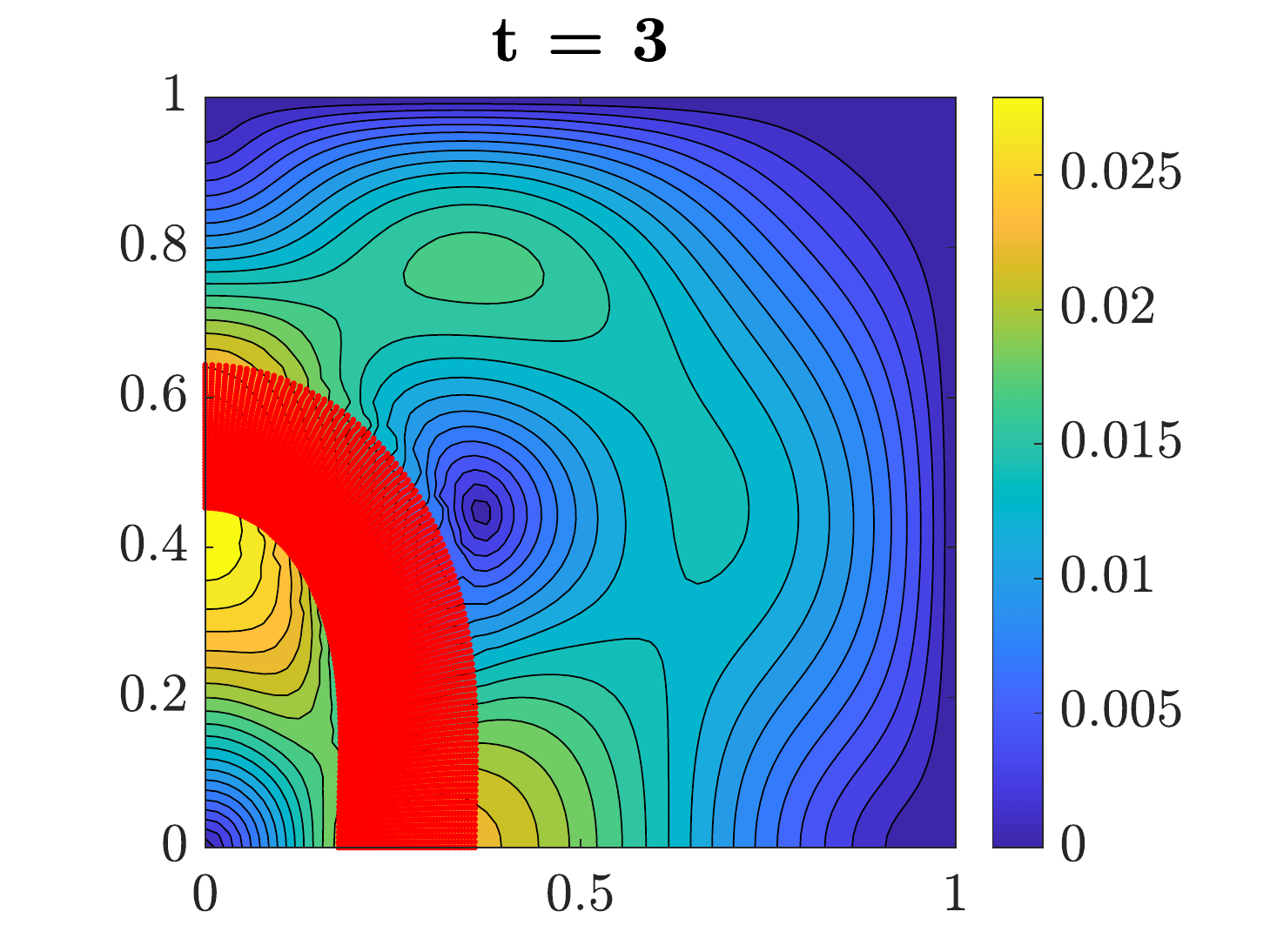}
		\includegraphics[trim = 30 10 20 0,clip, width=4cm]{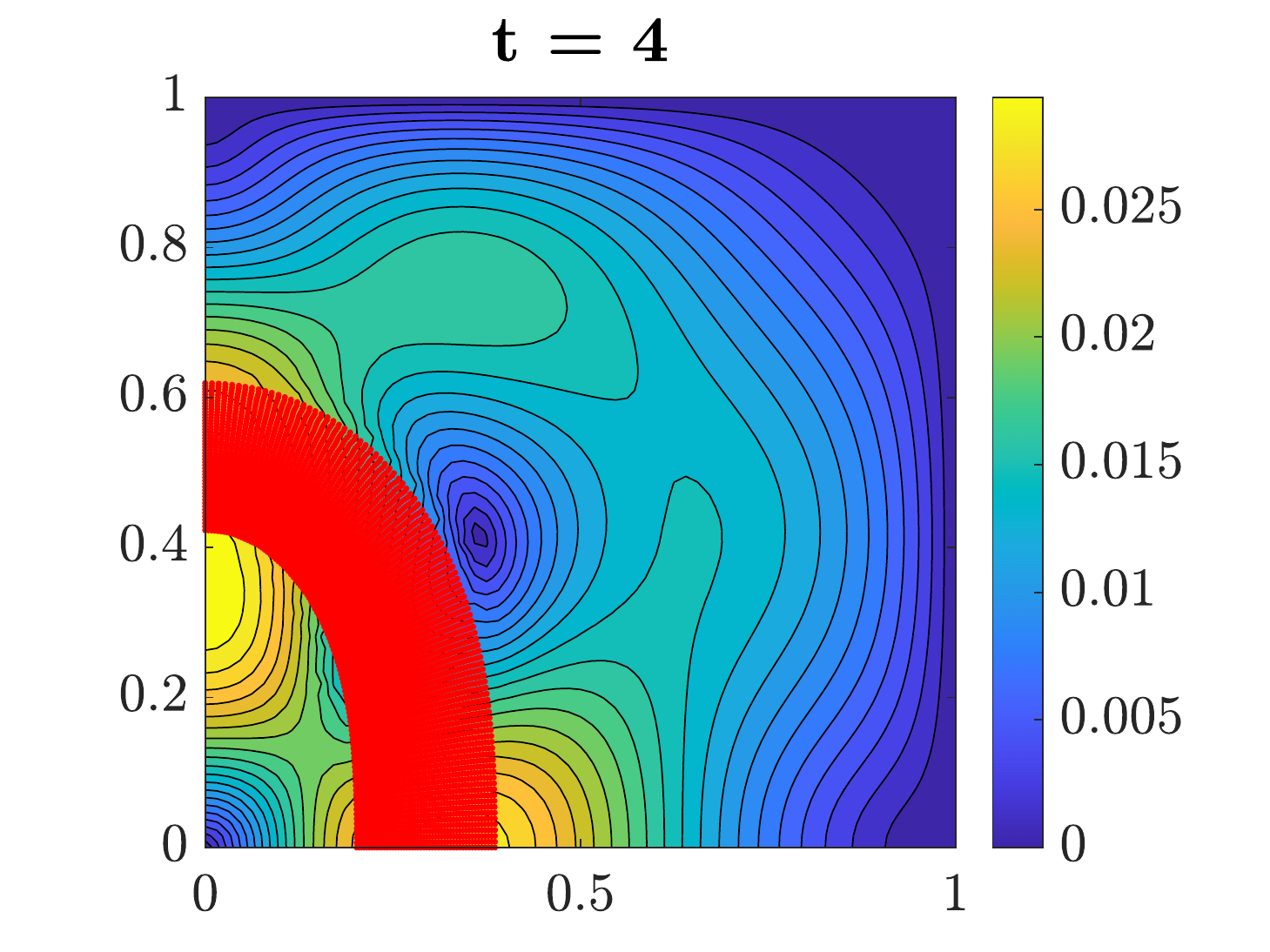}
		\includegraphics[trim = 30 10 20 0,clip, width=4cm]{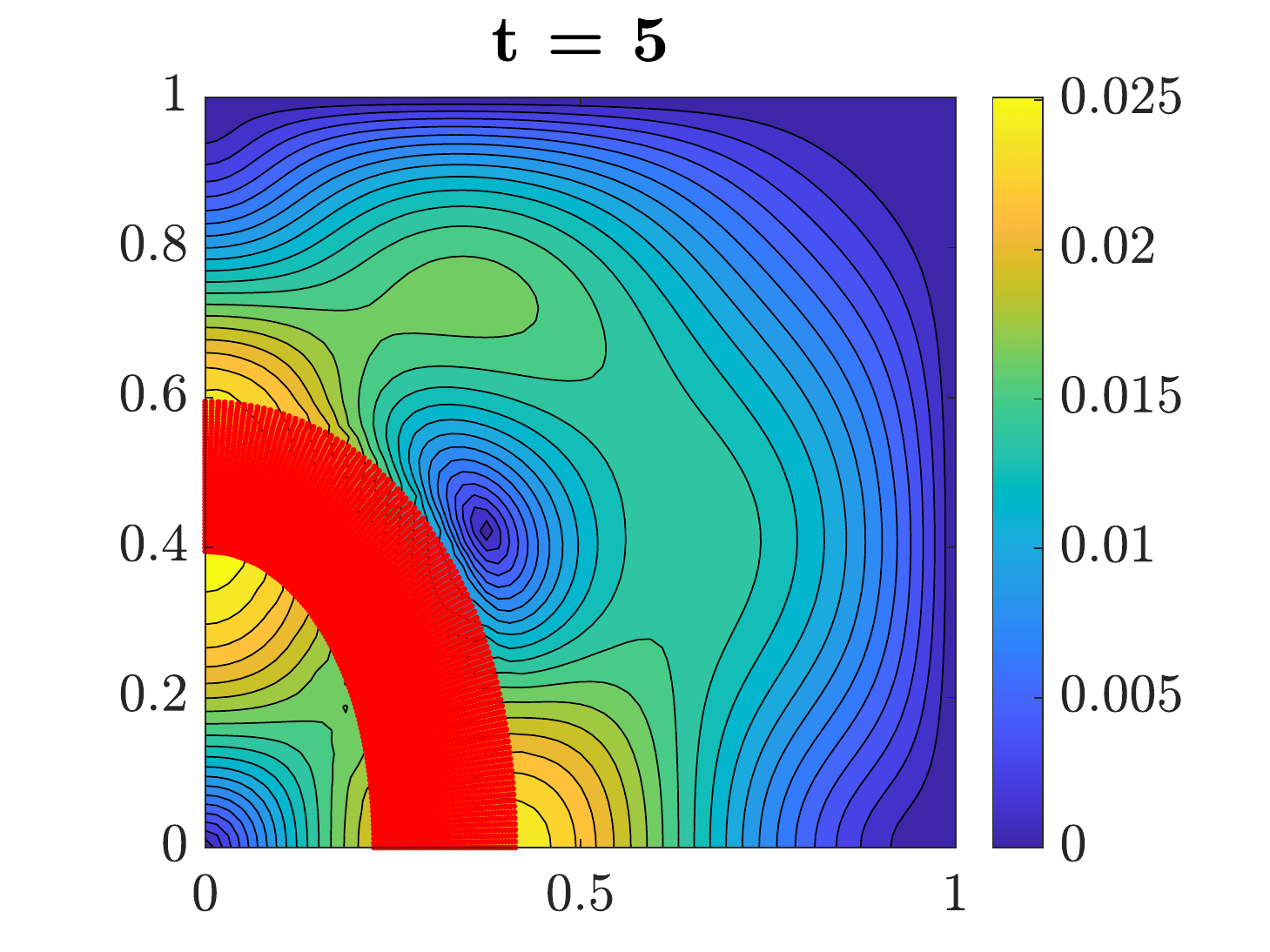}
		\includegraphics[trim = 30 10 20 0,clip, width=4cm]{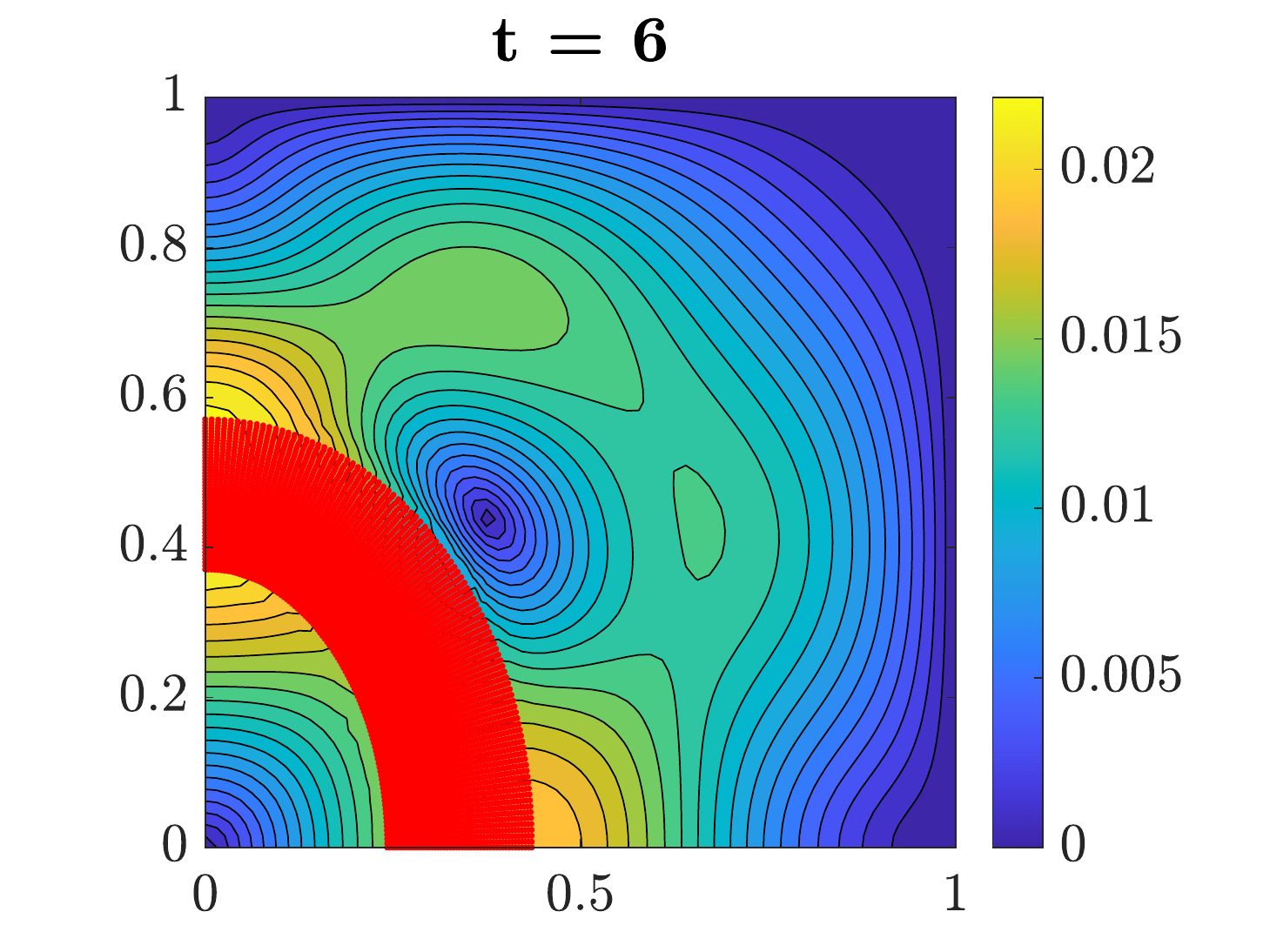}\\
		\medskip
		\includegraphics[trim = 30 10 20 0,clip, width=4cm]{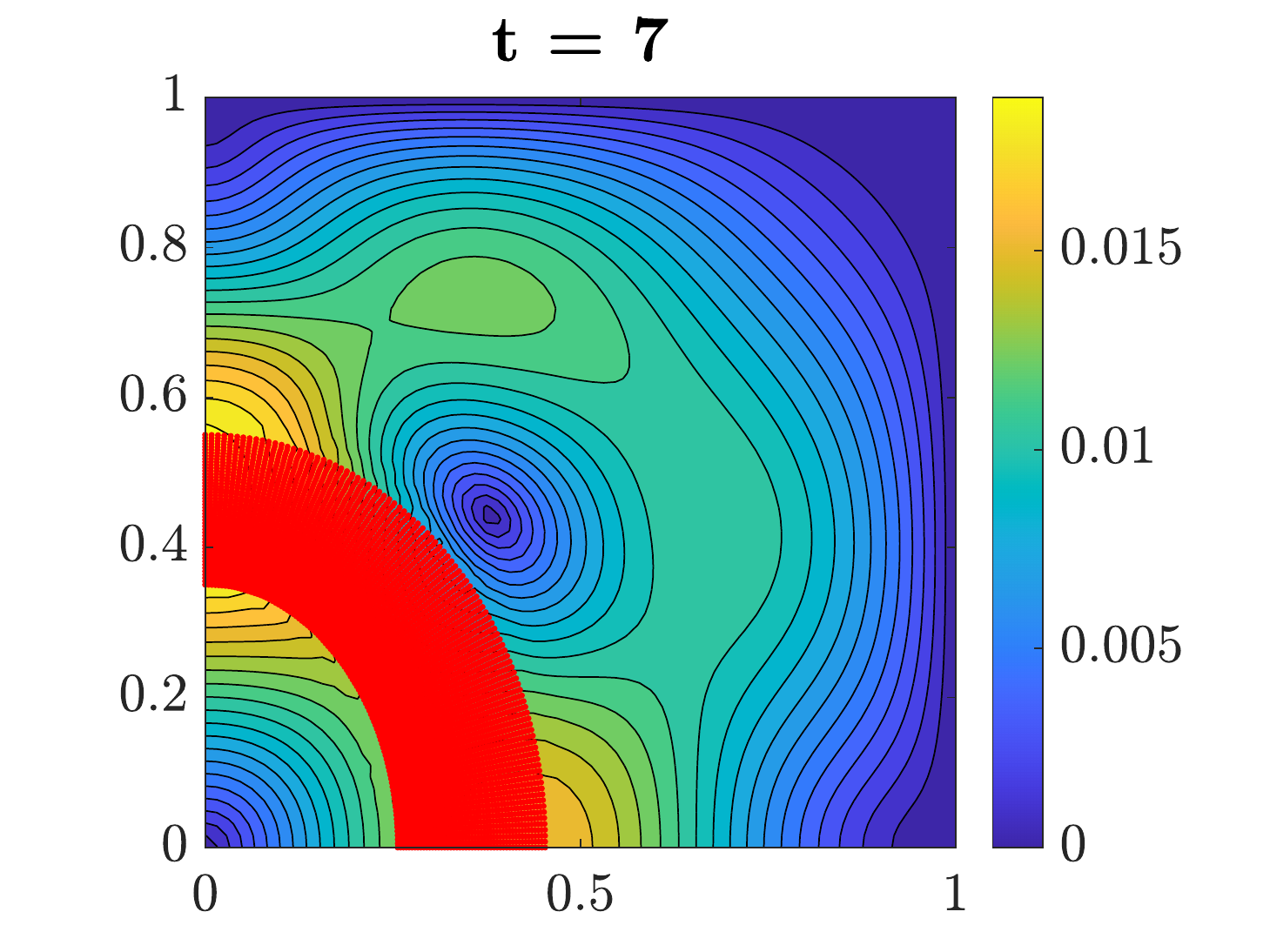}
		\includegraphics[trim = 30 10 20 0,clip, width=4cm]{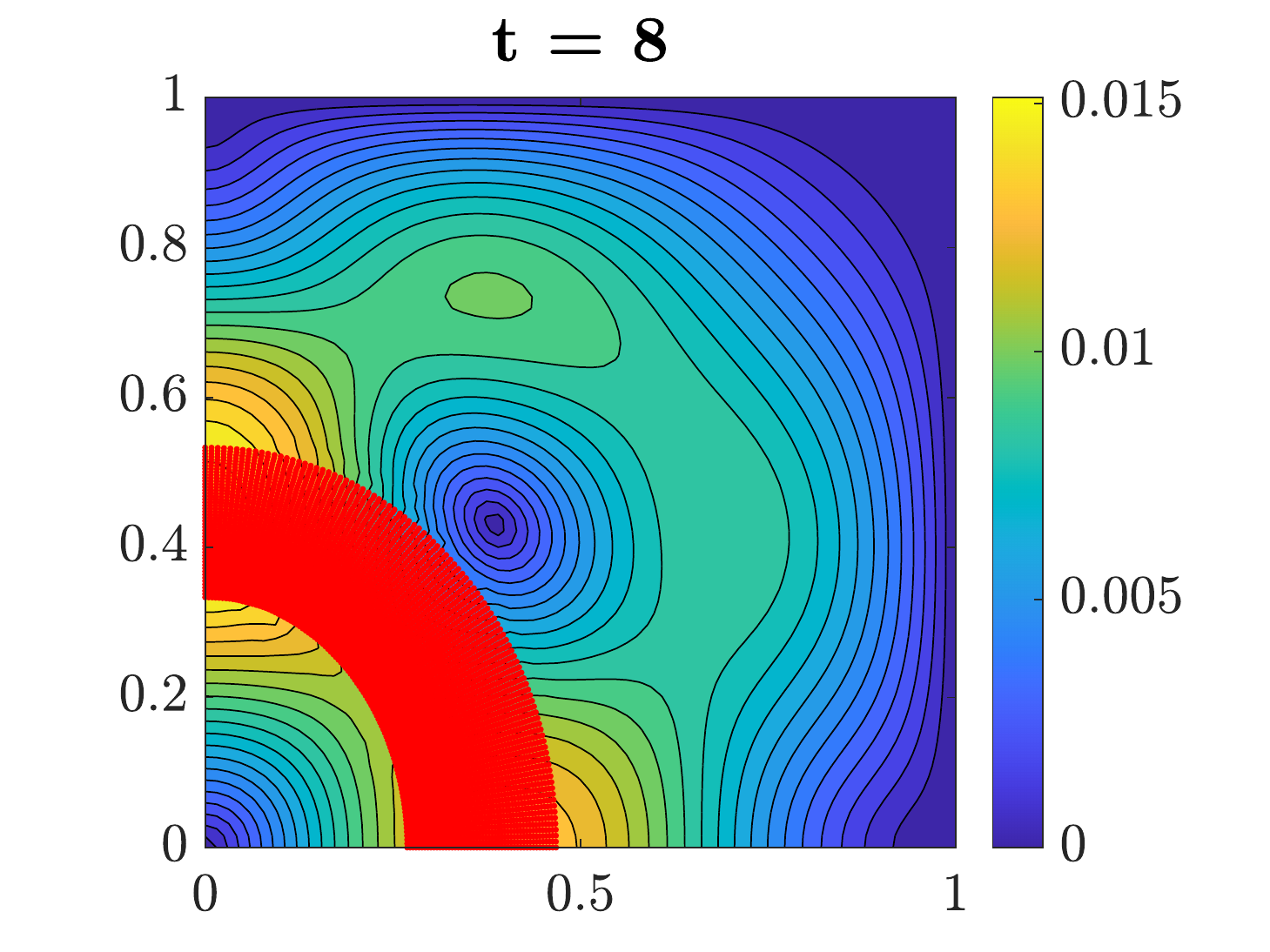}
		\includegraphics[trim = 30 10 20 0,clip, width=4cm]{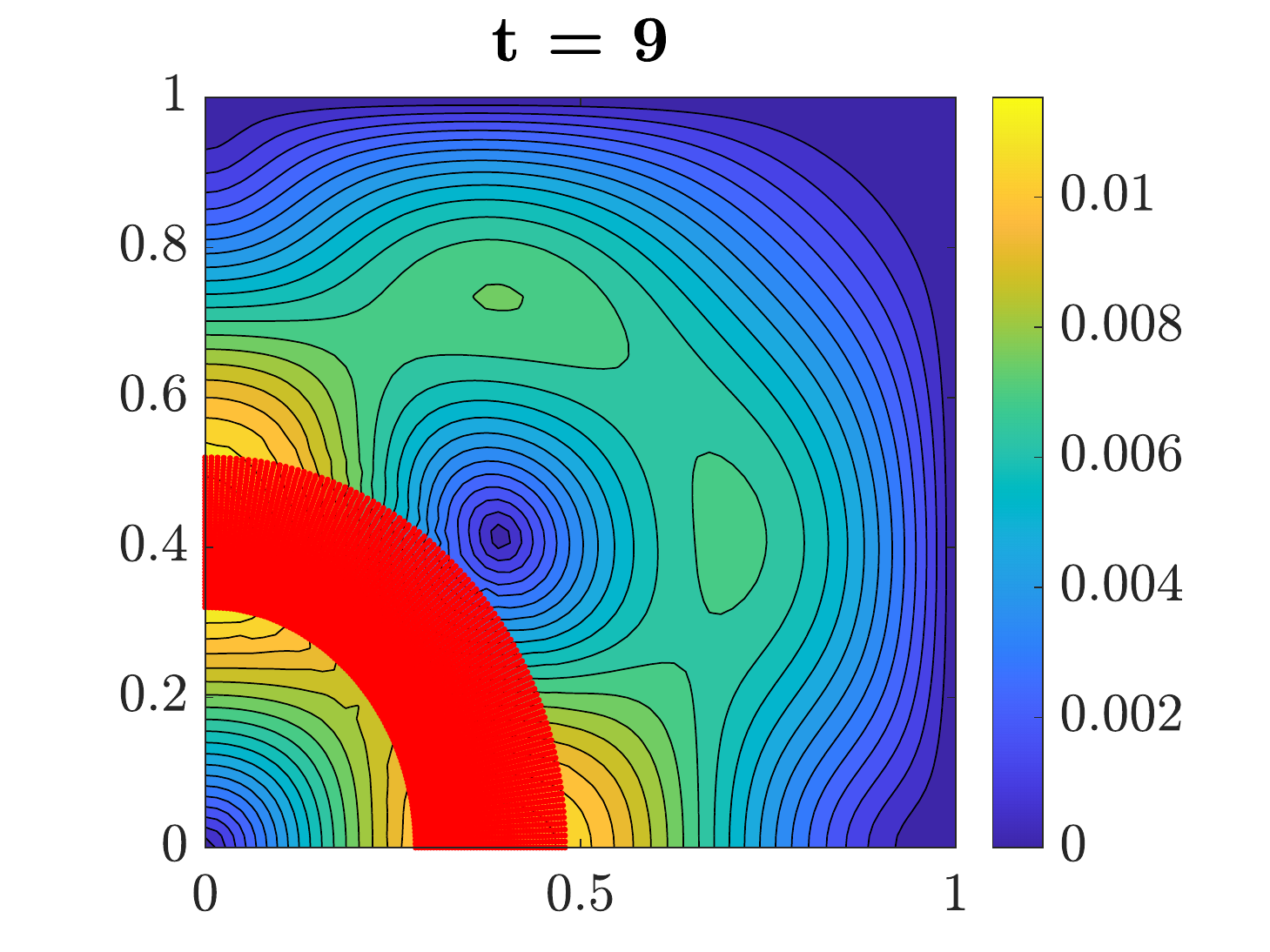}
		\includegraphics[trim = 30 10 20 0,clip, width=4cm]{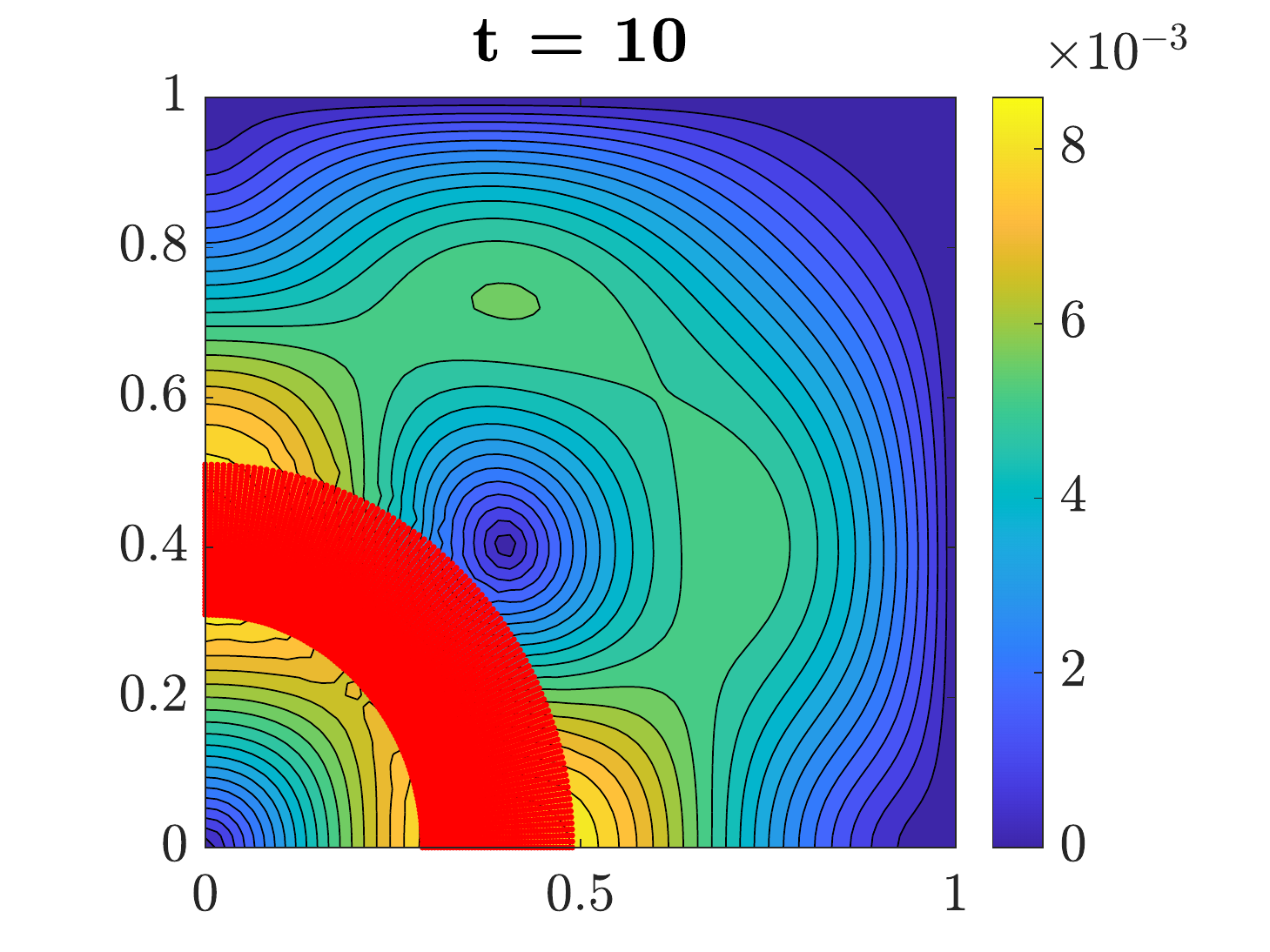}\\
		\medskip
		\includegraphics[trim = 30 10 20 0,clip, width=4cm]{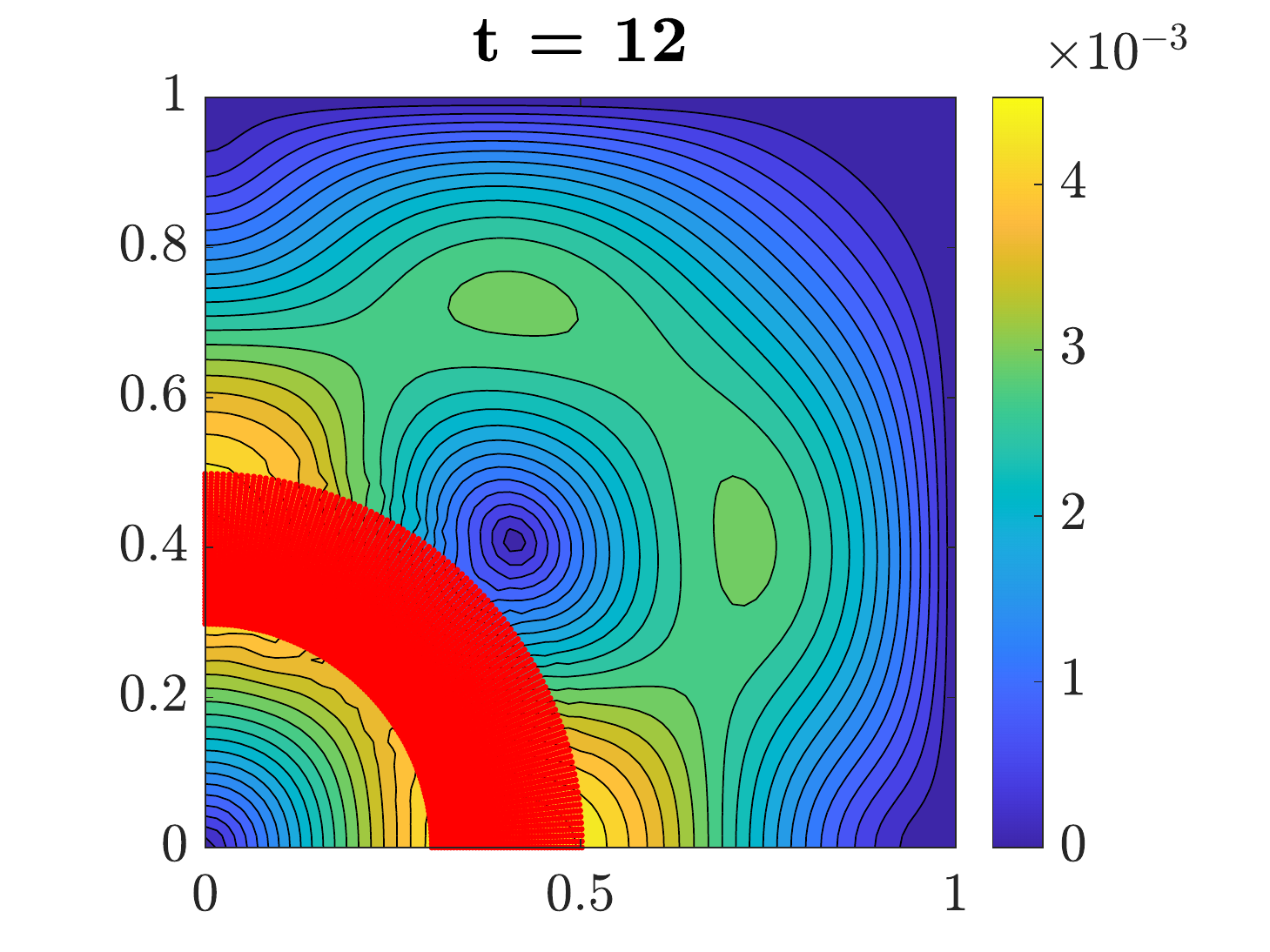}
		\includegraphics[trim = 30 10 20 0,clip, width=4cm]{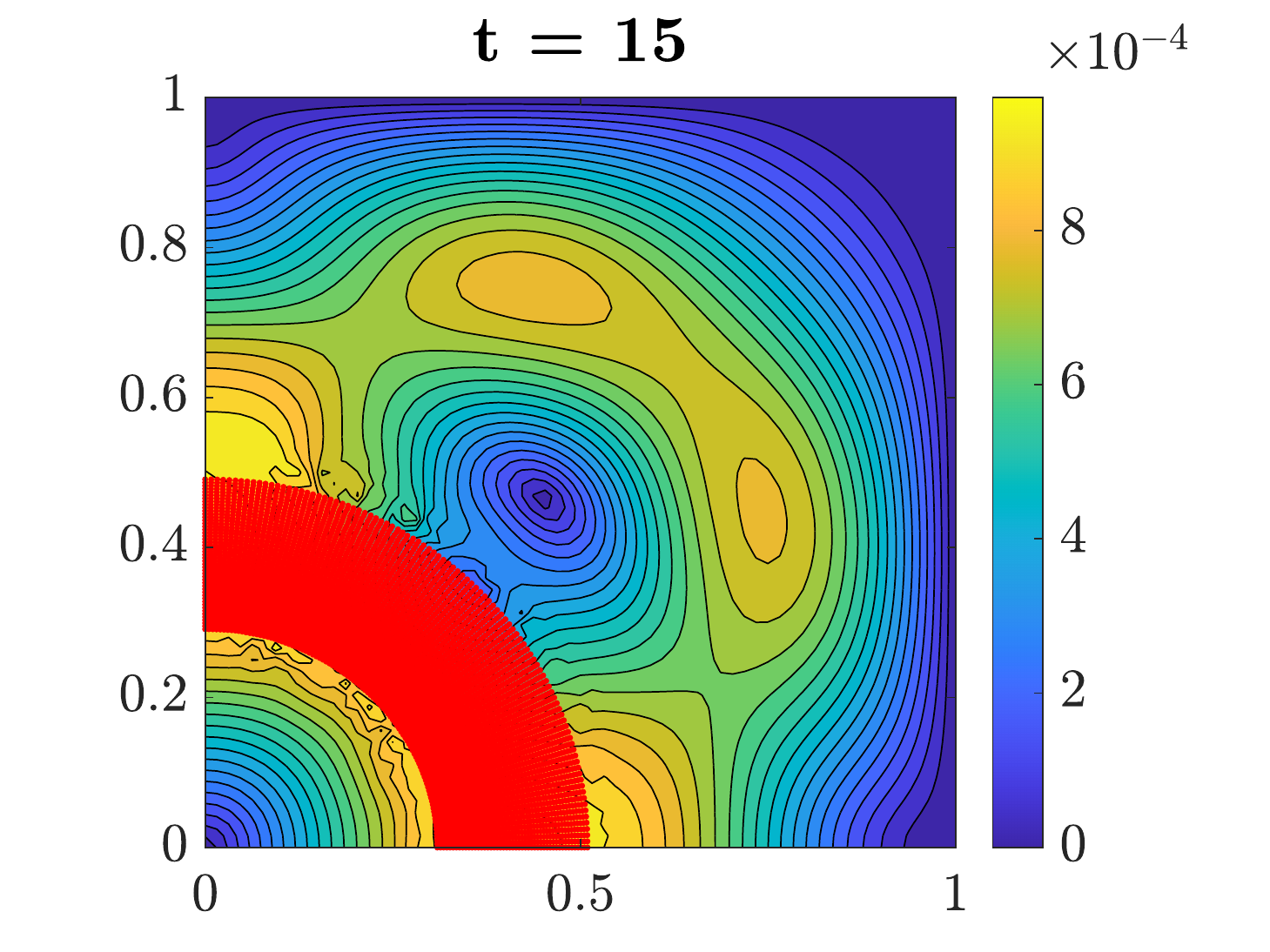}
		\includegraphics[trim = 30 10 20 0,clip, width=4cm]{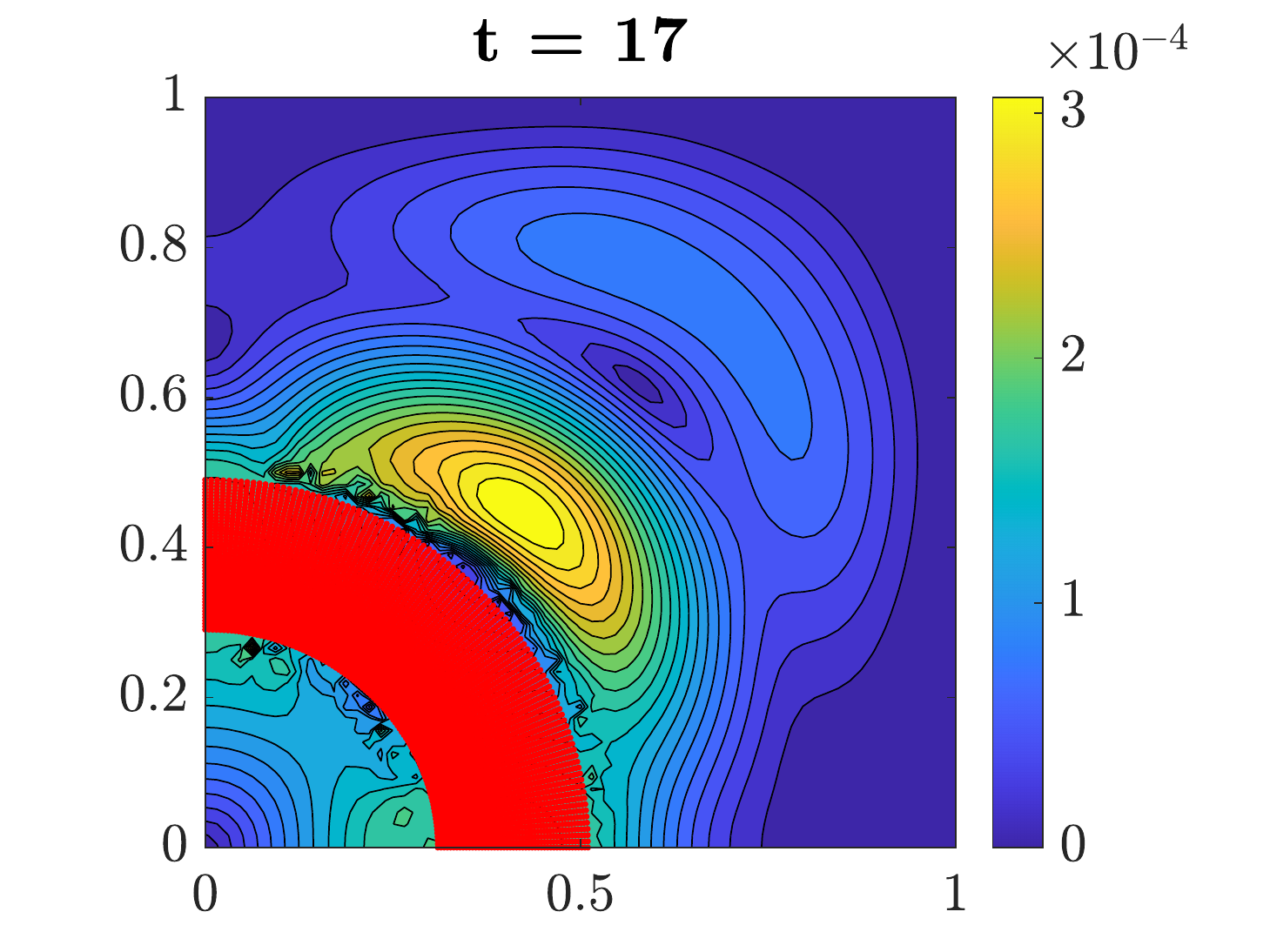}
		\includegraphics[trim = 30 10 20 0,clip, width=4cm]{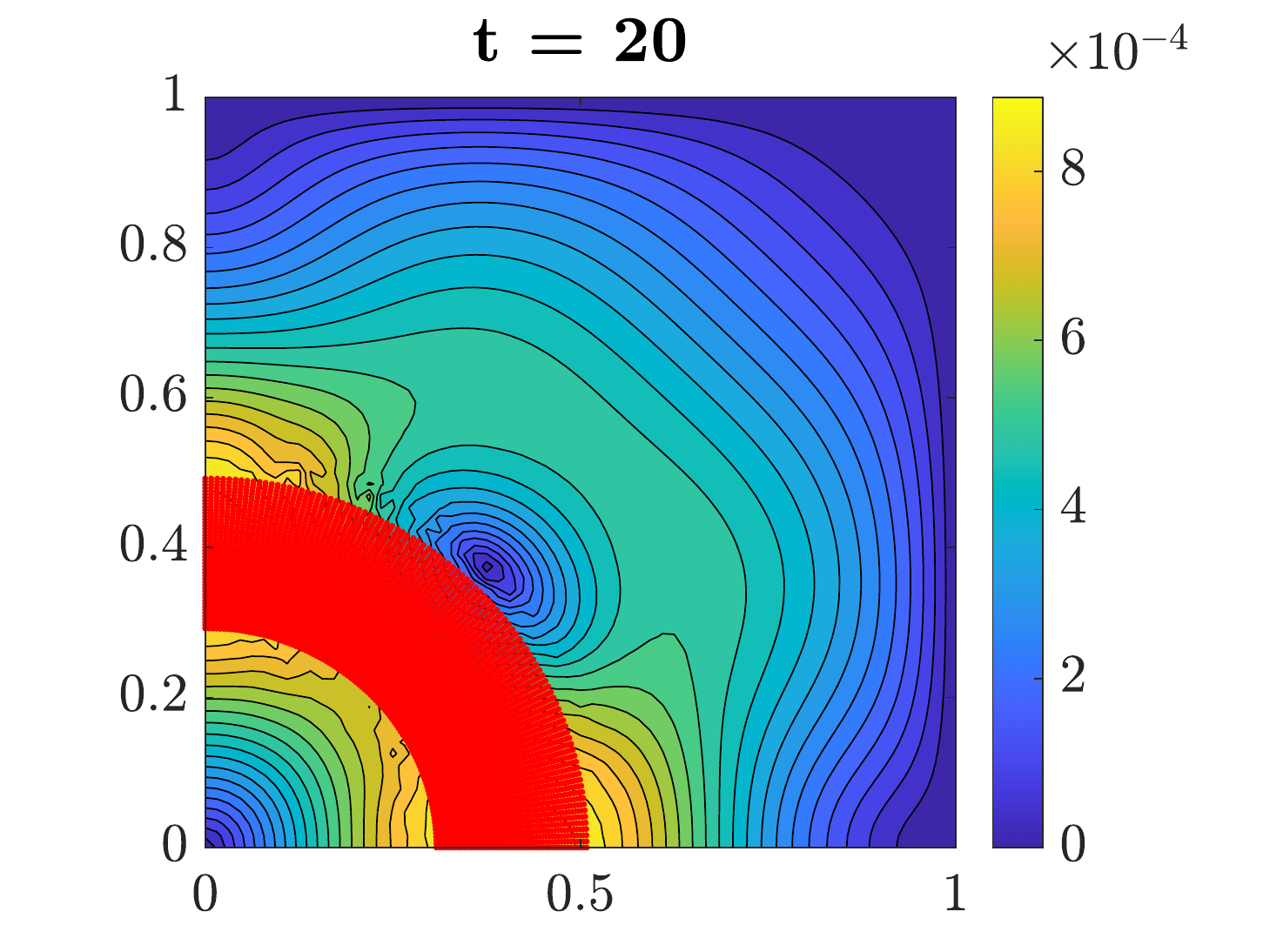}\\
		\caption{\rv Simulation of the annulus with linear constitutive law: some snapshots. The structure position is represented in red, while streamlines and color bars refer to the velocity. At the initial condition, the annulus is stretched. When the simulation starts, the structure is released and internal elastic forces bring it back to its resting configuration. \vr}%{Linear solid model: snapshots of structure evolution.}
		\label{fig:annulus_evolution}
	\end{center}
\end{figure}

We recall that both solid and fluid are incompressible so that the structure and the fluid inside it should theoretically maintain a constant volume. We study the evolution of the system during the time interval $[0,2]$. Some snapshots are presented in Figure \ref{fig:annulus_evolution}.

\begin{table}
	\begin{center}
		\begin{tabular}{r|r|r|r|r|r|r|r|r|r}
			\hline
			\multicolumn{10}{c}{{\bf Linear solid model -- Mesh refinement test}} \\
			\multicolumn{10}{c}{{\textit{Coupling with mesh intersection}}} \\
			\hline
			\multicolumn{10}{c}{procs = 128, T = 2, $\dt$ = 0.01} \\
			\hline
			dofs    & vol. loss (\%)        & $T_{ass}(s)$  & $T_{coup}(s)$ & \multicolumn{3}{c|}{block-diag} 		& \multicolumn{3}{c}{block-tri} \\
			&                       &               &               & its           & $T_{sol}(s)$	& $T_{tot}(s)$	& its           & $T_{sol}(s)$	& $T_{tot}(s)$ \\
			\hline
			30534   & 1.17		& 2.56e-3	& 2.09e-1	& 12 	& 1.38 		& 3.18e+2	& 7	& 8.08e-1 	& 2.03e+2 \\
			120454  & 5.06e-1  	& 1.02e-2	& 8.94e-1	& 31	& 4.69 		& 1.12e+3	& 9	& 1.34		& 4.43e+2 \\
			269766  & 3.18e-1  	& 2.35e-2	& 3.10		& 86	& 16.12		& 3.85e+3 	& 11& 2.06		& 1.04e+3\\
			478470  & 2.33e-1  	& 4.01e-2	& 8.71		& 160	& 37.53		& 9.25e+3	& 12& 3.02		& 2.30e+3\\
			746566  & 1.85e-1  	& 6.50e-2	& 19.52		& 394	& 1.20e+2	& 2.79e+4	& 13& 4.33		& 4.72e+3 \\
			1074054 & 1.54e-1	& 9.44e-2	& 37.72		& -		& -			& -			& 14& 5.51		& 8.65e+3\\
			1460934	& 1.27e-1	& 1.27e-1	& 67.10		& -		& -			& -			& 16& 7.18			& 1.49e+4\\
			\hline
		\end{tabular}
		\vspace*{2mm}
		\caption{Test 1, refining the mesh in the linear solid model, coupling with mesh intersection. The simulations are run on the Shaheen cluster. procs = number of processors; dofs = degrees of freedom; vol. loss = loss of structure volume in percentage; $T_{ass}$ = CPU time to assemble the stiffness and mass matrices; $T_{coup}$ = CPU time to assemble the coupling term; its = GMRES iterations; $T_{sol}$ = CPU time to solve the linear system; $T_{tot}$ = total simulation CPU time. The quantities $T_{coup}$, its and $T_{sol}$ are averaged over the time steps. All CPU times are reported in seconds.}
		\label{lin_tab_opti_inters}
	\end{center}
\end{table}

\begin{table}
	\begin{center}
		\begin{tabular}{r|r|r|r|r|r|r|r|r|r}
			\hline
			\multicolumn{10}{c}{{\bf Linear solid model -- Mesh refinement test}} \\
			\multicolumn{10}{c}{{\textit{Coupling without mesh intersection}}} \\
			\hline
			\multicolumn{10}{c}{procs = 128, T = 2, $\dt$ = 0.01} \\
			\hline
			dofs    & vol. loss (\%)        & $T_{ass}(s)$  & $T_{coup}(s)$ & \multicolumn{3}{c|}{block-diag} 		& \multicolumn{3}{c}{block-tri} \\
			&                       &               &               & its           & $T_{sol}(s)$	& $T_{tot}(s)$	& its           & $T_{sol}(s)$	& $T_{tot}(s)$ \\
			\hline
			30534   & 6.99e-2	& 2.50e-3	& 1.68e-1	& 9	& 1.04 	& 2.42e+2	& 5	& 6.57e-1  	& 1.65e+2 \\
			120454  & 6.89e-2	& 9.06e-3	& 2.50e-1	& 9	& 1.53 	& 3.57e+2	& 6	& 1.01	 	& 2.53e+2 \\
			269766  & 4.87e-2	& 2.33e-2	& 9.96e-1	& 10& 2.10	& 6.19e+2	& 6	& 1.31	 	& 4.65e+2\\
			478470  & 4.24e-2	& 4.13e-2	& 3.70		& 10& 2.65	& 1.27e+3	& 6	& 1.63	 	& 1.06e+3 \\
			746566  & 4.09e-2	& 6.50e-2	& 9.90		& 10& 3.30	& 2.64e+3	& 6	& 1.97	 	& 2.25e+3\\
			1074054 & 3.69e-2	& 9.46e-2	& 20.68		& 10& 3.93	& 4.92e+3	& 6	& 2.55	 	& 4.66e+3 \\
			1460934	& 3.52e-2	& 1.29e-1	& 45.39		& 10& 4.66	& 1.00e+4	& 6	& 3.15		& 9.86e+3\\
			\hline
		\end{tabular}
		\vspace*{2mm}
		\caption{Test 1, refining the mesh in the linear solid model, coupling without mesh intersection. The simulations are run on the Shaheen cluster. Same format as in Table \ref{lin_tab_opti_inters}.}
		\label{lin_tab_opti_nointers}
	\end{center}
\end{table}

\subsection{Mesh refinement test}

We study the behavior of the parallel FSI solver with respect to mesh refinement, by keeping fixed the time step size $\Delta t = 0.01$, 
the number of processors $procs = 128$ and the final time $T = 2$.
The number of degrees of freedom (dofs) varies from 30534 to 1460934. 

Table \ref{lin_tab_opti_inters} reports the results in case of coupling with mesh intersection. 
We first observe that the volume loss reduces when the mesh is refined and, for the finest grids, it is very small, always below $0.3\%$.
With respect to the number of dofs, $T_{ass}$ exhibits a moderate increase, whereas the growth of $T_{coup}$ is superlinear. 
The block-diag solver is not robust with respect to mesh refinement, because the GMRES iteration count increases when the mesh is refined.
In the two finest cases, the solver fails.
Thus, $T_{sol}$ exhibits a significant growth with a large number of dofs.
The block-tri solver, instead, is robust with respect to mesh refinement, because the GMRES iteration count remains bounded when the mesh is refined.
As a consequence, $T_{sol}$ exhibits a moderate growth with the number of dofs: for dofs = 746566, $T_{sol}$ is about 27 times smaller than its counterpart for block-diag.%. In the particular case dofs = 746566, block-tri is about 27 times smaller than block-diag. 

Table \ref{lin_tab_opti_nointers} reports the results in case of coupling without mesh intersection.
Even in this case, the volume loss reduces when the mesh is refined, and it remains always below $0.1\%$.
As before, $T_{ass}$ exhibits a moderate increase with respect to the number of dofs, whereas the growth of $T_{coup}$ is superlinear.
Differently from the case of coupling with mesh intersection, the block-diag solver is robust with respect to mesh refinement,
since the GMRES iteration count remains bounded when the number of dofs increases. Consequently, $T_{sol}$ slightly increases with the number of dofs.
The block-tri solver shows a small reduction of GMRES iteration count with respect to block-diag. 
However, the two solvers are comparable in terms of CPU times.

\subsection{Strong scalability test}

We now study the strong scalability of the parallel FSI solver, by keeping fixed the fluid and solid meshes to $128\times 128$ and $384\times 192$ elements respectively. \lg This choice produces \gl a total amount of 478470 dofs. The time parameters are fixed at $\Delta t = 0.01$ and $T = 2$. 
The number of processors increases from 4 to 256. The parallel speedup $S^p$ is computed with respect to the 4 processors run.

Table \ref{lin_tab_scal} reports the results in case of coupling with mesh intersection. We first observe an excellent scalability of the 
assembly phase of the coupling term ($T_{coup}$). The block-diag solver exhibits a scalable behavior in terms of GMRES iterations, 
which remain bounded when the number of processors increases. However, the solution time $T_{sol}$ is not scalable. 
Consequently, the global performance of the solver is impaired, with very slow speedup values (about 7.5 instead of 64 in case of 256 processors). 
The block-tri solver exhibits very low GMRES iteration counts (always $<14$) and solution times $T_{sol}$. 
Also in this case, $T_{sol}$ is not scalable, but, being very small in comparison with $T_{coup}$, which instead is scalable,
the solver results to be globally scalable, with speedup values comparable or larger than the ideal ones. 

Table \ref{lin_tab_scal_nointers} reports the results in case of coupling without mesh intersection.
We observe again a scalable behavior of the assembly phase of the coupling term.
In this case, both the block-diag and block-tri solvers are scalable in terms of GMRES iterations, which remain bounded
when the number of processors increases. Although the solution times $T_{sol}$ do not exhibit a scalable behavior, 
thanks to the excellent scalability of $T_{coup}$, the solvers result to be globally scalable, 
with speedup comparable of larger than the ideal ones.

\begin{figure}[h!]
	\vspace{5mm}
	\begin{center}
		%		\includegraphics[trim = 30 10 30 0,clip, width=8cm]{figures/fig_lin/lin_int_1}
		%		\includegraphics[trim = 30 10 30 0,clip, width=8cm]{figures/fig_lin/lin_int_2}\\
		%		\medskip
		%		\includegraphics[trim = 30 10 30 0,clip, width=8cm]{figures/fig_lin/lin_int_3}
		\begin{overpic}[width=7.95cm]{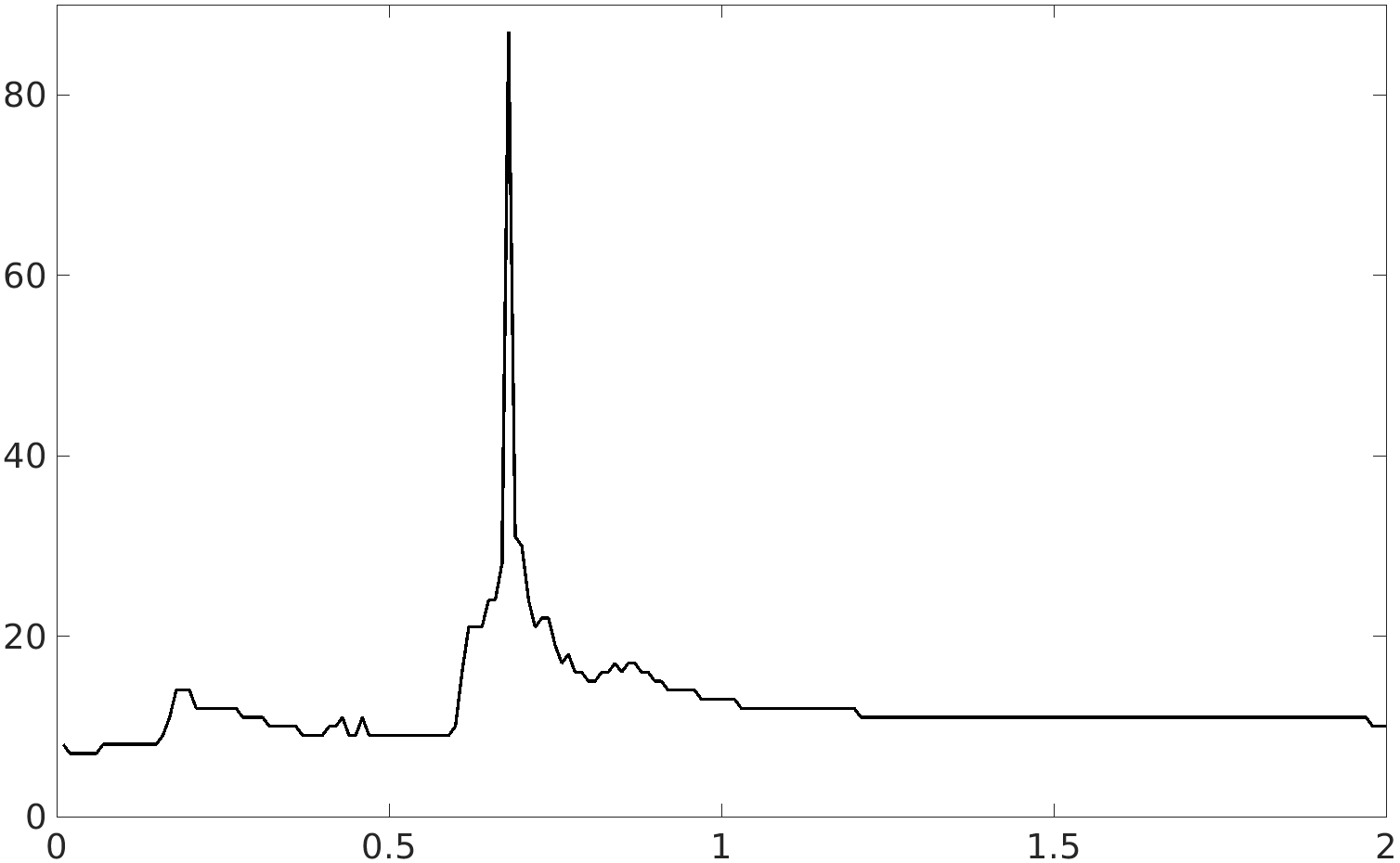}
			\put (48,-5) {\scriptsize time}
			\put (37,64) {\scriptsize \textbf{Linear iterations}}
			\put (-5,30) {\begin{sideways}
					\scriptsize $its$
			\end{sideways}}
		\end{overpic}
		\quad
		\begin{overpic}[width=8cm]{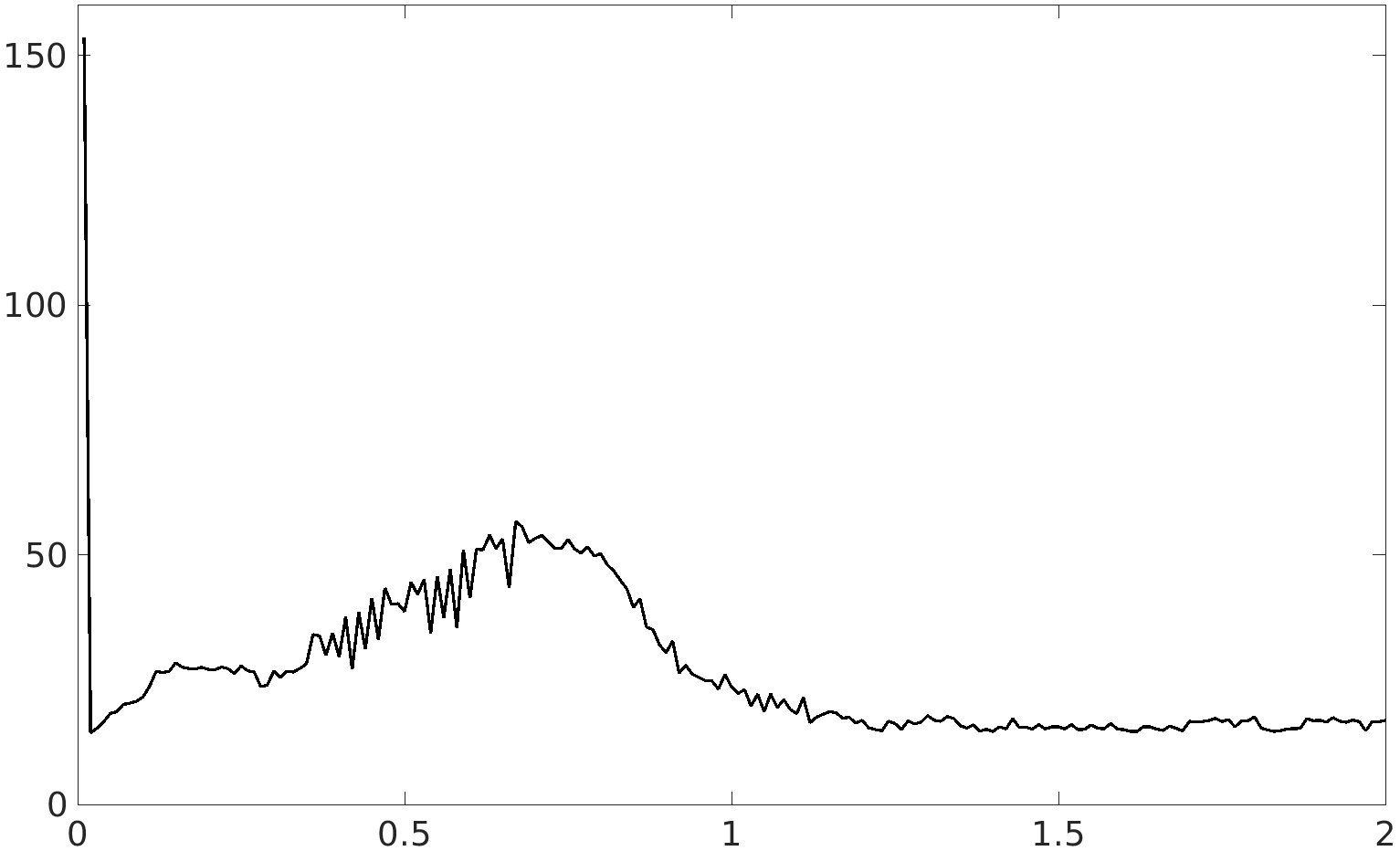}
			\put (48,-5) {\scriptsize time}
			\put (15,63) {\scriptsize \textbf{CPU time to assemble the coupling term}}
			\put (-5,25) {\begin{sideways}
					\scriptsize $T_{coup}(s)$
			\end{sideways}}
		\end{overpic}\\
		\vspace{10mm}
		\begin{overpic}[width=8cm]{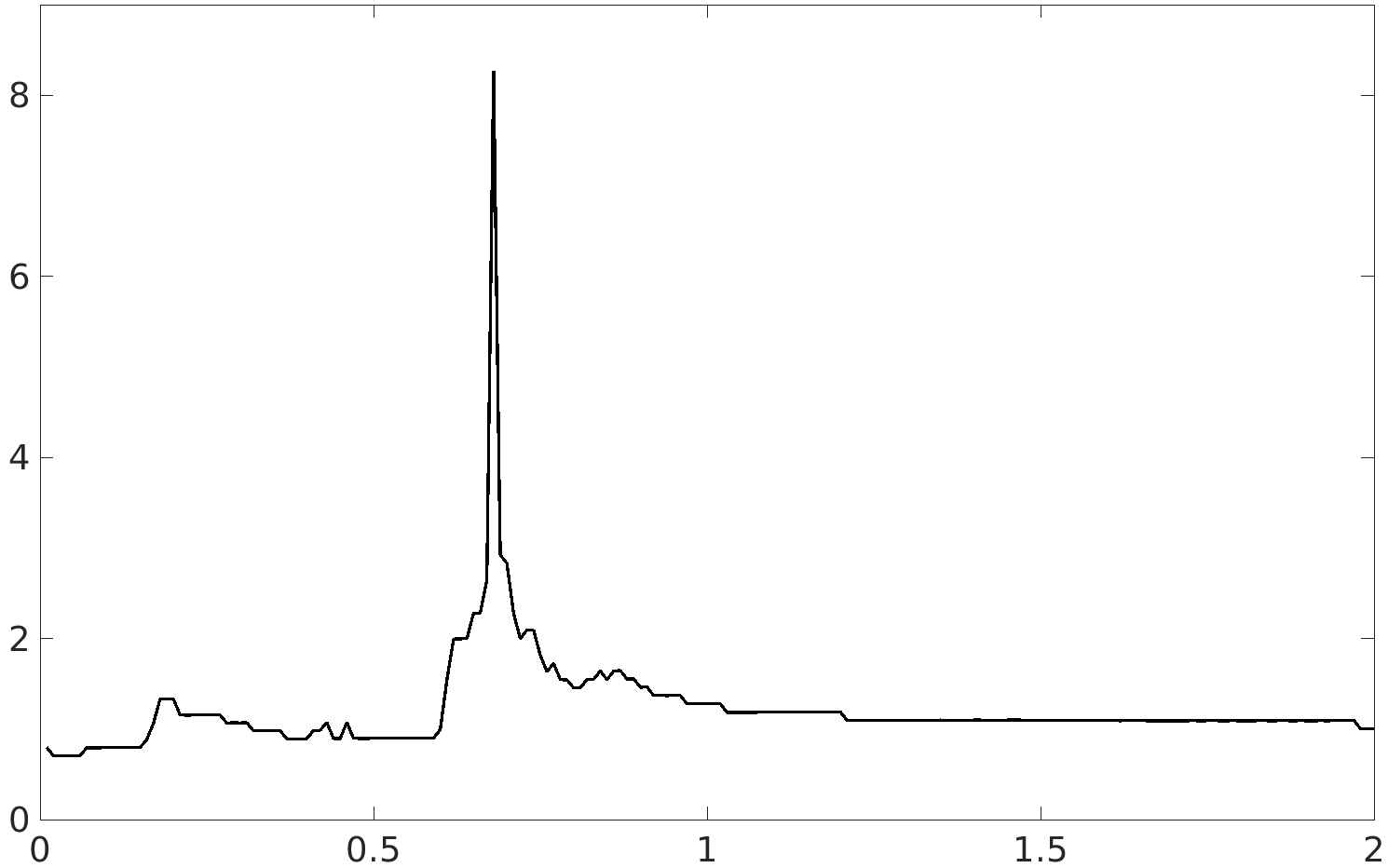}
			\put (48,-5) {\scriptsize time}
			\put (17,65) {\scriptsize \textbf{CPU time to solve the linear system}}
			\put (-5,25) {\begin{sideways}
					\scriptsize $T_{sol}(s)$
			\end{sideways}}
		\end{overpic}
		\vspace{3mm}
		\caption{Test 2, strong scalability in the linear solid model, coupling with mesh intersection. Time evolution of linear iterations, CPU time to assemble the coupling term and to solve the linear system. Run on 32 processors of Shaheen cluster, 478470 dofs and block-tri preconditioner.}
		\label{fig:time_evol_inter}
	\end{center}
\end{figure}

\begin{figure}[h!]
	\vspace{5mm}
	\begin{center}
		%		\includegraphics[trim = 30 10 30 0,clip, width=8cm]{figures/fig_lin/lin_noint_1}
		%		\includegraphics[trim = 30 10 30 0,clip, width=8cm]{figures/fig_lin/lin_noint_2}\\
		%		\medskip
		%		\includegraphics[trim = 30 10 30 0,clip, width=8cm]{figures/fig_lin/lin_noint_3}
		%
		%		\includegraphics[width=8cm]{figures/fig_lin/lin_noint_1.png}
		%		\includegraphics[width=8.2cm]{figures/fig_lin/lin_noint_2.png}\\
		%		\medskip
		%		\includegraphics[width=8cm]{figures/fig_lin/lin_noint_3.png}
		%
		\begin{overpic}[width=7.95cm]{figures/fig_lin/lin_int_1.png}
			\put (48,-5) {\scriptsize time}
			\put (37,64) {\scriptsize \textbf{Linear iterations}}
			\put (-5,30) {\begin{sideways}
					\scriptsize $its$
			\end{sideways}}
		\end{overpic}
		\quad
		\begin{overpic}[width=8cm]{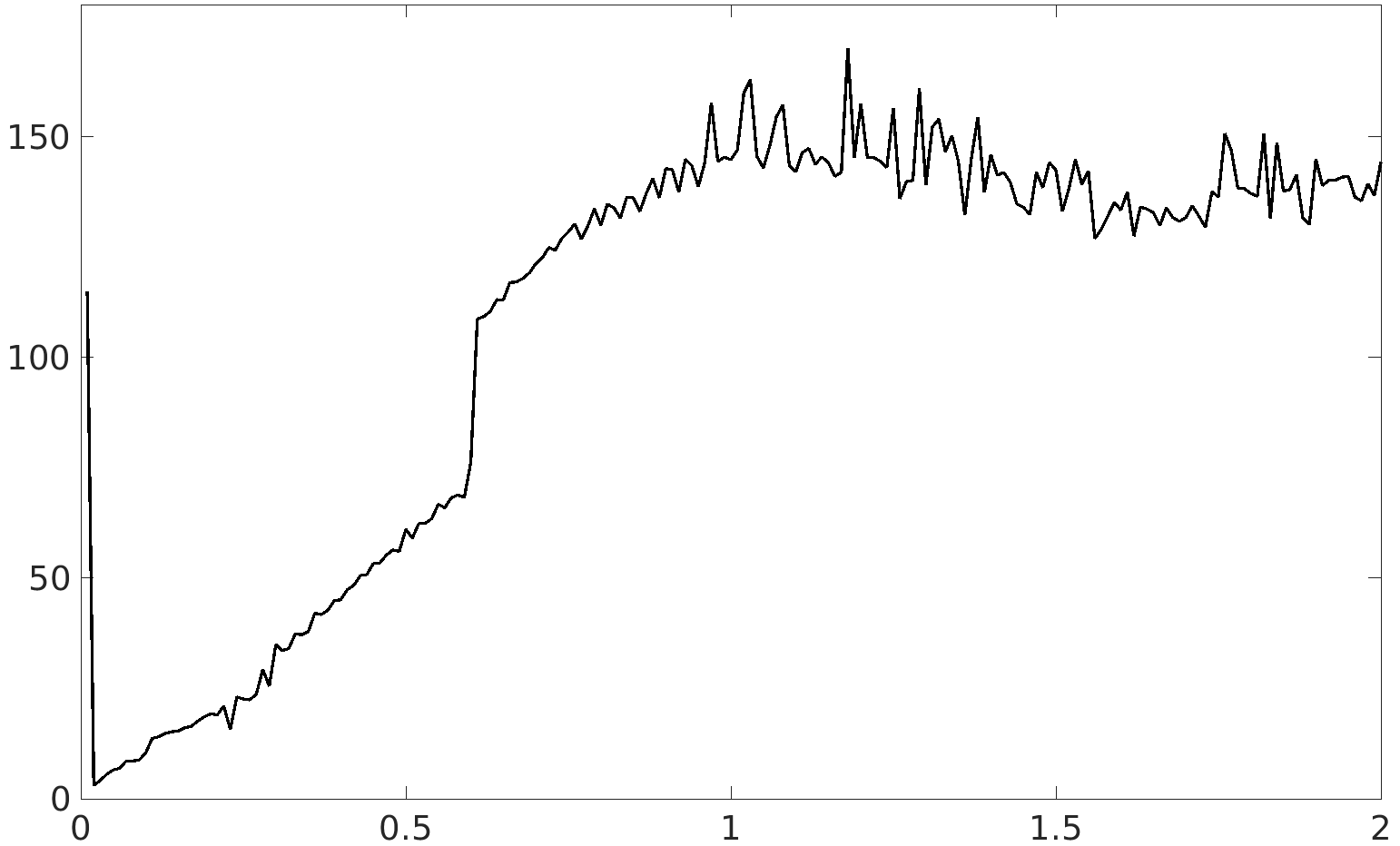}
			\put (48,-5) {\scriptsize time}
			\put (15,63) {\scriptsize \textbf{CPU time to assemble the coupling term}}
			\put (-5,25) {\begin{sideways}
					\scriptsize $T_{coup}(s)$
			\end{sideways}}
		\end{overpic}\\
		\vspace{10mm}
		\begin{overpic}[width=8cm]{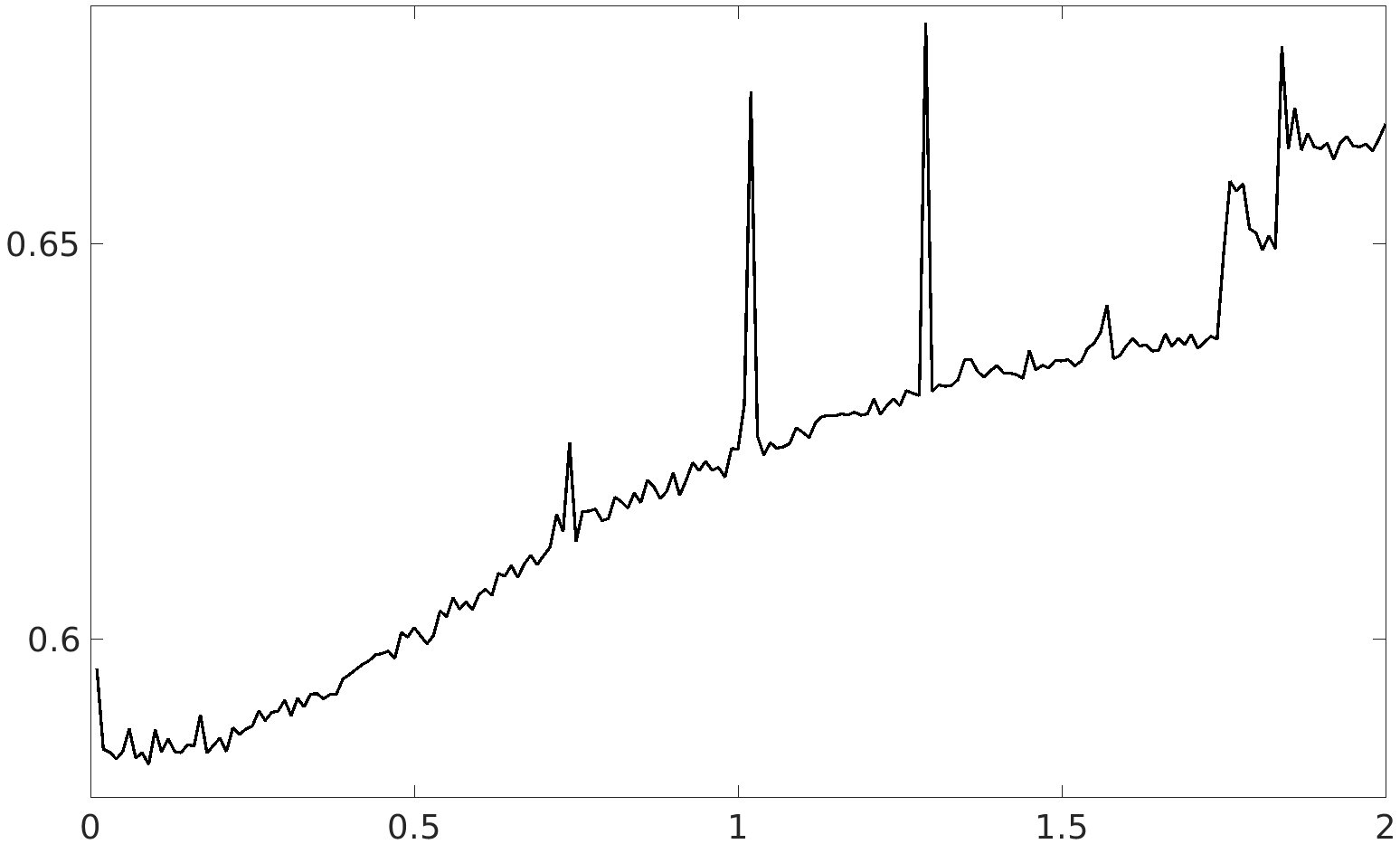}
			\put (48,-5) {\scriptsize time}
			\put (17,63) {\scriptsize \textbf{CPU time to solve the linear system}}
			\put (-5,25) {\begin{sideways}
					\scriptsize $T_{sol}(s)$
			\end{sideways}}
		\end{overpic}
		\vspace{3mm}
		\caption{Test 2, strong scalability in the linear solid model, coupling without mesh intersection. Time evolution of linear iterations, CPU time to assemble the coupling term and to solve the linear system. Run on 32 processors of Shaheen cluster, 478470 dofs and block-tri preconditioner.}
		\label{fig:time_evol_nointer}
	\end{center}
\end{figure}

In Figures \ref{fig:time_evol_inter} and \ref{fig:time_evol_nointer}, we report the evolution in time of the number of linear iterations and CPU time (in seconds) to assemble the coupling matrix and solve the linear system respectively in the case of the linear solid model discretized with 478470 dofs and solved with block-tri preconditioner on 32 processors of Shaheen cluster.

% Test 2 - Marconi200 case on Shaheen - intersection
\begin{table}
	\begin{center}
		\begin{tabular}{r|r|r|r|r|r|r|r|r|r|r}
			\hline
			\multicolumn{11}{c}{{\bf Linear solid model -- Strong scalability test}} \\
			\multicolumn{11}{c}{{\textit{Coupling with mesh intersection}}} \\
			\hline
			\multicolumn{11}{c}{dofs = 478470, T = 2, $\dt$ = 0.01} \\
			\hline
			procs   & $T_{ass}(s)$  & $T_{coup}(s)$ & \multicolumn{4}{c|}{block-diag}                               & \multicolumn{4}{c}{block-tri} \\
			&               &               & its           & $T_{sol}(s)$  & $T_{tot}$ (s) & $S^p$         & its           & $T_{sol}(s)$  & $T_{tot}$ (s) & $S^p$\\
			\hline
			4		& 6.41e-1	& 498.08	& 192	& 41.04		& 1.08e+5	& -			& 12	& 2.52		& 9.85e+4	& -		\\
			8		& 3.47e-1	& 169.69	& 168	& 21.60		& 3.83e+4	& 2.82 (2)	& 13	& 1.49		& 3.42e+4	& 2.88 (2)		\\
			16		& 1.78e-1	& 89.18		& 180	& 18.09		& 2.15e+4	& 5.02 (4)	& 12	& 1.20		& 1.79e+4	& 5.50 (4)		\\
			32		& 1.65e-1	& 26.45		& 192	& 19.42		& 9.17e+3	& 11.78 (8)	& 12	& 1.22		& 5.53e+3	& 17.81 (8)		\\
			64		& 8.32e-2	& 17.38 	& 165	& 23.09		& 8.09e+3	& 13.35 (16)& 13	& 1.57		& 3.64e+3	& 27.06 (16)		\\
			128		& 4.12e-2	& 8.52		& 170	& 40.58		& 9.82e+3	& 11.00 (32)& 12	& 3.02		& 2.35e+3	& 41.91 (32)		\\
			256 &	2.03e-2		& 4.08		& 144	& 68.31		& 1.45e+4	& 7.45 (64)	& 12	& 5.91	& 2.00e+3		& 49.25 (64)
			\\
			\hline
		\end{tabular}
		\vspace*{2mm}
		\caption{Test 2, strong scalability in the linear solid model, coupling with mesh intersection. The simulations are run on the Shaheen cluster. dofs = degrees of freedom; procs = number of processors; $T_{ass}$ = CPU time to assemble the stiffness and mass matrices; $T_{coup}$ = CPU time to assemble the coupling term; its = GMRES iterations; $T_{sol}$ = CPU time to solve the linear system; $T_{tot}$ = total simulation time; $S^p$ = parallel speedup computed with respect to the 4 processors run. The quantities $T_{coup}$, its and $T_{sol}$ are averaged over the time steps. All CPU times are reported in seconds.}
		\label{lin_tab_scal}
	\end{center}
\end{table}

% FIGURA

% Test 2 - Marconi200 case on Shaheen - no intersection
\begin{table}
	\begin{center}
		\begin{tabular}{r|r|r|r|r|r|r|r|r|r|r}
			\hline
			\multicolumn{11}{c}{{\bf Linear solid model -- Strong scalability test}} \\
			\multicolumn{11}{c}{{\textit{Coupling with mesh intersection}}} \\
			\hline
			\multicolumn{11}{c}{dofs = 478470, T = 2, $\dt$ = 0.01} \\
			\hline
			procs   & $T_{ass}(s)$  & $T_{coup}(s)$ & \multicolumn{4}{c|}{block-diag}                               & \multicolumn{4}{c}{block-tri} \\
			&               &               & its           & $T_{sol}(s)$  & $T_{tot}$ (s) & $S^p$         & its           & $T_{sol}(s)$  & $T_{tot}$ (s) & $S^p$\\
			\hline
			4	& 6.77e-1	& 1.68e+3	& 10	& 2.60	& 3.36e+5	&	-		& 6		& 1.69		& 3.37e+5	& - 		\\
			8	& 3.89e-1	& 741.99	& 10	& 1.87		& 1.50e+5	& 2.24 (2)	& 6		& 1.22		& 1.49e+5	& 2.26 (2)		\\
			16	& 1.79e-1	& 287.02	& 10	& 9.76e-1	& 5.76e+4	& 5.83 (4)	& 6		& 6.35e-1	& 5.75e+4	& 5.86 (4)		\\
			32	& 1.65e-1	& 108.92	& 10	& 9.95e-1	& 2.20e+4	& 15.27 (8)	& 6		& 6.21e-1	& 2.17e+4	& 15.53 (8)		\\
			64	& 7.70e-2	& 28.66		& 10	& 1.29		& 5.99e+3	& 56.09 (16)& 6		& 8.68e-1	& 6.02e+3	& 55.98 (16)		\\
			128	& 4.04e-2	& 4.32		& 10	& 87.20		& 1.83e+4	& 18.36 (32)& 6		& 58.70		& 1.26e+4	& 26.75 (32)		\\
			256 & 2.05e-2	& 1.63		& 10	& 5.04		& 1.33e+3	& 252.63(64)& 6		& 3.05		& 951.11	& 354.32 (64)
			\\
			\hline
		\end{tabular}
		\vspace*{2mm}
		\caption{Test 2, strong scalability in the linear solid model, coupling without mesh intersection. The simulations are run on the Shaheen cluster. Same format as in Table \ref{lin_tab_scal}.}
		\label{lin_tab_scal_nointers}
	\end{center}
\end{table}

% FIGURA

%\subsection{\lg Time step refinement \gl}
\subsection{Time step refinement}

In this test, we study the behavior of our solver with respect to the refinement of the time step.

We fix the number of processors at 64 and $T=2$. Also the mesh sizes are fixed and we consider the same spatial discretizations we used for the strong scalability test: $128\times 128$ elements for the fluid and $384\times 192$ for the solid.

Table \ref{lin_tab_time_inters} reports the results in the case of coupling with mesh intersection. We can see that the time step refinement affects the volume loss, which decreases. Also $T_{coup}$ decreases with respect to the time step refinement because the position of the solid body is affected by the precision we use in time. On the other hand, since the physical parameters are fixed, $T_{ass}$ is constant, as expected. For a comparison between the two preconditioners, we can say that block-tri is more robust than block-diag: the number of iterations is lower in the first case and it is bounded by 19, while in the second case, when $\dt=0.02$, the number of iterations is higher than 1000. Clearly, this fact has consequences on the time $T_{sol}$ we need to solve the linear system.

In Table \ref{lin_tab_time_nointers}, we summarize the results for $L_f(\X_h^n)$ assembled without mesh intersection. As in the previous case, the volume loss decreases when $\dt$ is refined, $T_{ass}$ is constant and $T_{coup}$ decreases. In particular, we can also notice that the assembly with mesh intersection is faster than the procedure with inexact integration. Furthermore, the behavior of the two preconditioners is comparable: indeed, the number of iterations for block-diag is bounded by 12, while for block-tri is bounded by~7. In both cases, this number decreases when we refine the time step. As a consequence, the values of $T_{sol}$ are very similar in the two cases.

\begin{table}
	\begin{center}
		\begin{tabular}{r|r|r|r|r|r|r|r|r|r}
			\hline
			\multicolumn{10}{c}{{\bf Linear solid model -- Time step refinement test}} \\
			\multicolumn{10}{c}{{\textit{Coupling with mesh intersection }}} \\
			\hline
			\multicolumn{10}{c}{dofs = 478470, procs = 64, T = 2, $\dt$ = 0.01} \\
			\hline
			$\dt$    & vol. loss (\%)        & $T_{ass}(s)$  & $T_{coup}(s)$ & \multicolumn{3}{c|}{block-diag} 		& \multicolumn{3}{c}{block-tri} \\
			&                       &               &               & its           & $T_{sol}(s)$	& $T_{tot}(s)$	& its           & $T_{sol}(s)$	& $T_{tot}(s)$ \\
			\hline
			0.02 & 2.55e-1 	& 1.01e-1 	& 15.06	& 1364 & 190.66	& 2.05e+4			& 19& 3.56	& 1.86e+3 \\
			0.01 & 2.33e-1	& 1.01e-1	& 11.83 & 170& 23.69& 7.48e+3	& 12& 2.33	& 2.83e+3 \\
			0.005& 2.04e-1	& 1.01e-1	& 10.86	& 30 & 5.70	& 6.68e+3	& 9	& 1.91	& 5.12e+3 \\
			0.002& 1.88e-1	& 1.01e-1	& 9.90	& 12 & 2.61	& 1.26e+4	& 7	& 1.44	& 1.13e+4 \\
			0.005& 1.81e-1	& 1.01e-1	& 9.91	& 8  & 1.73	& 2.28e+4	& 5	& 1.09	& 2.20e+4 \\
			\hline
		\end{tabular}
		\vspace*{2mm}
		\caption{Test 3, refining the time step in the linear solid model, coupling with mesh intersection. The simulations are run on the Shaheen cluster. Same format as in Tables \ref{lin_tab_opti_inters} and \ref{lin_tab_opti_nointers}.}
		\label{lin_tab_time_inters}
	\end{center}
\end{table}

\begin{table}
	\begin{center}
		\begin{tabular}{r|r|r|r|r|r|r|r|r|r}
			\hline
			\multicolumn{10}{c}{{\bf Linear solid model -- Time step refinement test}} \\
			\multicolumn{10}{c}{{\textit{Coupling without mesh intersection }}} \\
			\hline
			\multicolumn{10}{c}{dofs = 478470, procs = 64, T = 2, $\dt$ = 0.01} \\
			\hline
			$\dt$    & vol. loss (\%)        & $T_{ass}(s)$  & $T_{coup}(s)$ & \multicolumn{3}{c|}{block-diag} 		& \multicolumn{3}{c}{block-tri} \\
			&                       &               &               & its           & $T_{sol}(s)$	& $T_{tot}(s)$	& its           & $T_{sol}(s)$	& $T_{tot}(s)$ \\
			\hline
			0.02 & 6.27e-2 	& 1.01e-1 	& 45.78	& 12  & 2.41	& 4.79e+3	& 7 & 1.45		& 4.72e+3 \\
			0.01 & 4.24e-2	& 1.01e-1	& 32.52 & 10  & 2.15	& 7.50e+3	& 6 & 1.36		& 6.78e+3 \\
			0.005& 3.23e-2	& 1.01e-1	& 23.55	& 9   & 1.85	& 1.04e+4	& 5	& 1.17		& 9.91e+3 \\
			0.002& 2.58e-2	& 1.02e-1	& 14.57	& 7   & 1.57	& 1.60e+4	& 4	& 9.62e-1	& 1.55e+4 \\
			0.005& 2.37e-2	& 1.01e-1	& 10.07	& 6   & 1.42	& 2.32e+4	& 4	& 9.60e-1	& 2.21e+4 \\
			\hline
		\end{tabular}
		\vspace*{2mm}
		\caption{Test 3, refining the time step in the linear solid model, coupling without mesh intersection. The simulations are run on the Shaheen cluster. Same format as in Tables \ref{lin_tab_opti_inters} and \ref{lin_tab_opti_nointers}.}
		\label{lin_tab_time_nointers}
	\end{center}
\end{table}

\subsubsection{Volume loss}
In Table \ref{lin_tab_fluidref_big} and \ref{lin_tab_fluidref_small} we report the percentage of volume loss varying the fluid mesh and keeping fixed the solid one. In particular, in the first table, the solid mesh is made up of $384 \times 192$ elements, while in the second table is made up of $192 \times 96$ elements. In both cases, the fluid mesh varies from $512\times 512$ elements to $64 \times 64$ elements.

In both tables, we can see that the variation is higher when the coupling matrix is assembled computing the intersection between the two meshes; anyway, we can highlight two particular aspects.

In Table \ref{lin_tab_fluidref_big}, the assembly with intersection combined with the block-diag preconditioner does not converge when the fluid mesh is coarse ($64\times 64$ and $128\times 128$) since the variation of area is $99.94\%$ and $72.14\%$.

On the other hand, in Table \ref{lin_tab_fluidref_small}, there are not degenerate situations and conversely we can notice an increase in the volume loss when the coupling matrix is assembled without mesh intersection and the fluid mesh is significantly finer than the solid one ($512\times 512$ elements for the fluid mesh).

\begin{table}
	\begin{center}
		\begin{tabular}{r|r|r|r|r}
			\hline
			\multicolumn{5}{c}{{\bf Linear solid model -- Volume loss (\%)}} \\
			\hline
			\multicolumn{5}{c}{solid mesh $384\times192$, procs = 64, T = 2, $\dt$ = 0.01} \\
			\hline
			fluid mesh & \multicolumn{2}{c|}{\textit{with mesh intersection}} & \multicolumn{2}{c}{\textit{without mesh intersection}}\\
			& block-diag & block-tri & block-diag & block-tri\\
			\hline
			$512\times512$ & 2.14e-1 & 2.14e-1 & 4.15e-2 & 4.15e-2\\
			$256\times256$ & 2.33e-1 & 2.34e-1 & 4.24e-2 & 4.24e-2\\
			$128\times128$ & 72.15	 & 2.24e-1 & 6.69e-2 & 6.69e-2\\
			$64\times64$   & 99.94	 & 2.51e-1 & 6.59e-2 & 6.59e-2\\
			\hline
		\end{tabular}
	\end{center}
	\vspace*{2mm}
	\caption{Test 4, refining the fluid mesh keeping fixed the solid one. Loss of the structure volume in percentage.}
	\label{lin_tab_fluidref_big}
\end{table}

\begin{table}
	\begin{center}
		\begin{tabular}{r|r|r|r|r}
			\hline
			\multicolumn{5}{c}{{\bf Linear solid model -- Volume loss (\%)}} \\
			\hline
			\multicolumn{5}{c}{solid mesh $192\times96$, procs = 64, T = 2, $\dt$ = 0.01} \\
			\hline
			fluid mesh & \multicolumn{2}{c|}{\textit{with mesh intersection}} & \multicolumn{2}{c}{\textit{without mesh intersection}}\\
			& block-diag & block-tri & block-diag & block-tri\\
			\hline
			$512\times512$ & 5.48e-1 & 5.48e-1 & 2.66e-1 & 2.67e-1\\
			$256\times256$ & 4.60e-1 & 4.60e-1 & 5.76e-2 & 5.76e-2\\
			$128\times128$ & 5.06e-1 & 5.06e-1 & 6.89e-2 & 6.89e-2\\
			$64\times64$   & 4.95e-1 & 4.95e-1 & 6.61e-2 & 6.62e-2\\
			\hline
		\end{tabular}
	\end{center}
	\vspace*{2mm}
	\caption{Test 4, refining the fluid mesh keeping fixed the solid one. Loss of the structure volume in percentage.}
	\label{lin_tab_fluidref_small}
\end{table}

\subsection{Nonlinear solid model}
We now study the performance of the proposed solver in the case of a nonlinear model for the solid.

We assume that the energy density of the solid is given by
\begin{equation}
	W(\F) = \frac{\gamma}{2\eta} \exp\big(\eta \trace(\F^\top\F) - 2\big),
\end{equation}
where $\trace(\F^\top\F)$ is the trace of the tensor $\F^\top\F$ and $\gamma$ and $\eta$ are two constant parameters. $W$ is the exponential strain energy function of an isotropic hyperelastic material.

For our simulation, we consider a bar immersed in a fluid. At resting configuration, the structure is the rectangle ${\Omega_0^s = \B = [0,0.4]\times [0.45,0.55]}$, while the fluid domain $\Omega$ is the unit square. 

We impose homogeneous Dirichlet boundary conditions to $\u$ on $\partial\Omega$ and we set the following initial conditions
\begin{equation}
	\u(\x,0) = 0, \qquad \X(\s,0) = \s.
\end{equation}
The motion of the structure is generated by a force that pulls it down. This force is applied during the time interval $[0,1]$ at the midpoint of the right edge of the solid body.

In particular, we consider again the case where solid and fluid have same density and same viscosity, precisely, we set
\begin{equation}
	\rho_f=\rho_s=1 \quad\text{and}\quad \nu_f = \nu_s = 0.2.
\end{equation}
Moreover, we choose $\gamma = 1.333$ and $\eta=9.242$.

\begin{figure}
	\begin{center}
		\includegraphics[trim = 30 10 20 0,clip, width=4cm]{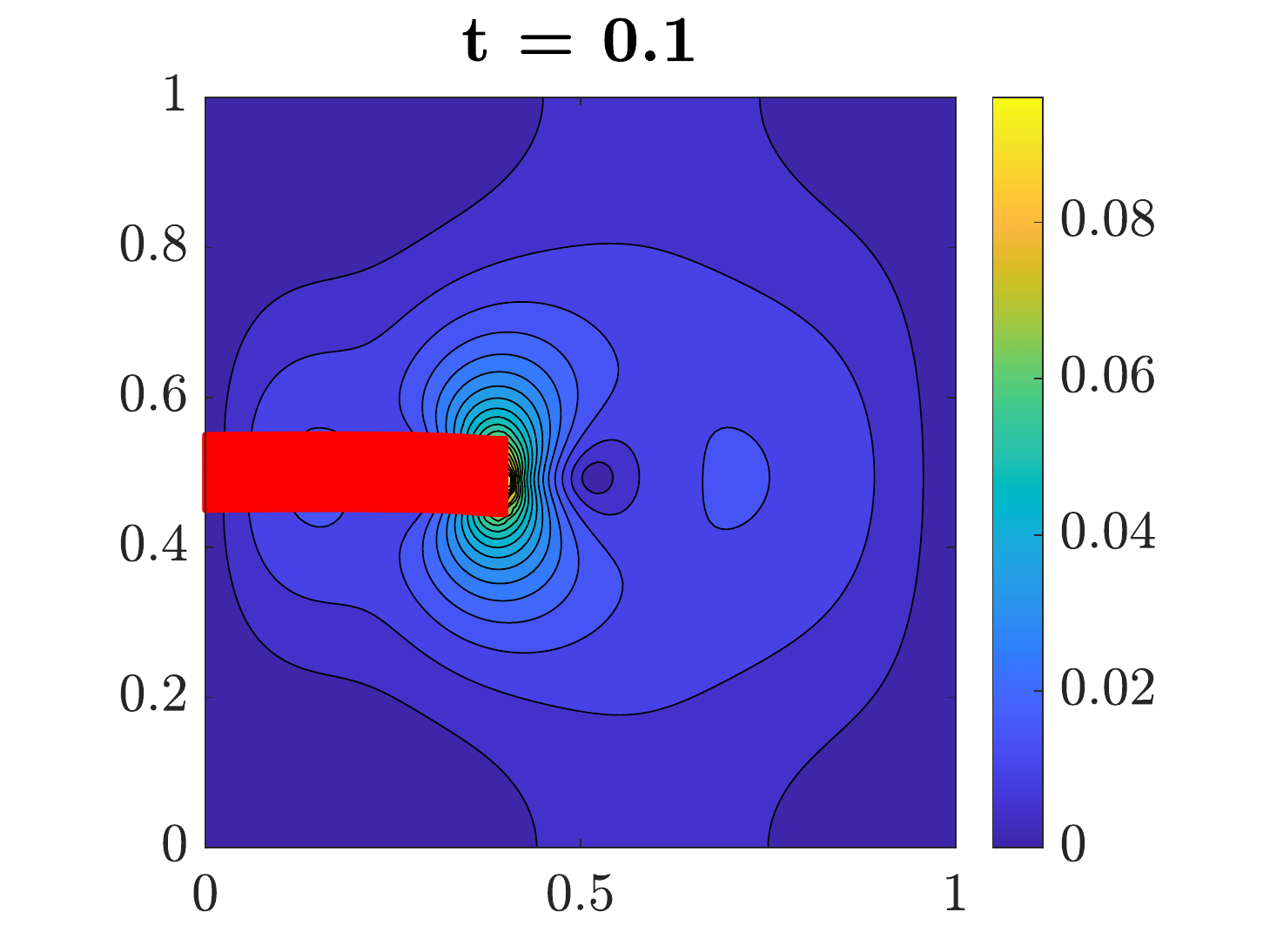}
		\includegraphics[trim = 30 10 20 0,clip, width=4cm]{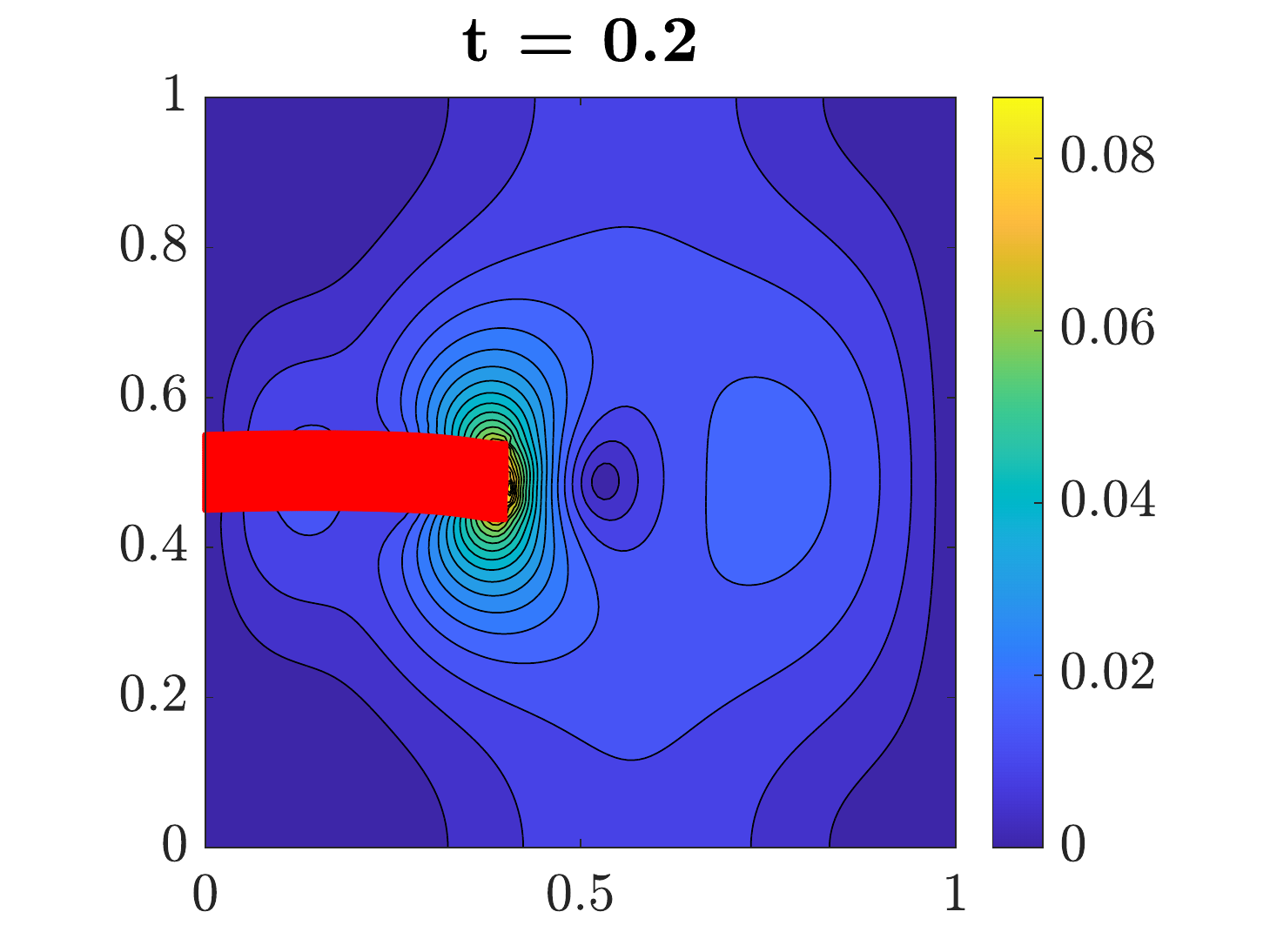}
		\includegraphics[trim = 30 10 20 0,clip, width=4cm]{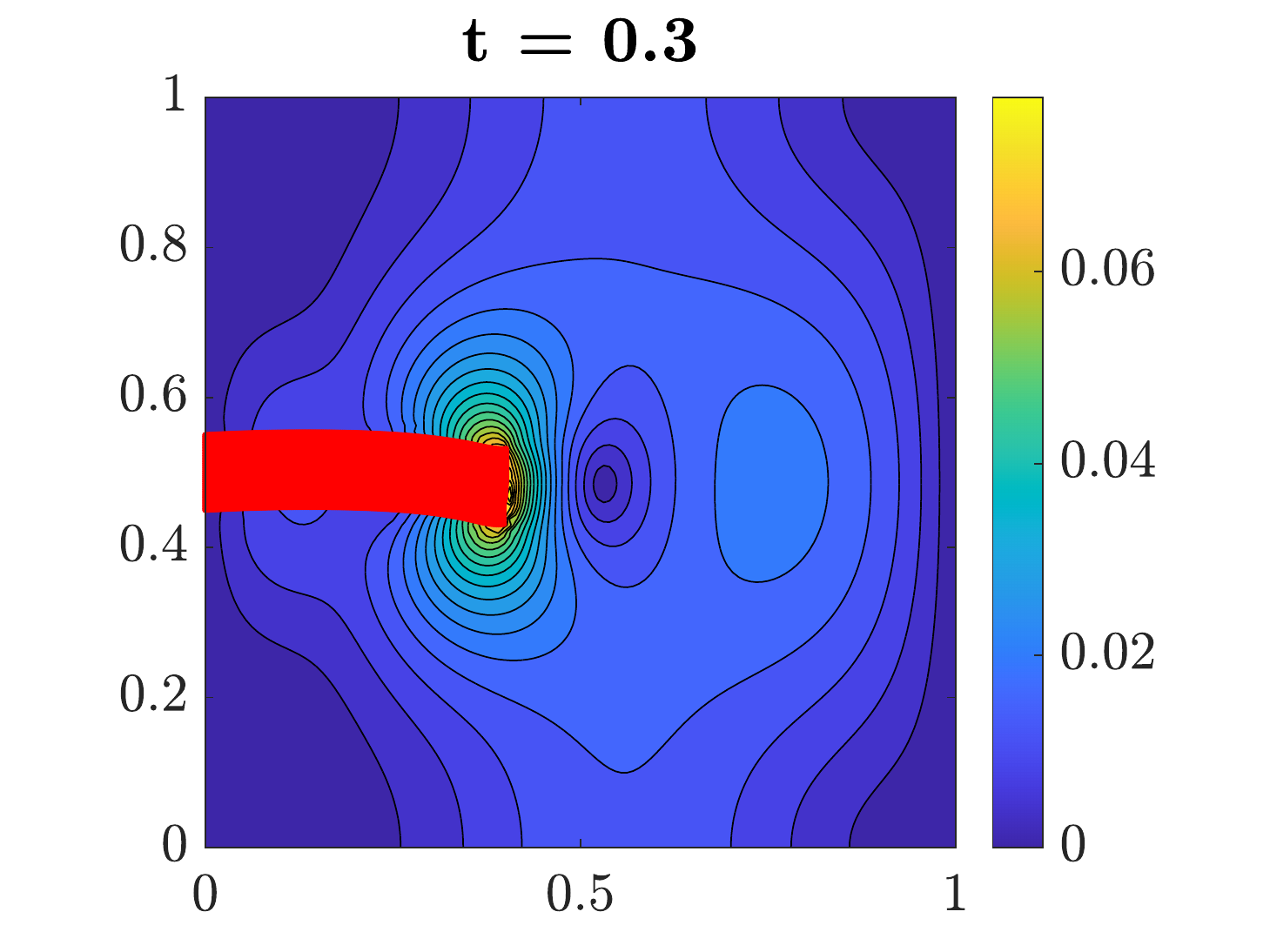}
		\includegraphics[trim = 30 10 20 0,clip, width=4cm]{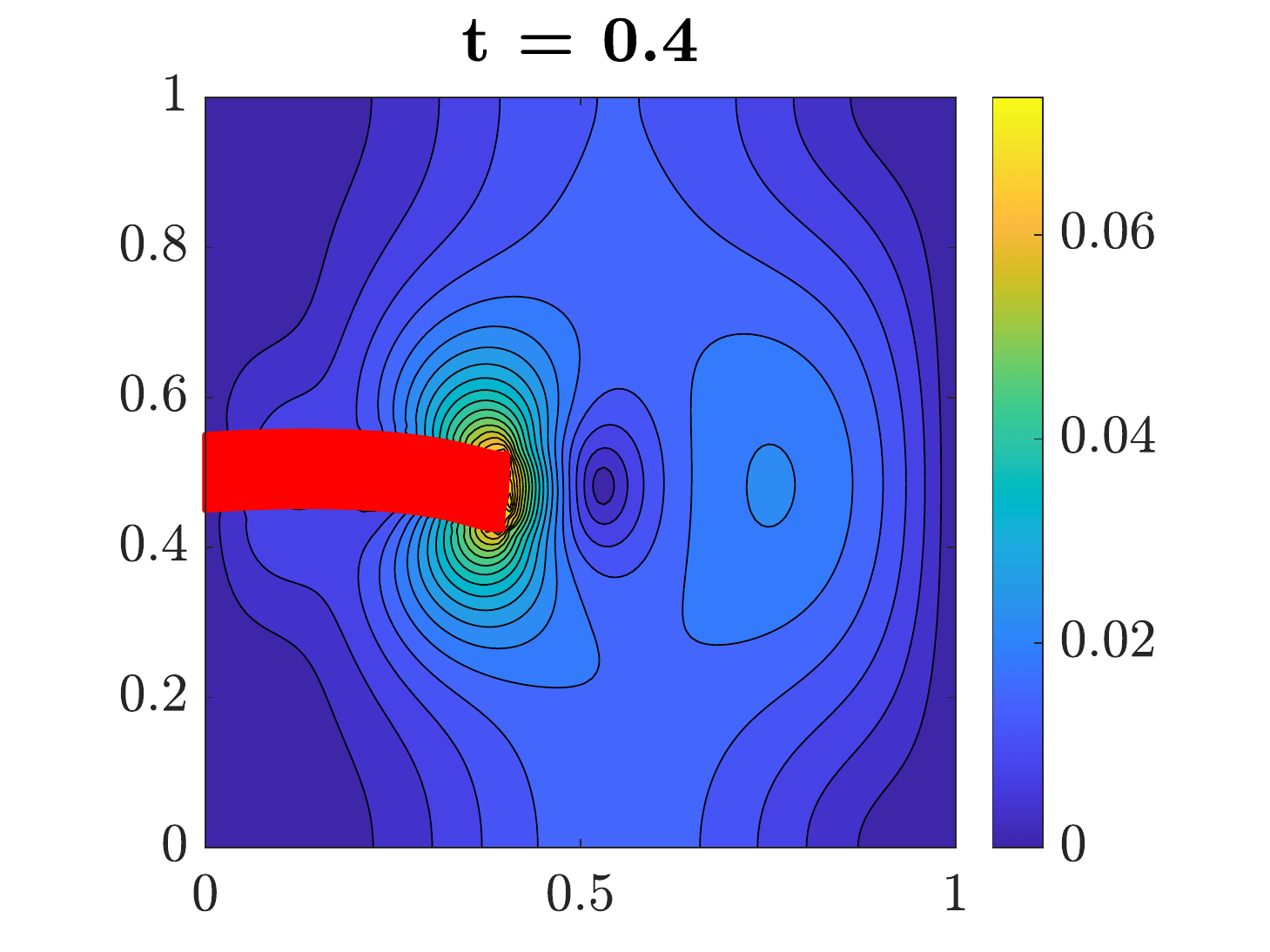}\\
		\medskip
		\includegraphics[trim = 30 10 20 0,clip, width=4cm]{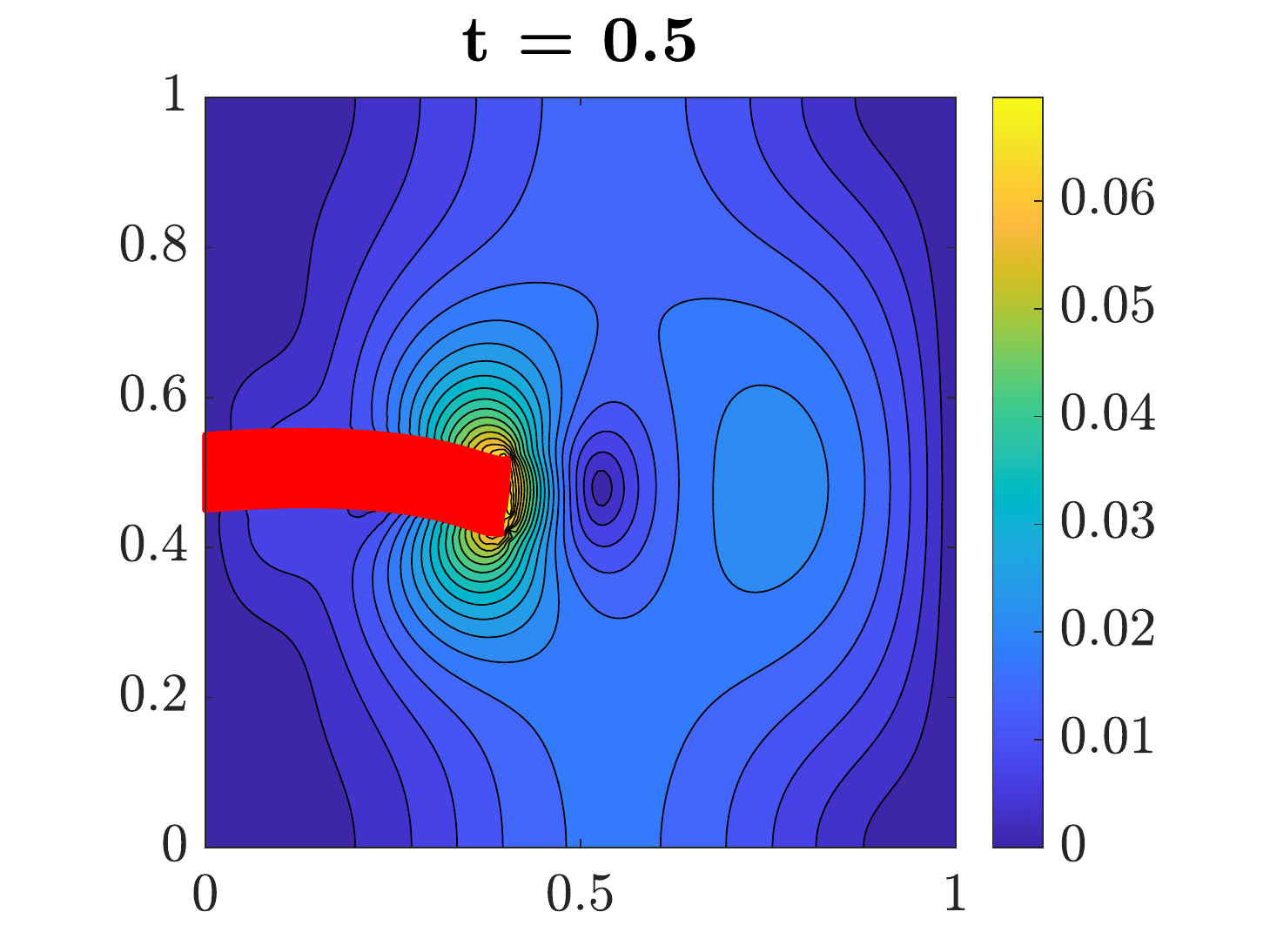}
		\includegraphics[trim = 30 10 20 0,clip, width=4cm]{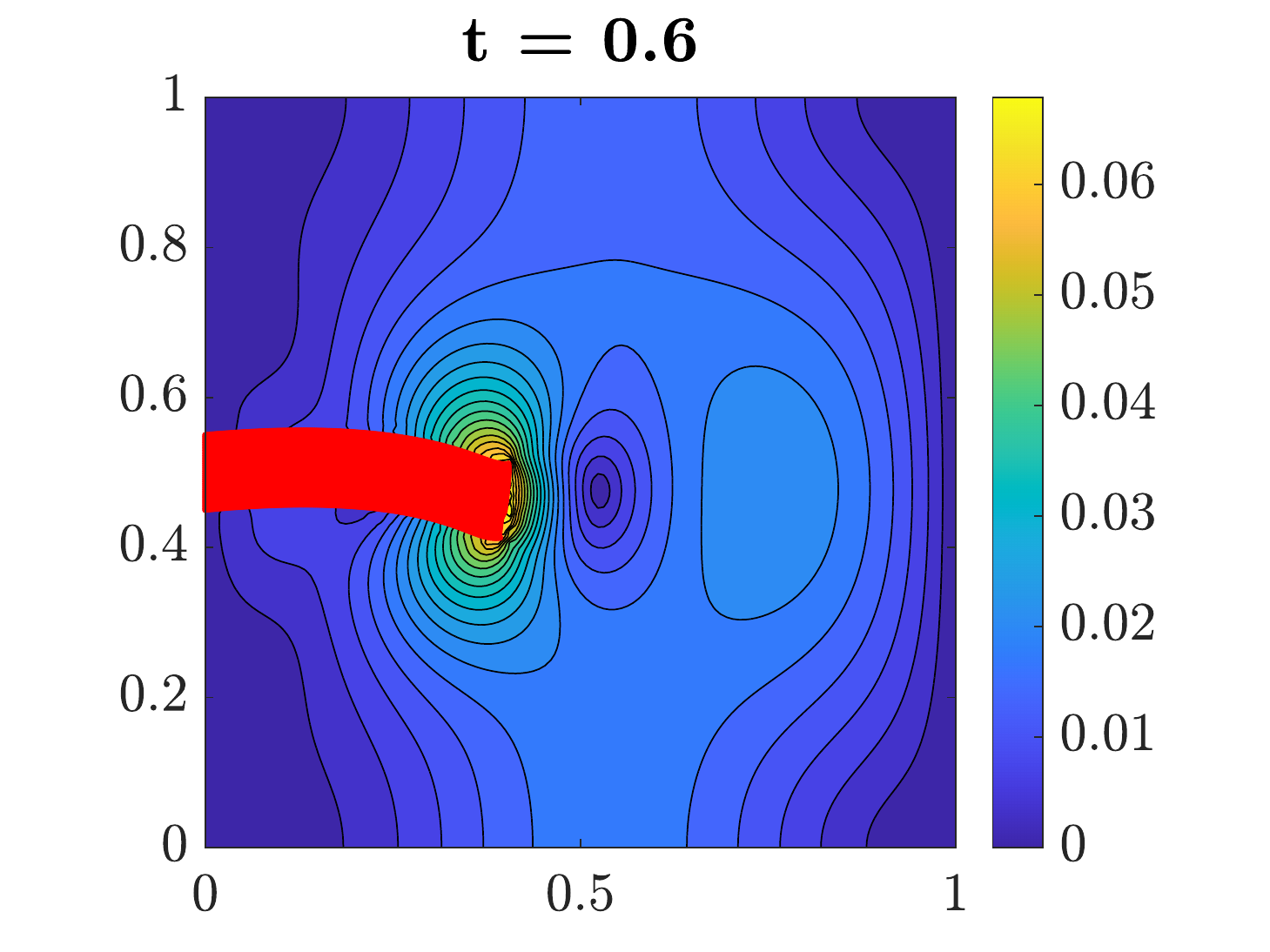}
		\includegraphics[trim = 30 10 20 0,clip, width=4cm]{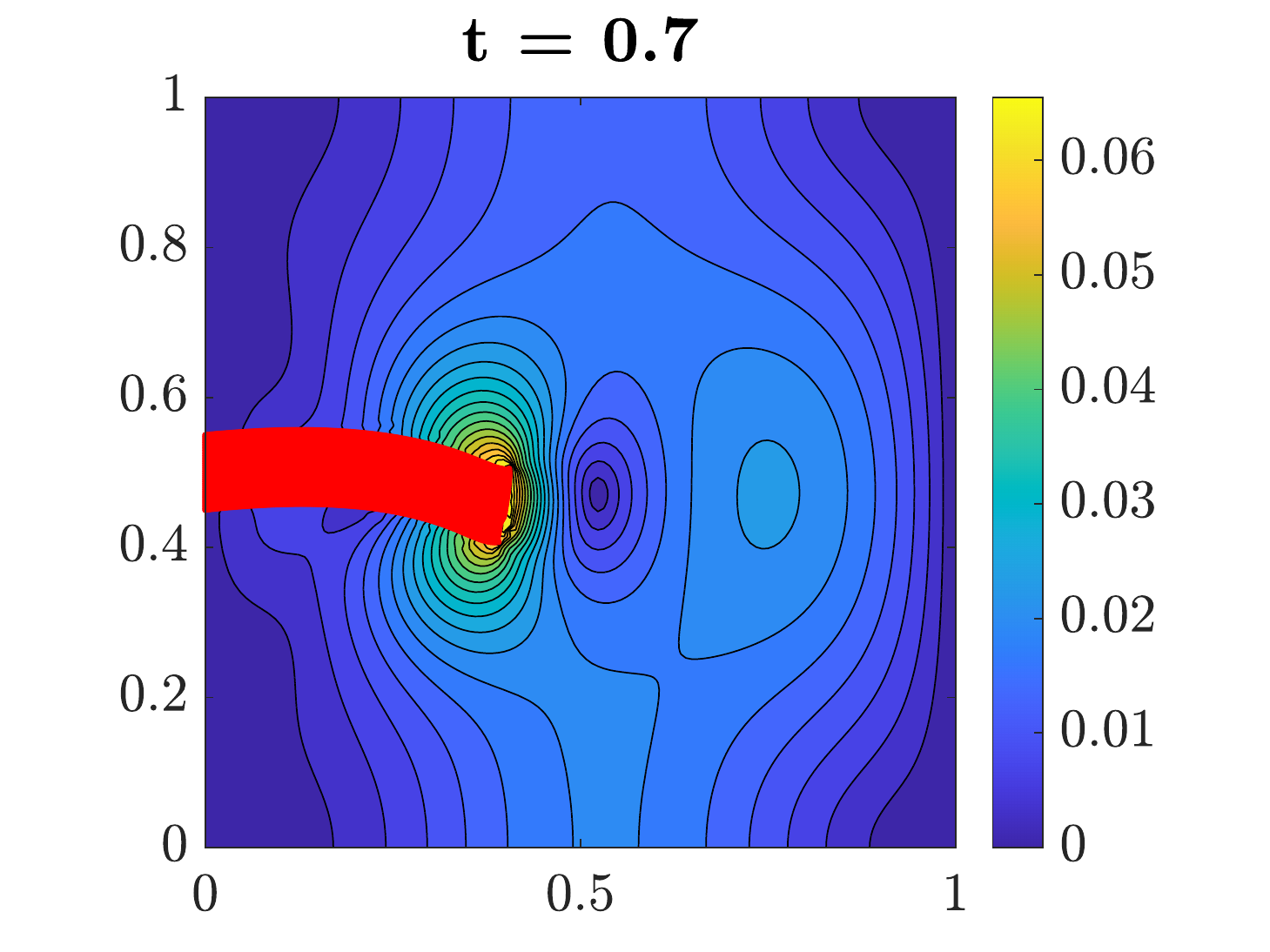}
		\includegraphics[trim = 30 10 20 0,clip, width=4cm]{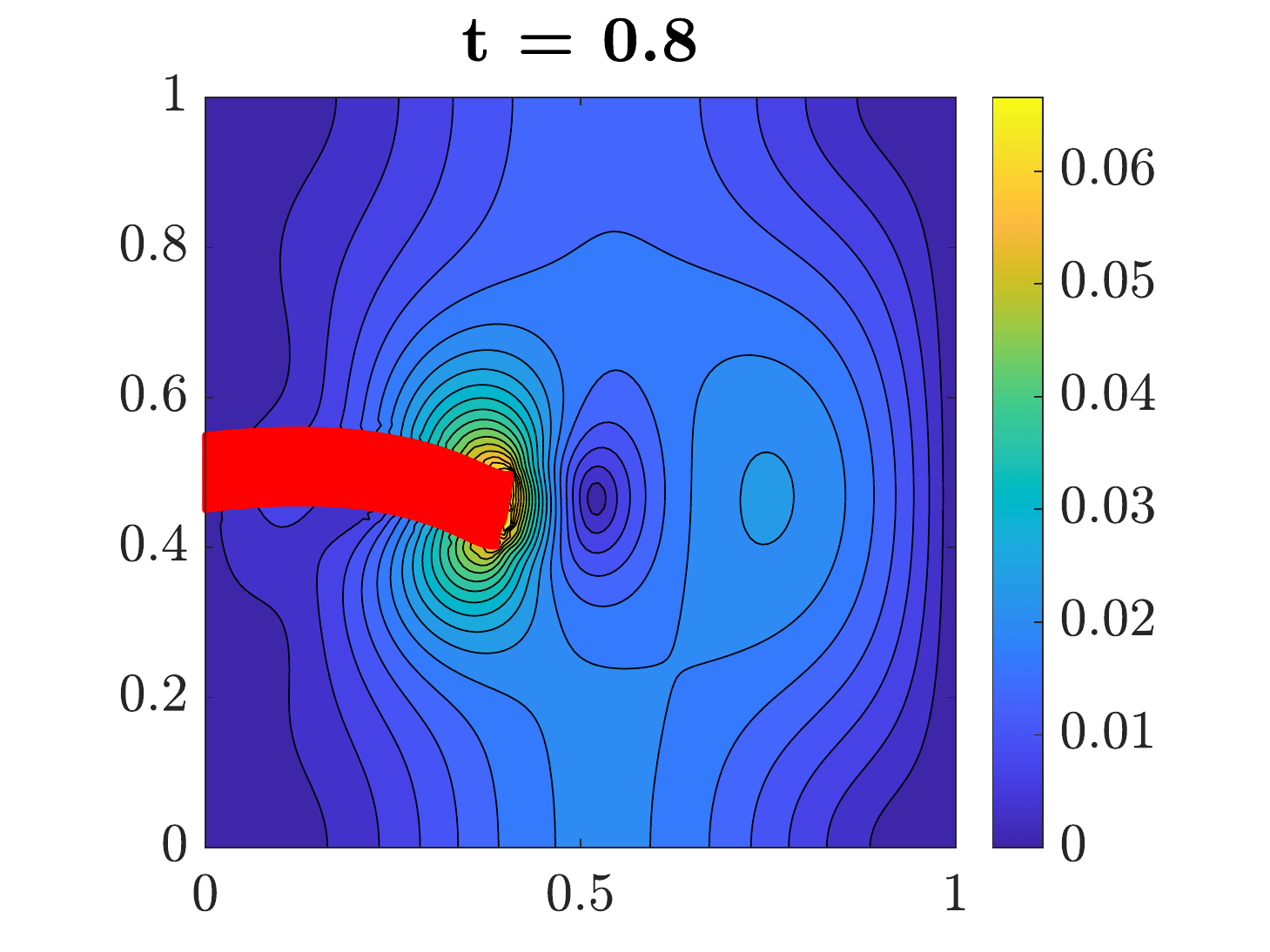}\\
		\medskip
		\includegraphics[trim = 30 10 20 0,clip, width=4cm]{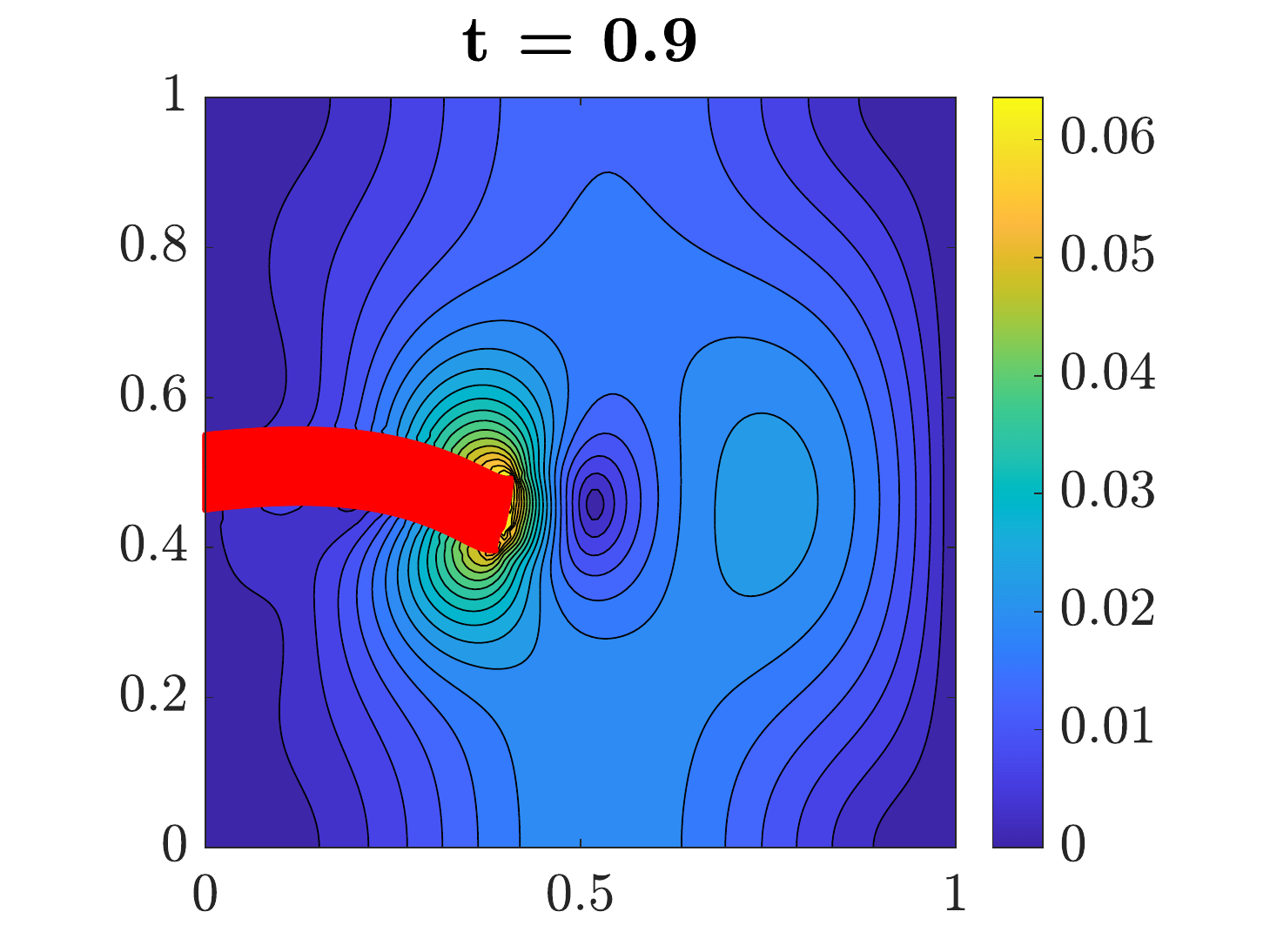}
		\includegraphics[trim = 30 10 20 0,clip, width=4cm]{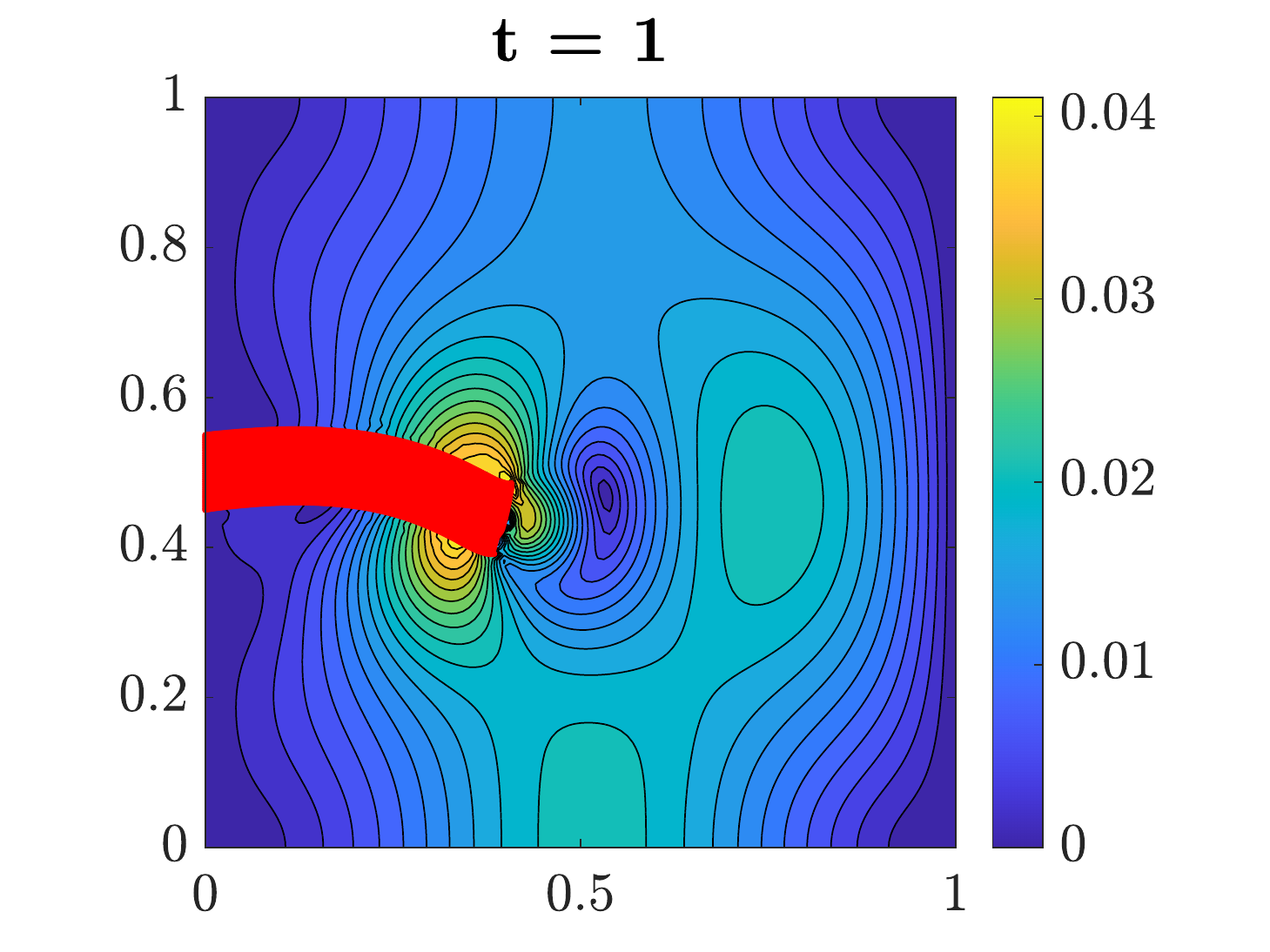}
		\includegraphics[trim = 30 10 20 0,clip, width=4cm]{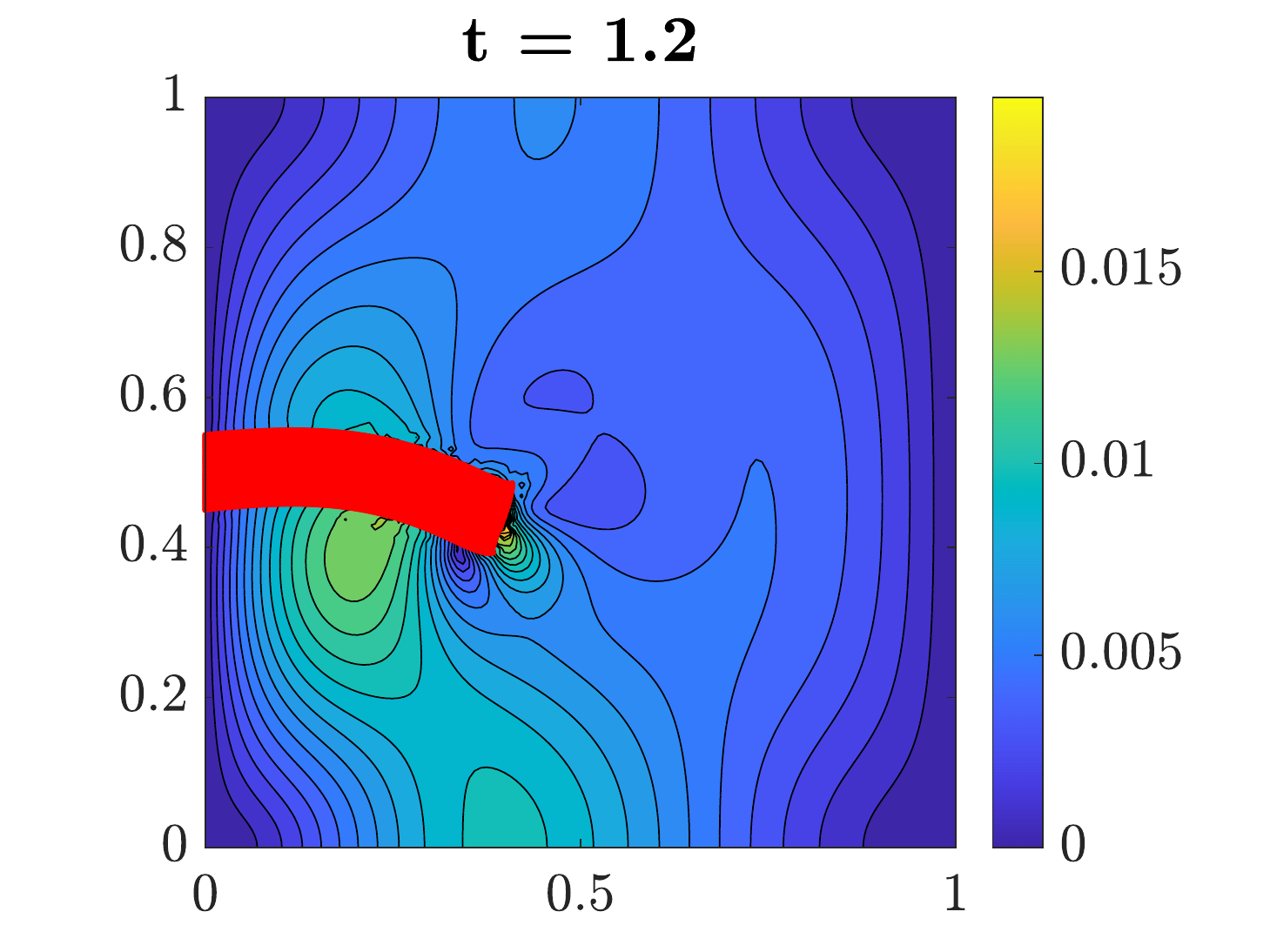}
		\includegraphics[trim = 30 10 20 0,clip, width=4cm]{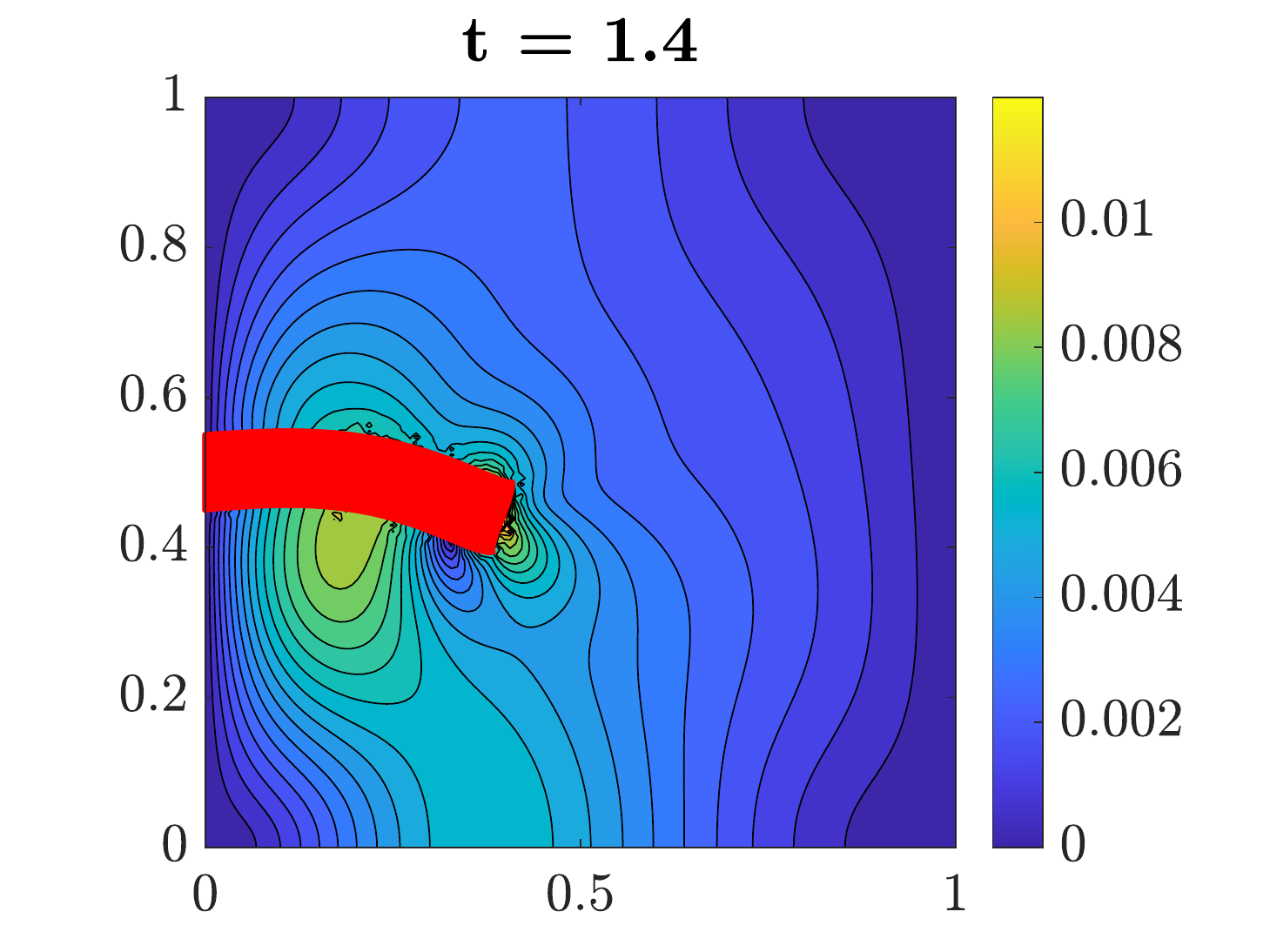}\\
		\medskip
		\includegraphics[trim = 30 10 20 0,clip, width=4cm]{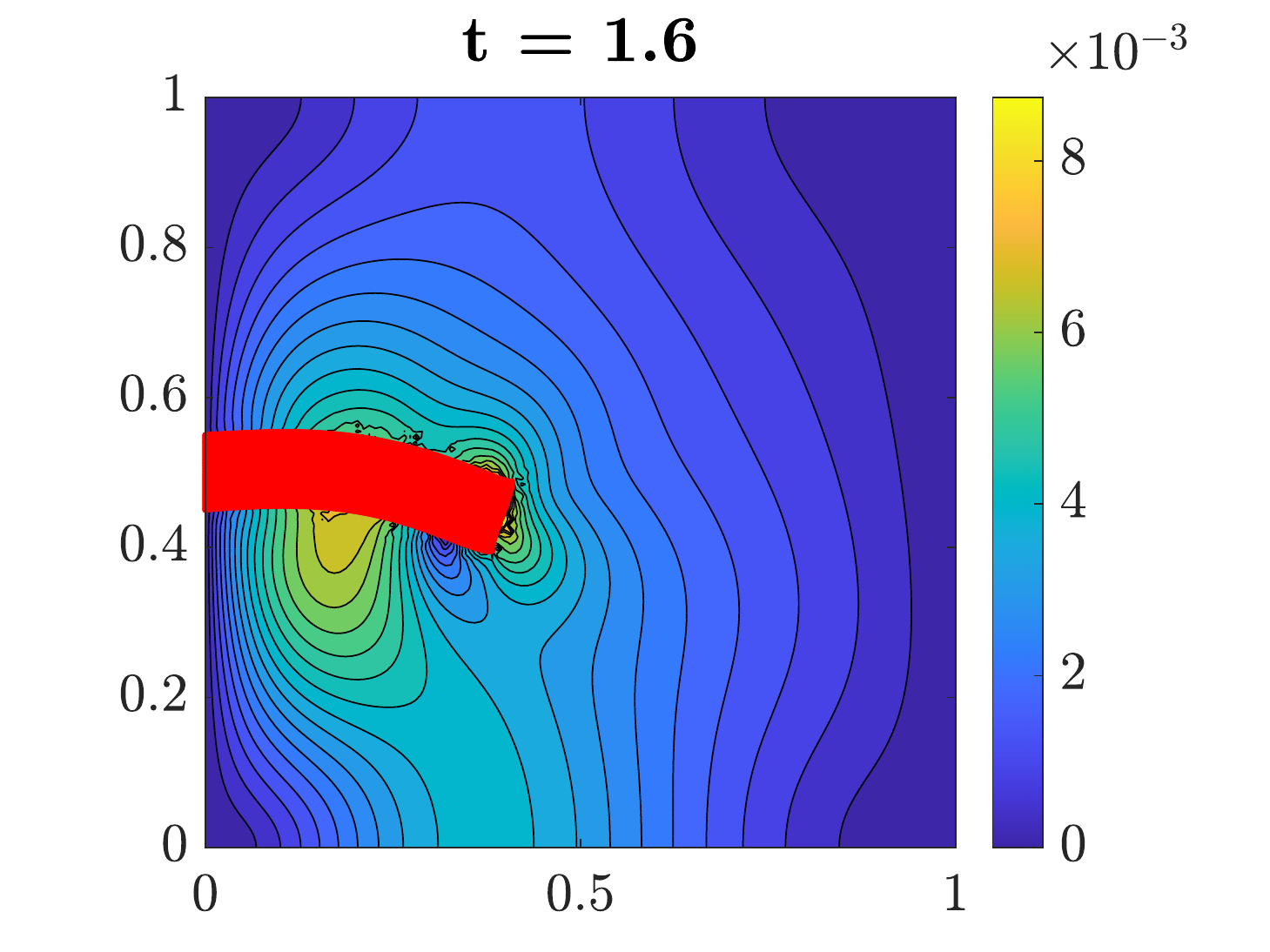}
		\includegraphics[trim = 30 10 20 0,clip, width=4cm]{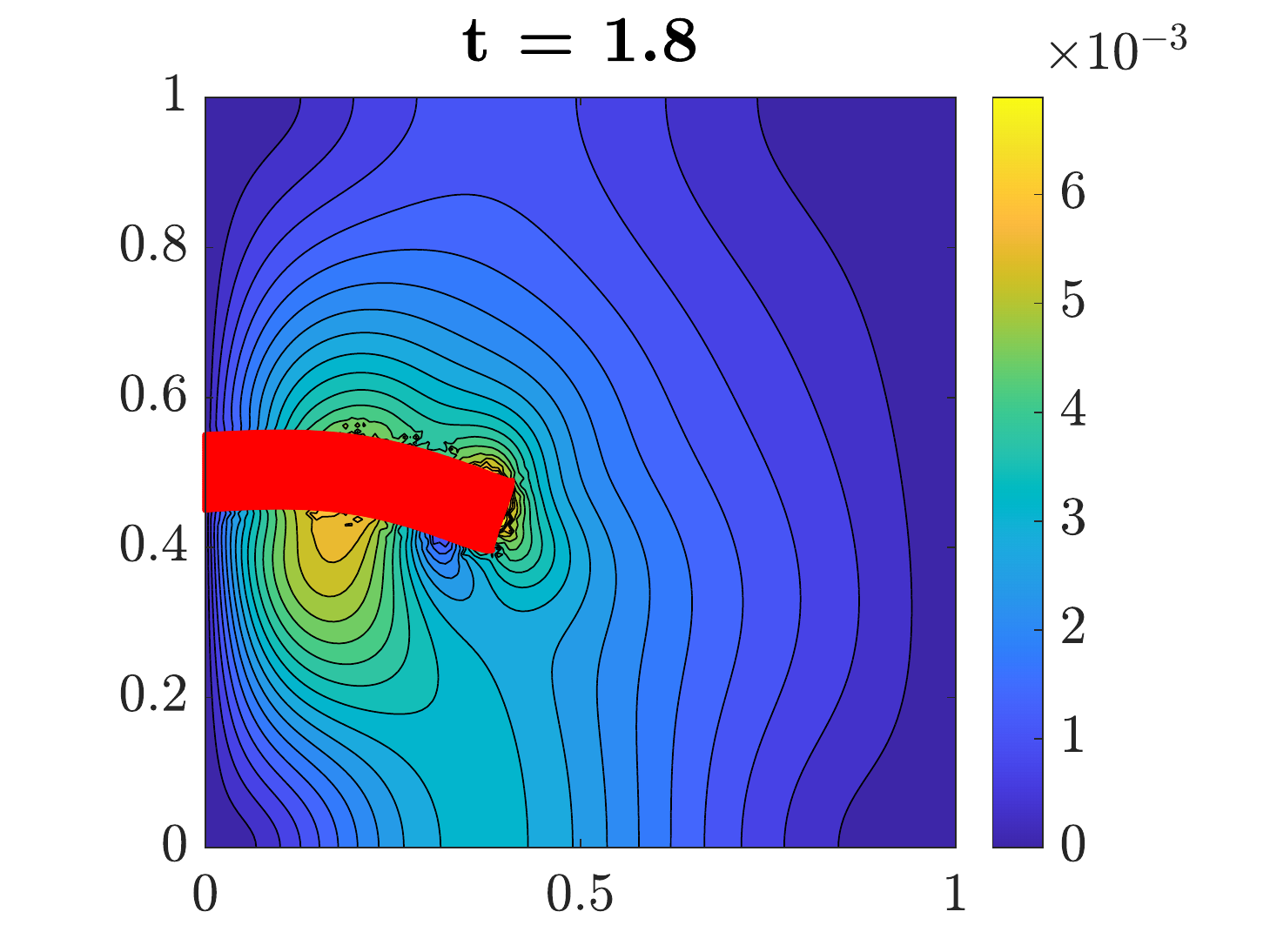}
		\includegraphics[trim = 30 10 20 0,clip, width=4cm]{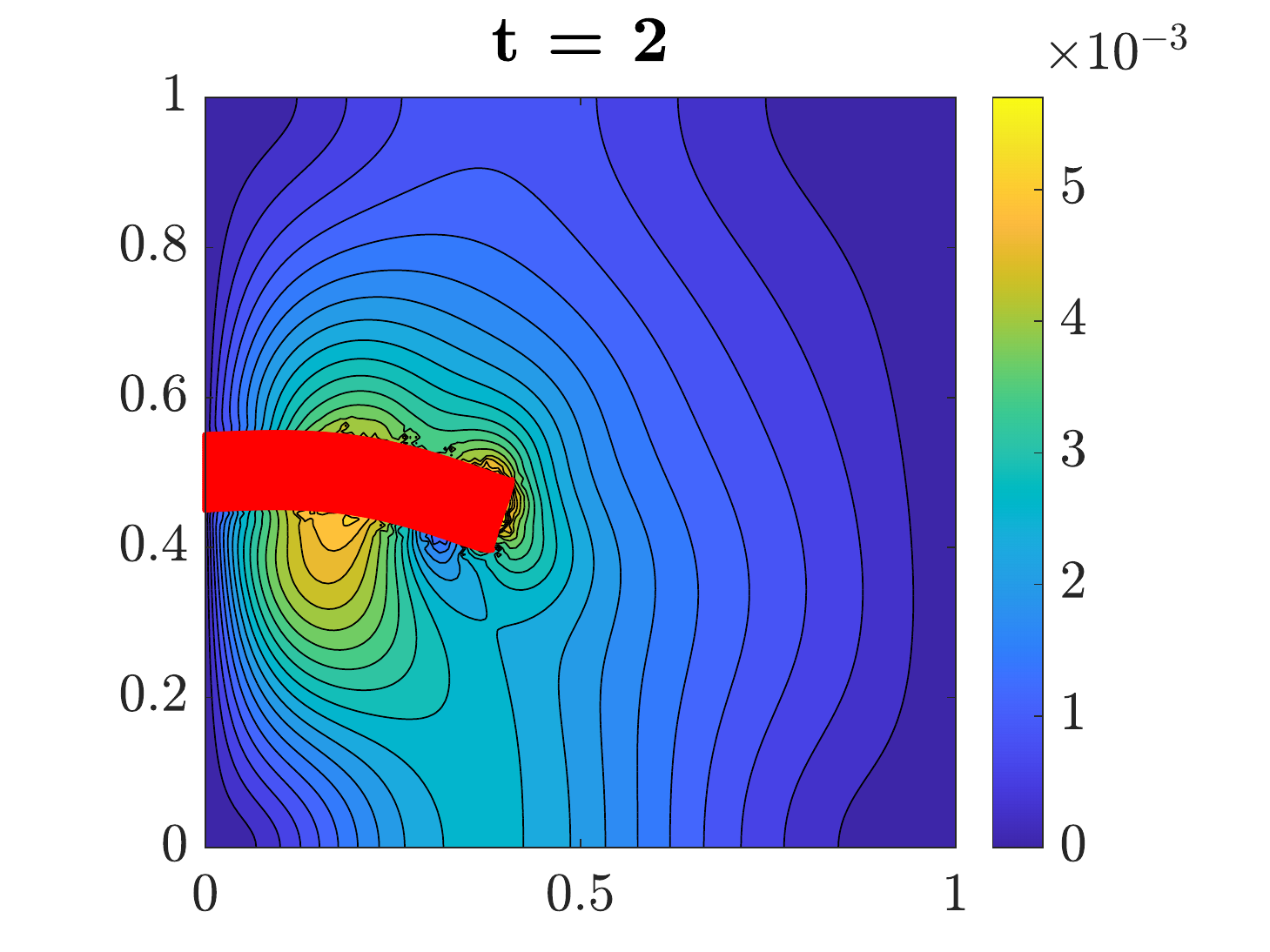}
		\includegraphics[trim = 30 10 20 0,clip, width=4cm]{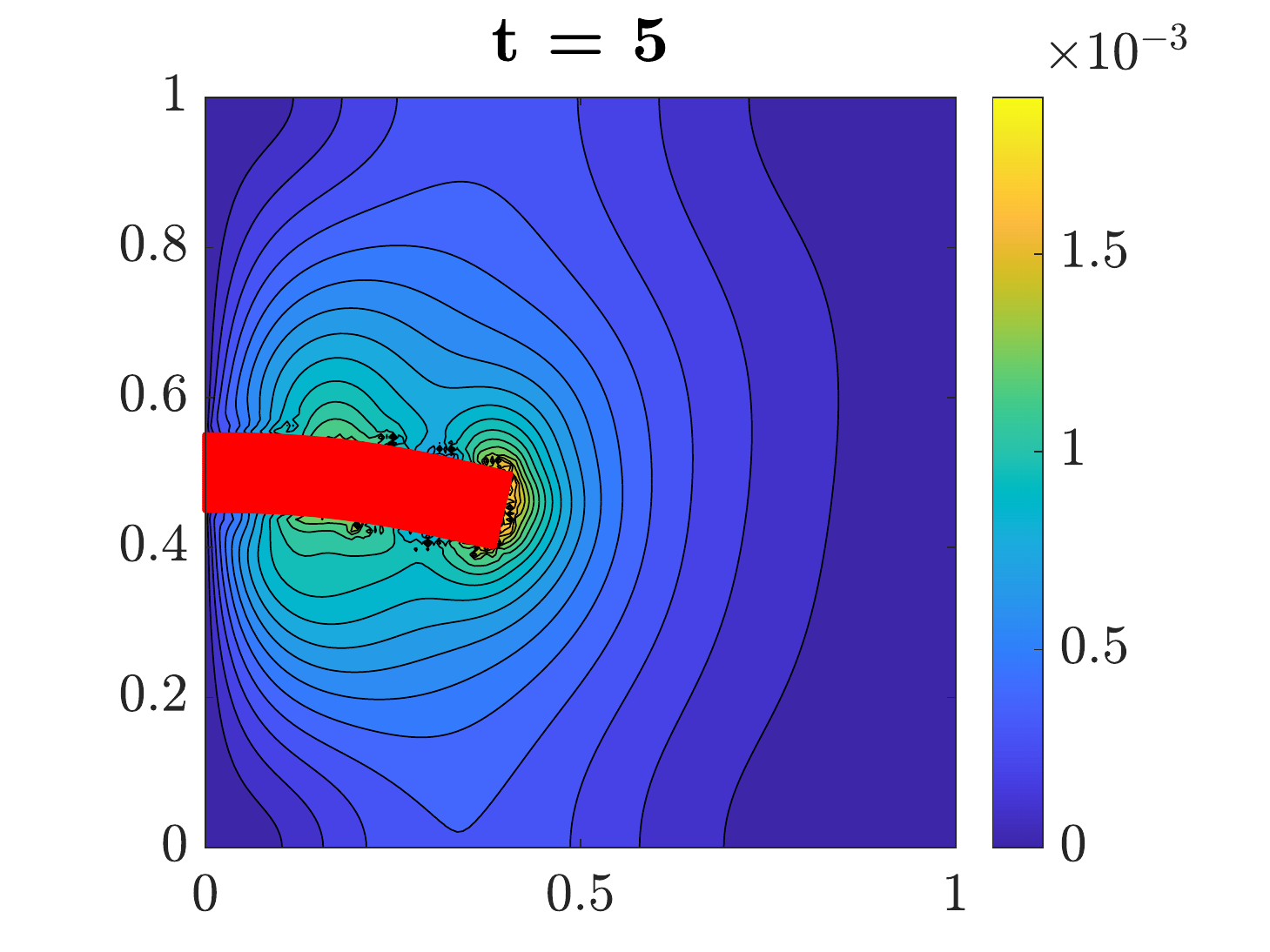}\\
		\caption{\rv Simulation of the bar with nonlinear constitutive law: some snapshots. The structure position is represented in red, while streamlines and color bars refer to the velocity. In the time interval $[0,1]$, the structure is pulled down by a force applied at the middle point of the right edge. Then, when the structure is released, internal forces bring it back to its resting configuration.\vr}%{Nonlinear solid model: snapshots of structure evolution.}
		\label{fig:bar_evolution}
	\end{center}
\end{figure}

Also in this case we are dealing with incompressible materials, therefore, during the simulation, the area of the structure should remain theoretically constant. We study the evolution of the system in the time interval $[0,2]$. Some snapshots are reported in Figure \ref{fig:bar_evolution}.

\begin{table}
\begin{center}
\begin{tabular}{r|r|r|r|r|r|r|r|r|r|r|r}
\hline
\multicolumn{12}{c}{{\bf Nonlinear solid model -- Mesh refinement test}} \\
\multicolumn{12}{c}{{\textit{Coupling with mesh intersection}}} \\
\hline
\multicolumn{12}{c}{procs = 32, T = 2, $\dt$ = 0.01} \\
\hline
dofs    & vol. loss (\%)        & $T_{ass}(s)$  & $T_{coup}(s)$ & \multicolumn{4}{c|}{block-diag}               & \multicolumn{4}{c}{block-tri} \\
	&                       &               &               & nit   & its   & $T_{sol}(s)$  & $T_{tot}(s)$  & nit   & its   & $T_{sol}(s)$  & $T_{tot}(s)$ \\
\hline
21222   & 1.70	& 7.24e-3  	 & 6.08e-2 	 & 2 	& 245		& 13.28		& 2.67e+3	& 2		& 21   & 1.56	& 3.24e+2 \\
83398   & 1.47	& 2.87e-2  	 & 8.41e-1 	 & 2 	& 269		& 47.78		& 9.72e+3	& 2		& 23   & 5.57	& 1.28e+3\\
186534  & 1.47	& 6.34e-2  	 & 4.10 	 & 3 	& 388		& 1.55e+2	& 3.19e+4	& 3		& 26   & 14.38	& 3.71e+3\\
330630  & 1.47	& 1.13e-1  	 & 9.05 	 & 3 	& 406		& 2.59e+2	& 5.37e+4	& 3		& 27   & 23.95	& 6.62e+3\\
515686	& -		& -  	 	 & - 	 	 & - 	& -			& -			& -			& -		& -    &-		&- \\
741702  & -		& -  	 	 & - 	 	 & -  	& -			&  - 		& -  		&	-	& -    &-	 	&- \\
\hline
\end{tabular}
\vspace*{2mm}
\caption{Test 5, refining the mesh in the nonlinear solid model, coupling with mesh intersection. The simulations are run on the Shaheen cluster. procs = number of processors; dofs = degrees of freedom; vol. loss = loss of structure volume in percentage; $T_{ass}$ = CPU time to assemble the stiffness and mass matrices; $T_{coup}$ = CPU time to assemble the coupling term; nit = Newton iterations; its = GMRES iterations to solve the Jacobian system; $T_{sol}$ = CPU time to solve the Jacobian system; $T_{tot}$ = total simulation CPU time. The quantities $T_{coup}$ and nit are averaged over the time steps, whereas the quantities its and $T_{sol}$ are averaged over the Newton iterations and the time steps. All CPU times are reported in seconds.}
\label{nonlin_tab_opti_inters}
\end{center}
\end{table}

\begin{table}
\begin{center}
\begin{tabular}{r|r|r|r|r|r|r|r|r|r|r|r}
\hline
\multicolumn{12}{c}{{\bf Nonlinear solid model -- Mesh refinement test}} \\
\multicolumn{12}{c}{{\textit{Coupling without mesh intersection}}} \\
\hline
\multicolumn{12}{c}{procs = 32, T = 2, $\dt$ = 0.01} \\
\hline
dofs    & vol. loss (\%)        & $T_{ass}(s)$  & $T_{coup}(s)$ & \multicolumn{4}{c|}{block-diag}               & \multicolumn{4}{c}{block-tri} \\
        &                       &               &               & nit   & its   & $T_{sol}(s)$  & $T_{tot}(s)$  & nit   & its   & $T_{sol}(s)$  & $T_{tot}(s)$ \\
\hline
21222   & 3.50e-1               & 4.90e-3       & 1.61e-2       & 3     & 245   & 7.26          & 4.49e+3       & 3     & 11    & 5.47e-1       & 331.91 \\
83398   & 2.07e-1               & 2.04e-2       & 6.89e-2       & 3     & 498   & 43.32         & 2.70e+4       & 3     & 14    & 1.99          & 1.24e+3 \\
186534  & 1.36e-1               & 4.56e-2       & 2.15e-1       & 3     & 1572  & 308.33        & 1.85e+5       & 3     & 16    & 4.92          & 3.18e+3 \\
330630  & 9.97e-2               & 8.20e-2       & 6.42e-1       & 3     & 3001  & 924.25        & 5.55e+5       & 3     & 18    & 8.63          & 5.95e+3 \\
515686  & 7.16e-2               & 1.26e-1       & 1.85          & -     & -     & -             & -             & 3     & 25    & 17.71         & 1.31e+4 \\
741702  & 5.34e-2               & 3.19e-1       & 5.45          & -     & -     & -             & -             & 3     & 30    & 28.59         & 2.19e+4 \\
\hline
\end{tabular}
\vspace*{2mm}
\caption{Test 5, refining the mesh in the nonlinear solid model, coupling without mesh intersection. The simulations are run on the Eos cluster. Same format as in Table \ref{nonlin_tab_opti_inters}.}
\label{nonlin_tab_opti_nointers}
\end{center}
\end{table}

\subsection{Mesh refinement test}

We study the behavior of the parallel FSI solver with respect to mesh refinement, by keeping fixed the time step size $\Delta t = 0.01$,
the number of processors $procs = 32$ and the final time $T = 2$.
The number of degrees of freedom (dofs) varies from 21222 to 741702.

Table \ref{nonlin_tab_opti_inters} reports the results in case of coupling with mesh intersection, run on the Shaheen cluster. First of all, we notice that the volume loss remains constant when the mesh is refined, indeed it stabilizes at $1.47\%$. With respect to the number of dofs, $T_{ass}$ exhibits a moderate increase, whereas the growth of $T_{coup}$ is superlinear. For both preconditioners, the number of Newton iterations remains bounded by 3 when we increase the dofs. However, in the two finest cases, the Newton's algorithm fails. \fc This might be due to the fact that we are using a large time step size with respect to the spatial meshes. Indeed in Section~\ref{sec:time_nonlin}, we show that with smaller time steps the method works\cf. In addition, we have that in the case of the block-diag preconditioner, the GMRES solver we use at each Newton step increases its number of iterations. As a consequence, $T_{sol}$ grows rapidly. On the other hand, the block-tri solver is more robust with respect to the mesh refinement, indeed the number of GMRES iterations grows moderately.

Table \ref{nonlin_tab_opti_nointers} reports the results in case of coupling without mesh intersection.
These simulations are run on the Eos cluster.
We first observe that the volume loss reduces when the mesh is refined and, for the fines grids, it is very small, always below $0.4\%$.
With respect to the number of dofs, $T_{ass}$ exhibits a moderate increase, whereas the growth of $T_{coup}$ is superlinear.
For both block-diag and block-tri solvers, the Newton iterations remain constant when the number of dofs increases.
However, the block-diag solver is not robust with respect to mesh refinement, because the GMRES iterations, 
needed to solve the Jacobian linear system at each Newton iteration, increase when the mesh is refined.
In the two finest cases, the solver fails \fc because it reaches the maximum number of allowed iterations\cf.
Thus, $T_{sol}$ exhibits a significant growth with a large number of dofs.
The block-tri solver, instead, is quite robust with respect to mesh refinement, because the GMRES iteration count exhibits a moderate growth when the mesh is refined.

\subsection{Strong scalability test}

Here we study the strong scalability of the parallel nonlinear FSI solver, by keeping fixed the fluid and solid meshes to $160\times 160$ and $480\times 120$ elements, respectively. \lg This choice produces \gl a total amount of 515686 dofs. The time parameters are fixed at $\Delta t = 0.002$ and $T = 2$ in the case with mesh intersection, and 
at $\Delta t = 0.01$ and $T = 2$ in the case without mesh intersection.
The number of processors increases from 4 to 64. The parallel speedup $S^p$ is computed with respect to the 4 processors run.
We consider only the block-tri preconditioner.

Table \ref{nonlin_scal_int} reports the results in case of coupling with mesh intersection. We first observe a good scalability of the
assembly phase of the coupling term ($T_{coup}$). The solver is scalable in terms of both Newton and GMRES iterations,
which remain bounded when the number of processors increases. However, the solution time $T_{sol}$ is not scalable.
Consequently, the global performance of the solver is impaired, with slow speedup values.
This is due to the poor performance of the direct solvers used to invert the diagonal blocks of the preconditioner.

Table \ref{nonlin_scal_noint} reports the results in case of coupling without mesh intersection.
We observe again a scalable behavior of the assembly phase of the coupling term and
of the Newton and GMRES iterations.
However, as in the case with mesh intersection, the non-scalable behavior of the solution time $T_{sol}$ impairs
the global performance of the solver when the number of processors increases.

\begin{figure}[h!]
	\vspace{5mm}
	\begin{center}
		%\includegraphics[width=8cm]{figures/fig_nl/fig_nonlin_int_1.png}
		%\includegraphics[width=8cm]{figures/fig_nl/fig_nonlin_int_2.png}\\
		%\medskip
		%\includegraphics[width=8cm]{figures/fig_nl/fig_nonlin_int_3.png}
		%\includegraphics[width=8cm]{figures/fig_nl/fig_nonlin_int_4.png}
		\begin{overpic}[width=7.95cm]{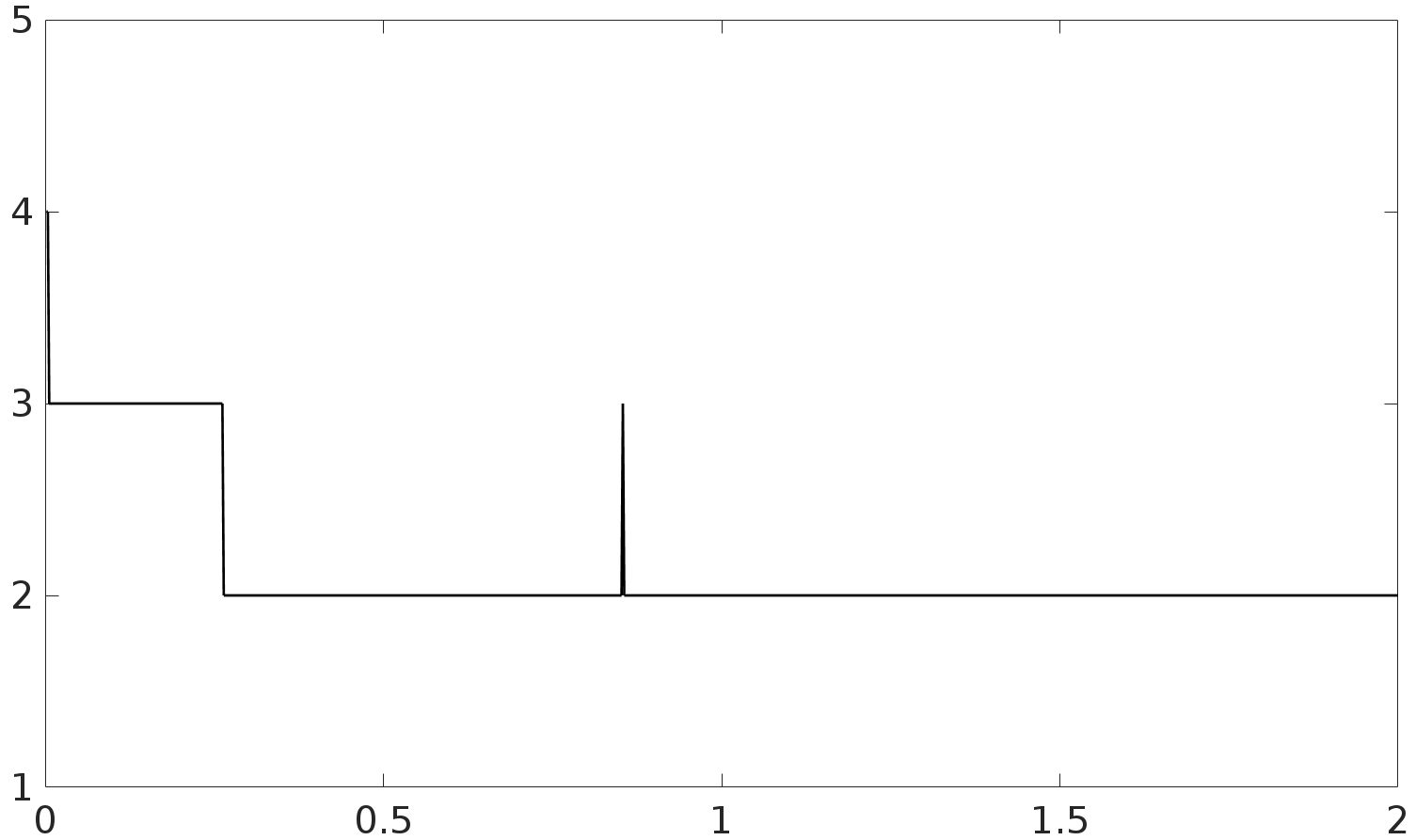}
			\put (47,-5) {\scriptsize time}
			\put (32.5,59) {\scriptsize \textbf{Nonlinear iterations}}
			\put (-5,25) {\begin{sideways}
					\scriptsize $nit$
			\end{sideways}}
		\end{overpic}
		\quad
		\begin{overpic}[width=8cm]{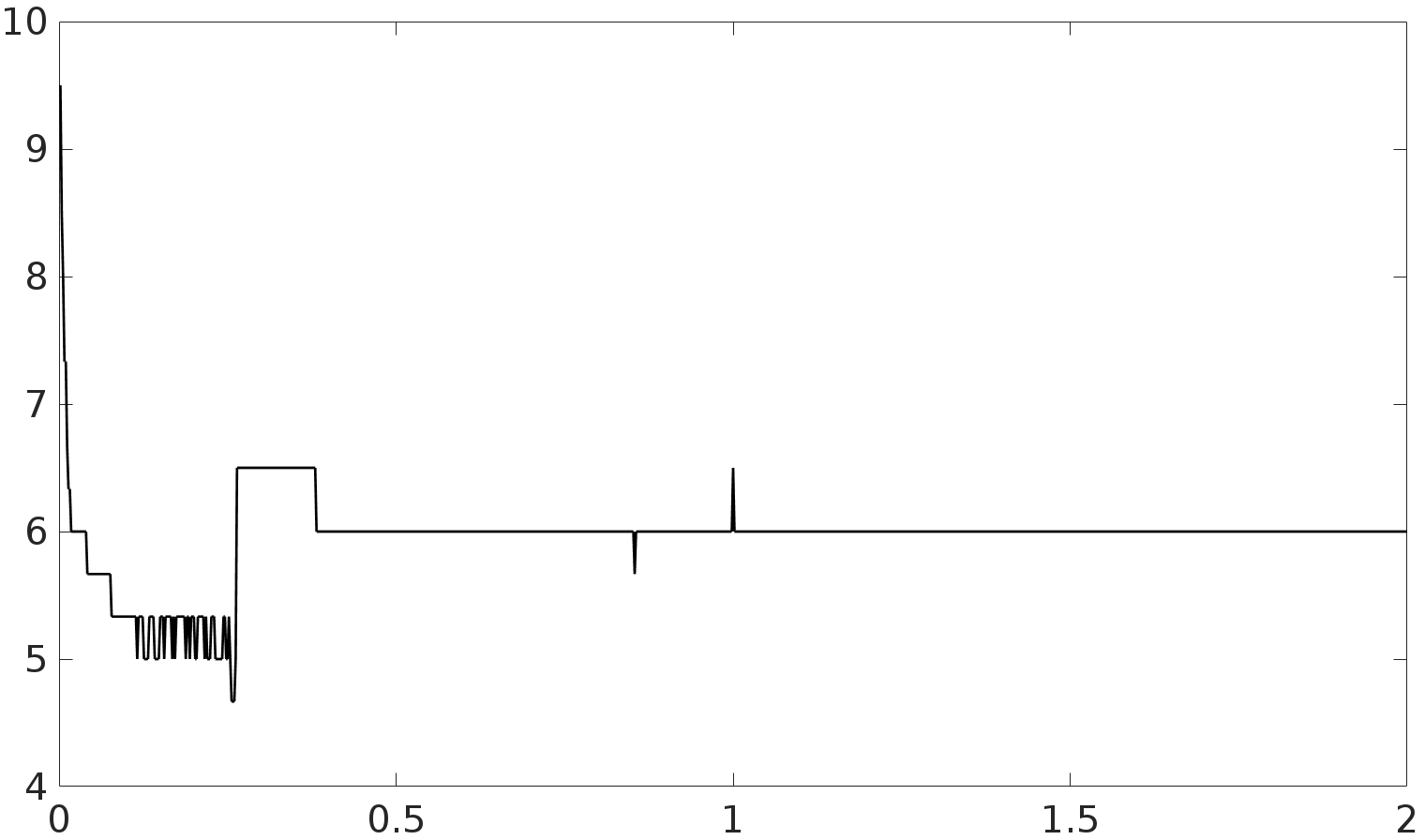}
			\put (48,-5) {\scriptsize time}
			\put (6.5,59) {\scriptsize \textbf{Average linear iterations per nonlinear iteration}}
			\put (-5,25) {\begin{sideways}
					\scriptsize $its$
			\end{sideways}}
		\end{overpic}\\
		\vspace{10mm}
		\begin{overpic}[width=8cm]{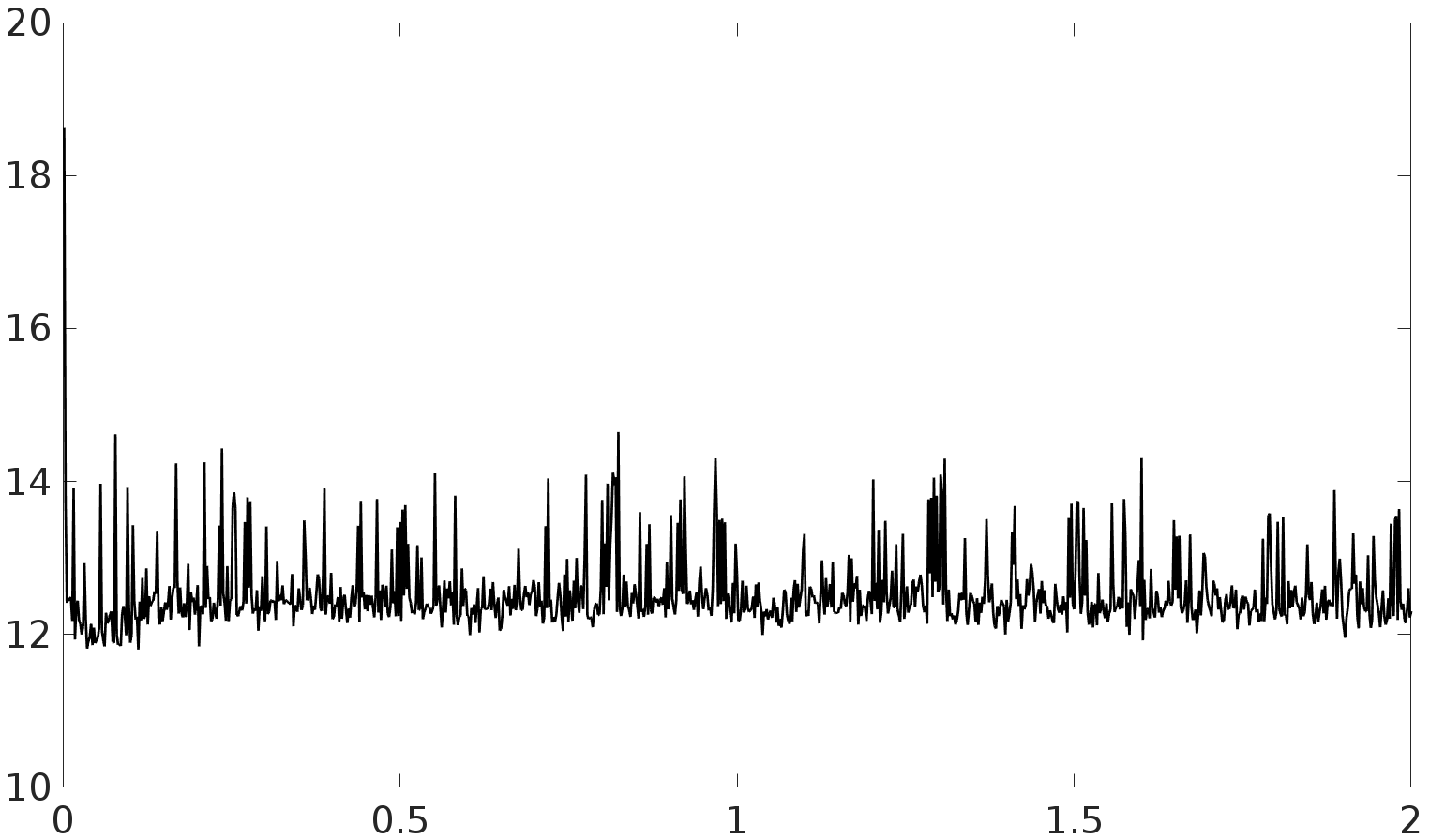}
			\put (47.5,-5) {\scriptsize time}
			\put (14,59) {\scriptsize \textbf{CPU time to assemble the coupling term}}
			\put (-5,25) {\begin{sideways}
					\scriptsize $T_{coup}(s)$
			\end{sideways}}
		\end{overpic}
		\quad
		\begin{overpic}[width=8cm]{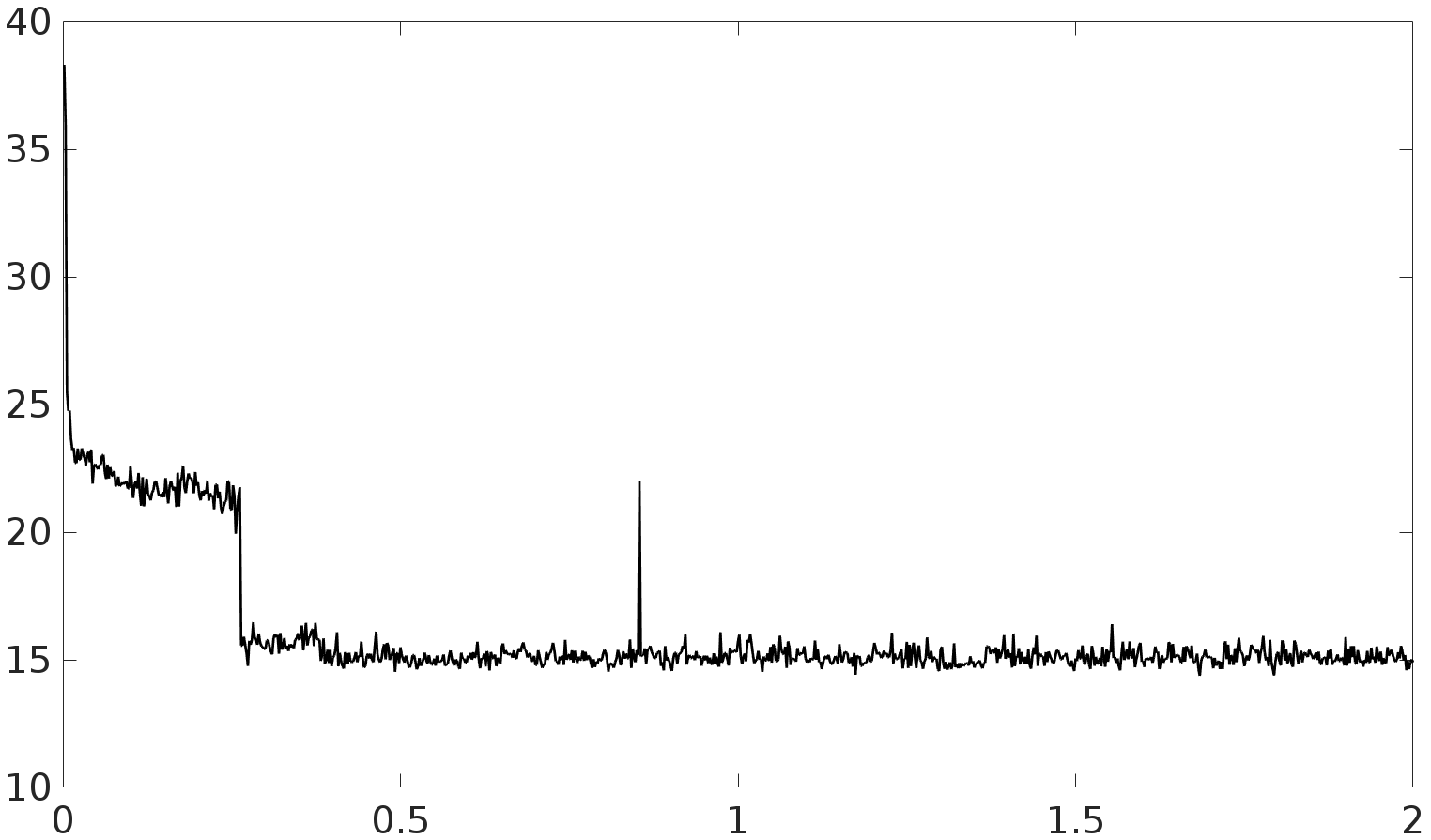}
			\put (48,-5) {\scriptsize time}
			\put (13.5,59) {\scriptsize \textbf{CPU time to solve the nonlinear system}}
			\put (-5,25) {\begin{sideways}
					\scriptsize $T_{sol}(s)$
			\end{sideways}}
		\end{overpic}
		\vspace{1mm}
		\caption{Nonlinear solid model -- Coupling with mesh intersection. 
			Time evolution of nonlinear iterations, average linear iterations per nonlinear iteration, CPU time to assemble
			the coupling term and to solve the nonlinear system. Run on Eos cluster with 32 processors.}
		\label{fig:time_evol_nonlin_inter}
	\end{center}
\end{figure}

\begin{figure}[h!]
	\vspace{5mm}
	\begin{center}
		%		\includegraphics[width=8cm]{figures/fig_nl/fig_nonlin_1.png}
		%		\includegraphics[width=8cm]{figures/fig_nl/fig_nonlin_2.png}\\
		%		\medskip
		%		\includegraphics[width=8cm]{figures/fig_nl/fig_nonlin_3.png}
		%		\includegraphics[width=8cm]{figures/fig_nl/fig_nonlin_4.png}
		\begin{overpic}[width=7.95cm]{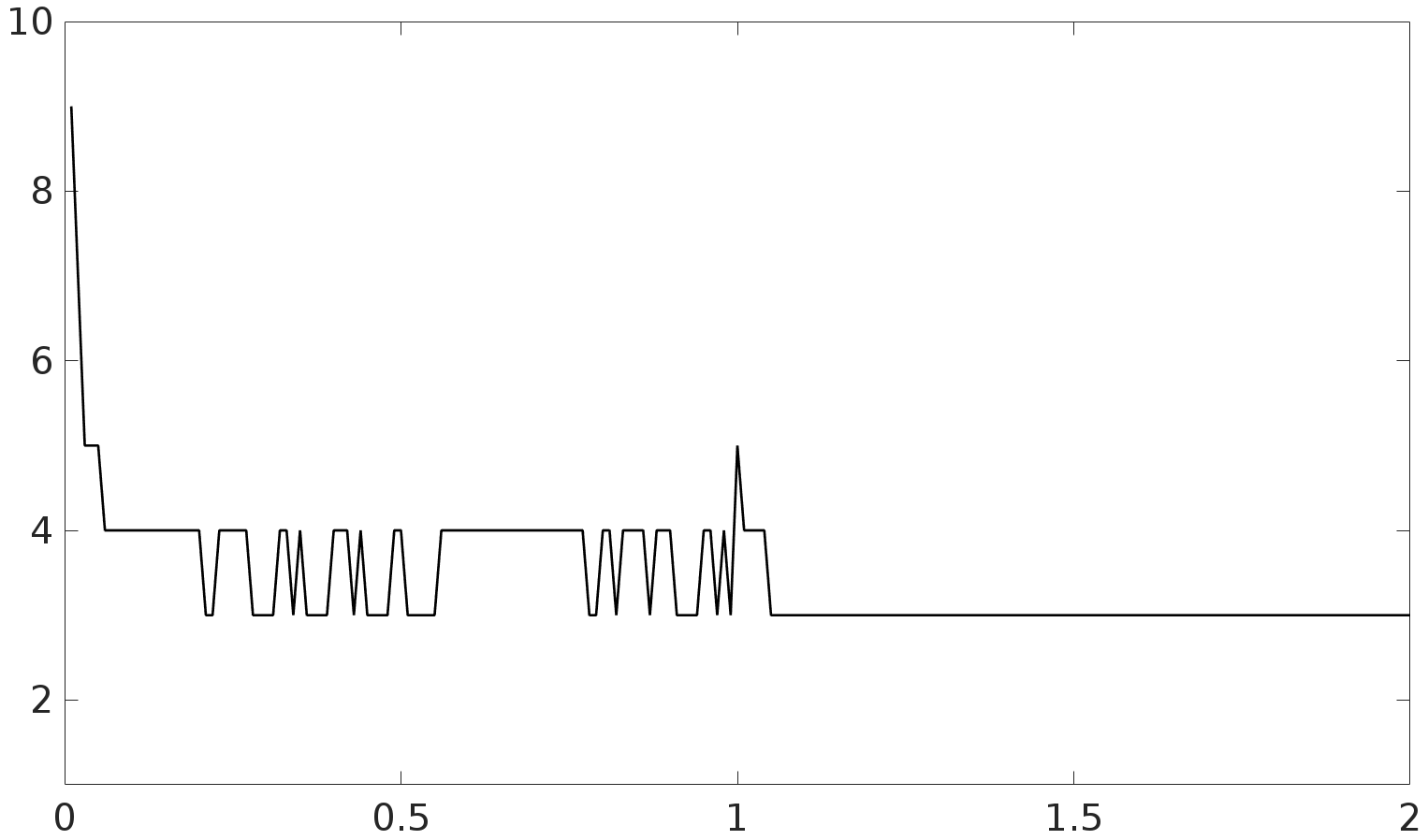}
			\put (47.5,-5) {\scriptsize time}
			\put (32.5,59) {\scriptsize \textbf{Nonlinear iterations}}
			\put (-5,25) {\begin{sideways}
					\scriptsize $nit$
			\end{sideways}}
		\end{overpic}
		\quad
		\begin{overpic}[width=8cm]{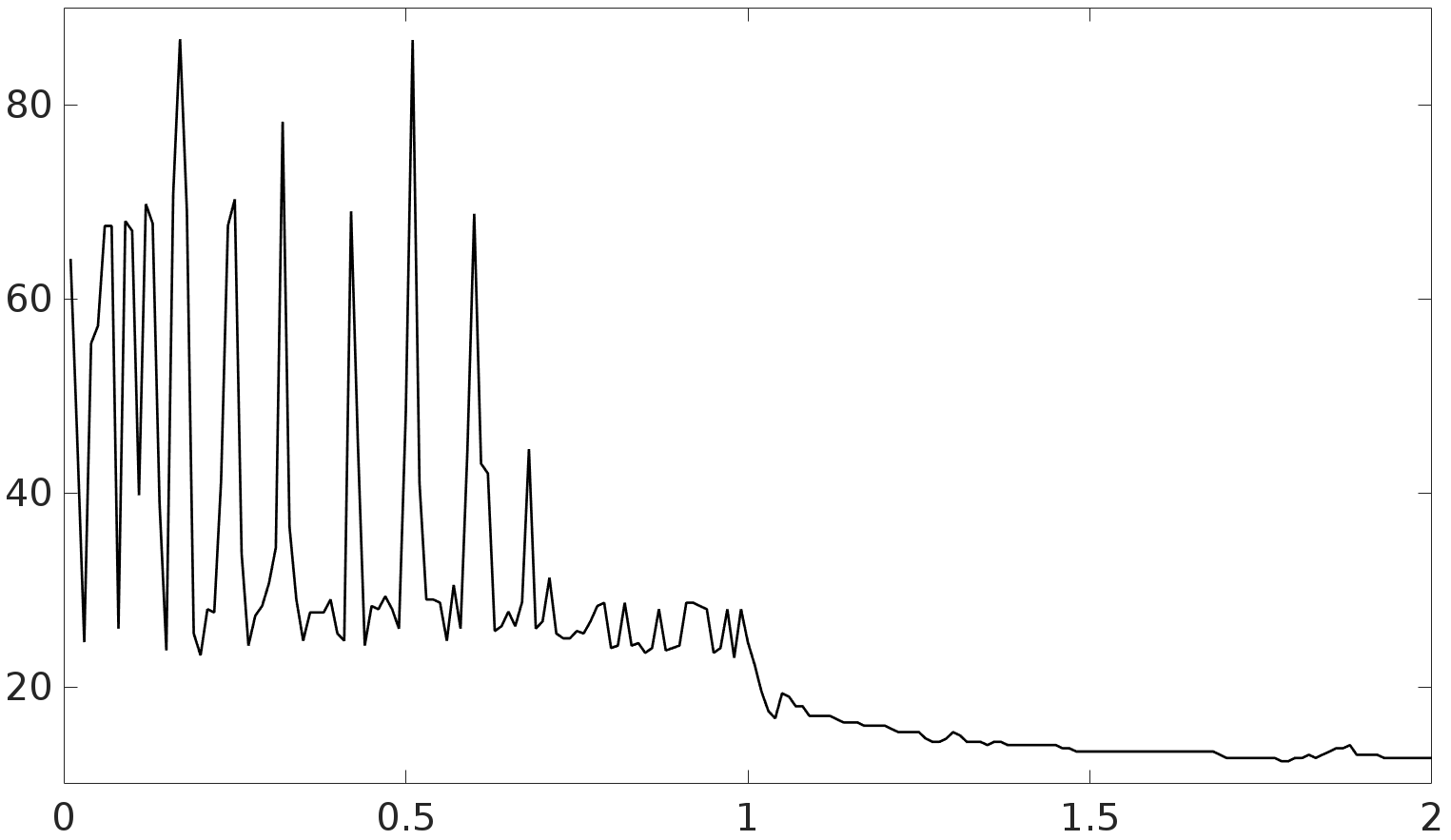}
			\put (48,-5) {\scriptsize time}
			\put (6.2,59) {\scriptsize \textbf{Average linear iterations per nonlinear iteration}}
			\put (-5,25) {\begin{sideways}
					\scriptsize $its$
			\end{sideways}}
		\end{overpic}\\
		\vspace{10mm}
		\begin{overpic}[width=8cm]{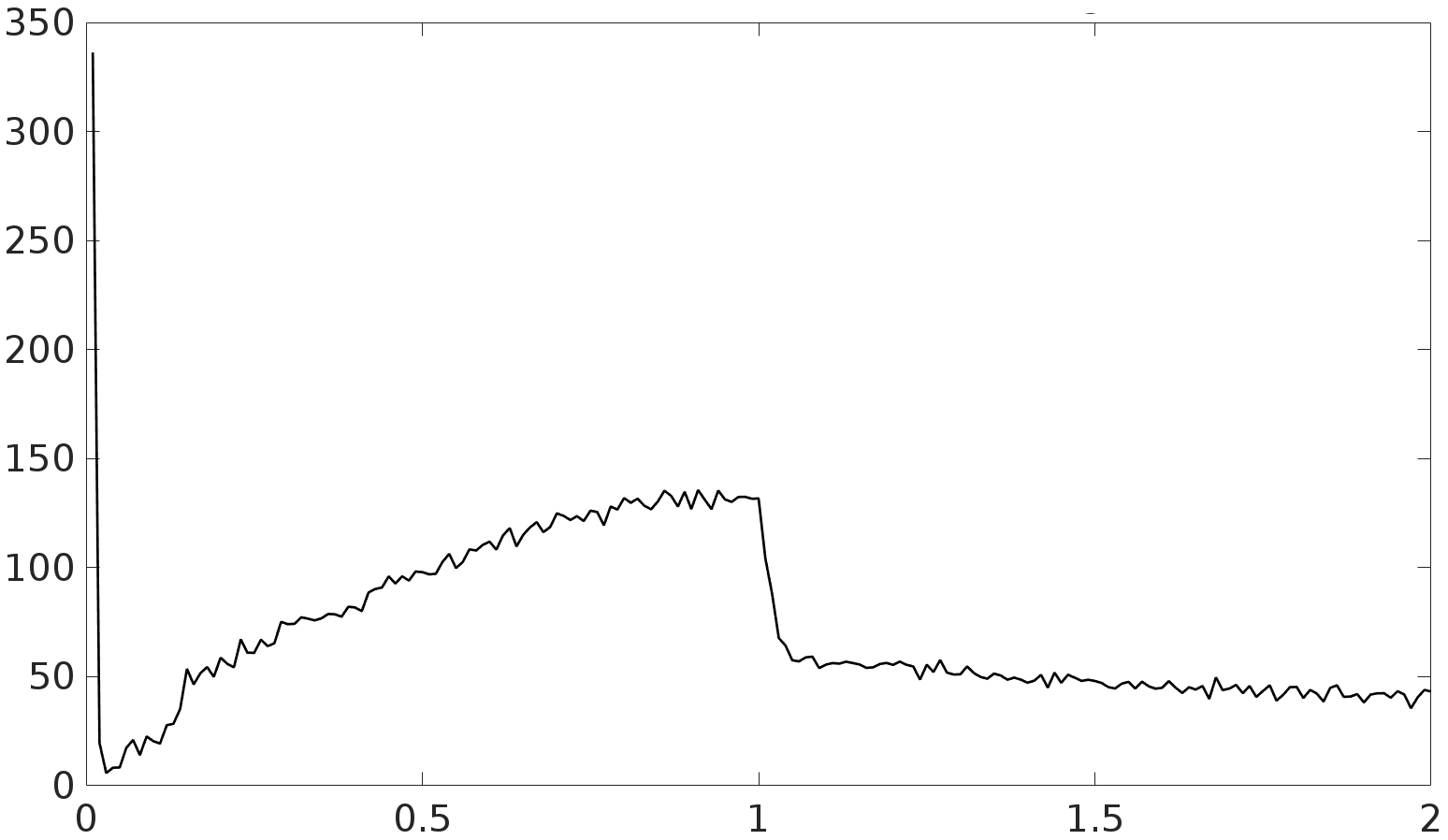}
			\put (48,-5) {\scriptsize time}
			\put (14,59) {\scriptsize \textbf{CPU time to assemble the coupling term}}
			\put (-5,25) {\begin{sideways}
					\scriptsize $T_{coup}(s)$
			\end{sideways}}
		\end{overpic}
		\quad
		\begin{overpic}[width=8cm]{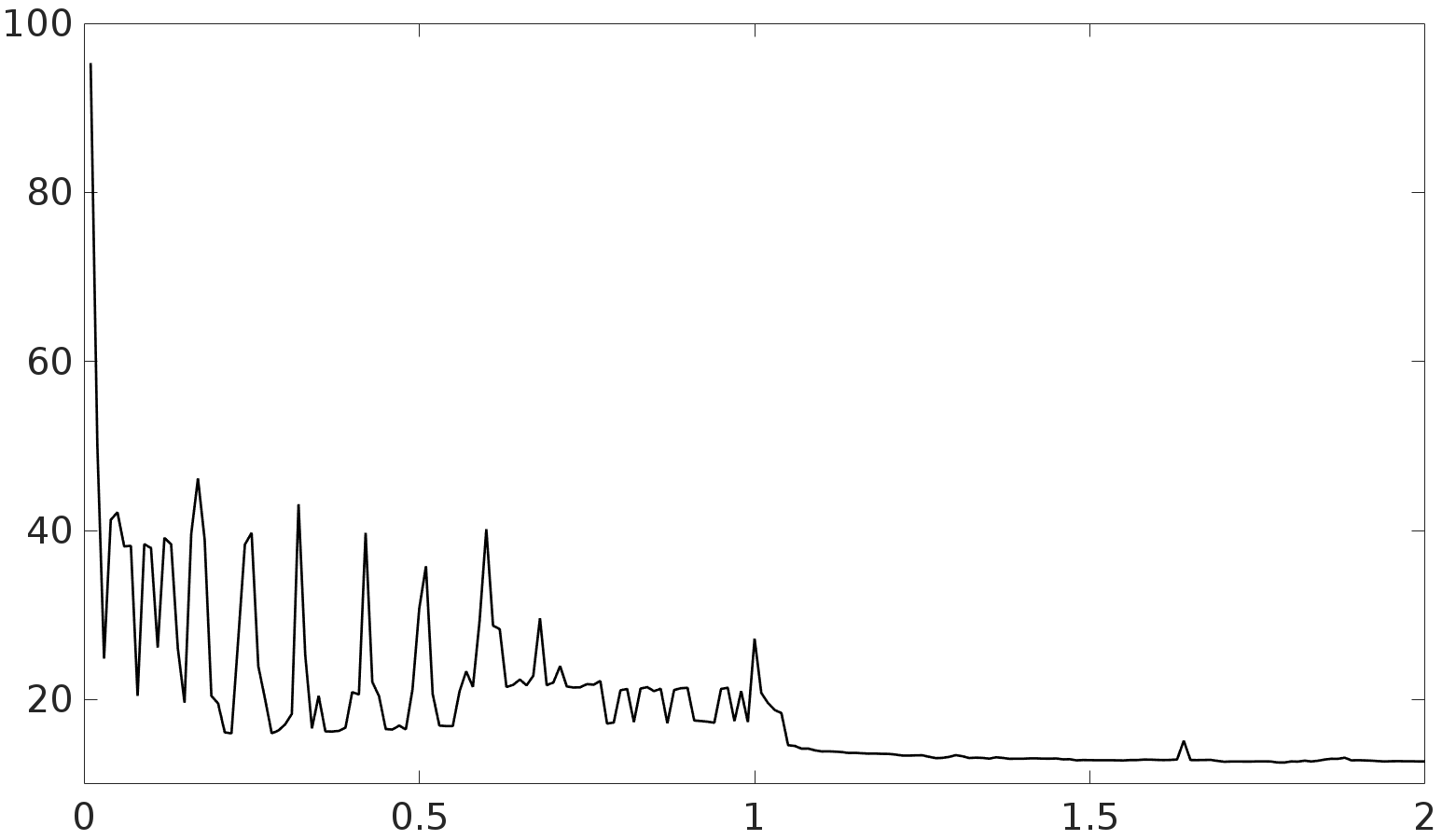}
			\put (48.5,-5) {\scriptsize time}
			\put (13.5,59) {\scriptsize \textbf{CPU time to solve the nonlinear system}}
			\put (-5,25) {\begin{sideways}
					\scriptsize $T_{sol}(s)$
			\end{sideways}}
		\end{overpic}
		\vspace{1mm}
		\caption{Nonlinear solid model -- Coupling without mesh intersection. 
			Time evolution of nonlinear iterations, average linear iterations per nonlinear iteration, CPU time to assemble
			the coupling term and to solve the nonlinear system. Run on Eos cluster with 32 processors.}
		\label{fig:time_evol_nonlin_nointer}
	\end{center}
\end{figure}

In Figures \ref{fig:time_evol_nonlin_inter} and \ref{fig:time_evol_nonlin_nointer}, we report the evolution in time of the number of Newton and GMRES iterations and CPU time (in seconds) to assemble the coupling matrix and solve the nonlinear system on 32 processors of Eos cluster.

\begin{table}
\begin{center}
\begin{tabular}{r|r|r|r|r|r|r|r}
\hline
\multicolumn{8}{c}{{\bf Nonlinear solid model -- Strong scalability test}} \\
\multicolumn{8}{c}{{\textit{Coupling with mesh intersection}}} \\
\hline
\multicolumn{8}{c}{dofs = 515686, T = 2, $\dt$ = 0.002} \\
\hline
procs   & $T_{ass}(s)$  & $T_{coup}(s)$ & \multicolumn{5}{c}{block-tri} \\
        &               &               & nit   & its           & $T_{sol}(s)$  & $T_{tot}$ (s) & $S^p$\\
\hline
4	& 6.79e-1	& 80.89 	& 2	& 6		& 11.17		& 1.05e+5	& - \\
8	& 3.96e-1	& 38.71		& 2	& 6		& 9.02		& 5.79e+4	& 1.81 (2) \\
16	& 2.71e-1	& 21.10		& 2	& 6		& 9.33		& 4.10e+4	& 2.56 (4) \\
32	& 1.71e-1	& 12.53		& 2	& 6		& 7.55		& 2.86e+4	& 3.67 (8) \\
64	& 1.00e-1	& 8.44		& 2	& 6		& 13.29		& 3.68e+4	& 2.85 (16) \\
\hline
\end{tabular}
\vspace*{2mm}
\caption{Test 6, strong scalability in the nonlinear solid model, coupling with mesh intersection. The simulations are run on the Eos cluster. dofs = degrees of freedom; procs = number of processors; $T_{ass}$ = CPU time to assemble the stiffness and mass matrices; $T_{coup}$ = CPU time to assemble the coupling term; nit = Newton iterations; its = GMRES iterations to solve the Jacobian system; $T_{sol}$ = CPU time to solve the Jacobian system; $T_{tot}$ = total simulation CPU time; $S^p$ = parallel speedup computed with respect to the 4 processors run. The quantities $T_{coup}$ and nit are averaged over the time steps, whereas the quantities its and $T_{sol}$ are averaged over the Newton iterations and the time steps. All CPU times are reported in seconds.}
\label{nonlin_scal_int}
\end{center}
\end{table}

\begin{table}
\begin{center}
\begin{tabular}{r|r|r|r|r|r|r|r}
\hline
\multicolumn{8}{c}{{\bf Nonlinear solid model -- Strong scalability test}} \\
\multicolumn{8}{c}{{\textit{Coupling without mesh intersection}}} \\
\hline
\multicolumn{8}{c}{dofs = 515686, T = 2, $\dt$ = 0.01} \\
\hline
procs   & $T_{ass}(s)$  & $T_{coup}(s)$ & \multicolumn{5}{c}{block-tri} \\
        &               &               & nit   & its           & $T_{sol}(s)$  & $T_{tot}$ (s) & $S^p$\\
\hline
4       & 6.31e-1	& 185.50 	& 3	& 33		& 30.49		& 5.93e+4	& - \\
8       & 3.75e-1	& 22.85		& 3	& 33		& 27.62		& 2.48e+4	& 2.39 (2) \\
16      & 1.79e-1	& 4.74		& 3	& 29		& 24.79		& 1.91e+4	& 3.10 (4) \\
32      & 1.19e-1	& 2.51		& 3	& 25		& 18.26		& 1.36e+4	& 4.36 (8) \\
64      & 1.05e-1	& 8.77e-1	& 3	& 22		& 25.24		& 1.81e+4	& 3.28 (16) \\
\hline
\end{tabular}
\vspace*{2mm}
\caption{Test 6, strong scalability in the nonlinear solid model, coupling without mesh intersection. The simulations are run on the Eos cluster.
Same format as in Table \ref{nonlin_scal_int}.}
\label{nonlin_scal_noint}
\end{center}
\end{table}

%\subsection{\lg Time step refinement \gl}
\subsection{Time step refinement}

We now study the behavior of our solver with respect to the refinement of the time marching step. We consider fluid and solid meshes made up of $320\times 320$ and $ 480\times 120$ elements respectively and we fix the number of processors $procs = 64$. The time step $\dt$ decreases from 0.02 to 0.001, while the final time $T$ is still 2.

Table \ref{nonlin_tab_timestep_inters} reports the results in the case of coupling with mesh intersection. The simulations have been performed on the Shaheen cluster. We can notice that for $\dt=0.02$ and $\dt=0.01$ both solvers fail. Conversely, when $\dt$ is refined, they improve their performance: the Newton iterations are always 2, but the GMRES iterations decrease, especially in the case of the block-diag preconditioner, where we have 227 iterations for  $\dt=0.005$ and 16 iterations when $\dt=0.001$. On the other hand, the number of GMRES iterations for the block-tri solver is always small, bounded by 19. Moreover, the volume loss, $T_{ass}$, and $T_{coup}$ remain constant when $\dt$ decreases.

In Table \ref{nonlin_tab_timestep_nointers}, we report the results of the same tests but in the case of coupling without mesh intersection. We can see that the block-tri solver is more robust than the block-diag one but also with respect to itself when applied to the case of coupling with mesh intersection; indeed, here it never fails. In addition, we can notice an improvement of the performance when $\dt$ is refined: even if the number of Newton iterations is small, the GMRES iterations decrease from 193 to 23. The behavior of block-diag is more critical: in the first three cases it fails and for $\dt=0.002$ and $\dt=0.001$ the number of GMRES iterations is greater than 700. Also in this case, the volume loss and  $T_{ass}$ remain substantially constant, while $T_{coup}$ decreases when $\dt$ is refined: \lg probably, this is due to the fact that, when $\Delta t$ becomes small, the pattern of the coupling matrix does not change significantly between subsequent time steps, thus due to the dynamic memory allocation the assembling time $T_{coup}$ reduces. \gl

\begin{table}
	\begin{center}
		\begin{tabular}{r|r|r|r|r|r|r|r|r|r|r|r}
			\hline
			\multicolumn{12}{c}{{\bf Nonlinear solid model -- Time step refinement test}} \\
			\multicolumn{12}{c}{{\textit{Coupling with mesh intersection }}} \\
			\hline
			\multicolumn{12}{c}{dofs = 515686, procs = 64, T = 2} \\
			\hline
			$\dt$    & vol. loss (\%)        & $T_{ass}(s)$  & $T_{coup}(s)$ & \multicolumn{4}{c|}{block-diag}               & \multicolumn{4}{c}{block-tri} \\
			&                       &               &               & nit   & its   & $T_{sol}(s)$  & $T_{tot}(s)$  & nit   & its   & $T_{sol}(s)$  & $T_{tot}(s)$ \\
			\hline
			0.02 & -	 & - & - & -  & - & -  & failed  & - & - & -	& failed \\
			0.01 & -	 & - & - & -  & - & -  & failed  & - & - & -	& failed \\
			0.005 & 1.44 & 1.06e-1 	& 12.04	& 2 & 227 	& 71.35   	& 3.34e+4   & 2	& 19 & 14.97 & 1.08e+4\\
			0.002 & 1.44 & 1.06e-1	& 11.43	& 2 & 59	& 25.49		& 3.70e+4	& 2	& 13 & 12.53 & 2.40e+4\\
			0.001 & 1.44 & 1.06e-1  & 11.24 & 2 & 16	& 13.14		& 4.91e+4	& 2 & 10 & 11.44 & 4.54e+4\\
			\hline
		\end{tabular}
		%\caption{Test 7, refining the time step in the nonlinear solid model. Shaheen cluster.}
		\vspace*{2mm}
		\caption{Test 7, refining the time step $\dt$ in the nonlinear solid model, coupling with mesh intersection. The simulations are run on Shaheen cluster. dofs = degrees of freedom; procs = number of processors; $T_{ass}$ = CPU time to assemble the stiffness and mass matrices; $T_{coup}$ = CPU time to assemble the coupling term; nit = Newton iterations; its = GMRES iterations to solve the Jacobian system; $T_{sol}$ = CPU time to solve the Jacobian system; $T_{tot}$ = total simulation CPU time. The quantities $T_{coup}$ and nit are averaged over the time steps, whereas the quantities its and $T_{sol}$ are averaged over the Newton iterations and the time steps. All CPU times are reported in seconds.}
		\label{nonlin_tab_timestep_inters}
	\end{center}
\end{table}

\begin{table}
	\begin{center}
		\begin{tabular}{r|r|r|r|r|r|r|r|r|r|r|r}
			\hline
			\multicolumn{12}{c}{{\bf Nonlinear solid model -- Time step refinement test}} \\
			\multicolumn{12}{c}{{\textit{Coupling without mesh intersection }}} \\
			\hline
			\multicolumn{12}{c}{dofs = 515686, procs = 64, T = 2} \\
			\hline
			$\dt$    & vol. loss (\%)        & $T_{ass}(s)$  & $T_{coup}(s)$ & \multicolumn{4}{c|}{block-diag}               & \multicolumn{4}{c}{block-tri} \\
			&                       &               &               & nit   & its   & $T_{sol}(s)$  & $T_{tot}(s)$  & nit   & its   & $T_{sol}(s)$  & $T_{tot}(s)$ \\
			\hline
			0.02 	& 1.08e-2	& 1.06e-1 & 51.20&   -	& 	-	&  - 	& failed& 4	& 193 &	65.66 	& 1.17e+4\\
			0.01 	& 7.16e-2	& 1.06e-1 & 35.86&   -	& 	-	&  - 	& failed& 3	& 80  &	37.16	& 1.46e+4\\
			0.005 	& 1.04e-1	& 1.06e-1 & 24.28&   -	& 	-	&  - 	& failed& 3	& 48  &	26.25	& 2.03e+4\\
			0.002 	& 1.24e-1	& 1.05e-1 & 13.49&   3	& 	955	& 283.82&2.95e+5& 3	& 30  &	20.69	& 3.42e+4\\
			0.001 	& 1.31e-1	& 8.66e-2 & 8.62 &   3	& 	710	& 214.62&4.46e+5& 3 & 23  &	17.89	& 5.30e+4\\
			\hline
		\end{tabular}
		\vspace*{2mm}
		\caption{Test 7, refining the time step $\dt$ in the nonlinear solid model. The simulations are run on Shaheen cluster. Same format as in Table~\ref{nonlin_tab_timestep_inters}.}
		\label{nonlin_tab_timestep_nointers}
	\end{center}
\end{table}

\section{\lg Conclusion \gl}
We presented and discussed two possible preconditioners for saddle point problems arising from the discretization of fluid-structure interaction systems modeled by a fictitious domain approach with distributed Lagrange multiplier. \fc\rv As a first step towards the design of an effective parallel solver for a large class of problems, this study has been conducted on 2D examples helpful to understand the challenges posed by the considered formulation and quite representative of the possible applications\vr. For this reason, we assumed that fluid and solid have same density and same viscosity. Moreover, we modeled the fluid dynamics with the Stokes equation, dropping the nonlinear convective term. For the solid we considered only isotropic constitutive laws. Consequently, several extensions are possible from the modeling point of view and they will be addressed in future works.\cf

A crucial aspect of this class of problems is the assembly procedure one can adopt to build the finite element matrix representing the coupling between fluid and structure: one can implement exact composite quadrature rules by computing the intersection of the two meshes or, alternatively, one can compute an approximated matrix by skipping this geometrical computation.

The first preconditioner is a block-diagonal matrix, while the second one is block-triangular: their action is performed by the exact inversion of the two diagonal blocks.
%An analysis of these two preconditioners has been performed in terms of robustness with respect to mesh refinement, strong scalability, and robustness with respect to the time step choice.
\rv At the moment we do not have any theoretical convergence result for the two block preconditioners applied to this highly nonlinear coupling problem. Therefore, an in depth numerical analysis of these two preconditioners has been performed in terms of robustness with respect to mesh refinement, strong scalability, and robustness with respect to the time step choice, considering both linear and nonlinear elastic model for the structure. \vr
Moreover, a comparison is performed in terms of the assembly procedures of the coupling matrix.

The block-tri preconditioner is optimal for both linear and nonlinear solid models, while it appears to be scalable only for the first class of problems. In addition, its behavior in terms of execution times does not change with respect the assembly strategy chosen for the coupling matrix. On the other hand, the block-diag suffers when applied to nonlinear problems or when the coupling matrix is assembled with mesh intersection.

We finally observe that our results have shown that the use of exact composite quadrature for assembling the coupling term does not yield clear improvements with respect to inexact quadrature formulas. A theoretical investigation of this issue is beyond the scope of this work and it will be addressed in a future study.

\section*{Acknowledgments}
The authors are members of INdAM Research group GNCS. The research of L. Gastaldi is partially supported by PRIN/MUR on grant No.20227K44ME and IMATI/CNR.

\bibliographystyle{abbrv}
\bibliography{fsi}

\end{document}